\documentclass[twoside,12pt]{report}% a4paper,12pt,onepage,onecolumn,leqno

\usepackage{cite}
\usepackage{mathrsfs}
\usepackage{helvet}         % selects Helvetica as sans-serif font
\usepackage{courier}        % selects Courier as typewriter font
\usepackage{type1cm}        % activate if the above 3 fonts are
\usepackage{graphicx}        % standard LaTeX graphics tool
\usepackage{multicol}        % used for the two-column index
\usepackage[bottom]{footmisc}% places footnotes at page bottom
\usepackage[total={15.5cm,22.5cm}, bindingoffset=0cm, top=3cm, left=3cm, includefoot, includehead]{geometry}
\usepackage[UglyObsolete]{diagrams}
\usepackage[leqno]{amsmath}
\usepackage{amssymb}
\usepackage{amsthm}
\usepackage[mathscr]{eucal}
\usepackage{amsfonts}
\usepackage{amscd}
\usepackage{url}
\usepackage{color}
\usepackage{graphicx}
\usepackage[%bookmarks=false,%pdfborder=0 0 0,%
colorlinks,linkcolor=blue,citecolor=green,urlcolor=red]{hyperref}
\usepackage{makeidx}         % allows index generation
\DeclareMathAlphabet{\eurm}{U}{eur}{m}{n}
\DeclareMathAlphabet{\eubf}{U}{eur}{b}{n}
\newcommand{\E}[1]{{\eurm{#1}}}

\newtheorem{Statement}{Statement}[section]
\newtheorem{Corollary}{Corollary}[section]
\newtheorem{Definition}{Definition}[section]
\newtheorem{Example}{Example}[section]
\newtheorem{Lemma}{Lemma}[section]
\newtheorem{Note}{Note}[section]
\newtheorem{Proposition}{Proposition}[section]
\newtheorem{Remark}{Remark}[section]
\newtheorem{Theorem}{Theorem}[section]
\newtheorem{Exercise}[Statement]{Exercise}
\newcommand{\bCr}{\begin{Corollary}\em}
\newcommand{\eCr}{\end{Corollary}}
\newcommand{\bDf}{\begin{Definition}\em}
\newcommand{\eDf}{\end{Definition}}
\newcommand{\bEx}{\begin{Example}\em}
\newcommand{\eEx}{\end{Example}}
\newcommand{\bLm}{\begin{Lemma}\em}
\newcommand{\eLm}{\end{Lemma}}
\newcommand{\bNt}{\begin{Note}\em}
\newcommand{\eNt}{\end{Note}}
\newcommand{\bPr}{\begin{Proposition}\em}
\newcommand{\ePr}{\end{Proposition}}
\newcommand{\bRm}{\begin{Remark}\em}
\newcommand{\eRm}{\end{Remark}}
\newcommand{\bTh}{\begin{Theorem}}
\newcommand{\eTh}{\end{Theorem}}
\newcommand{\bDr}{\begin{Exercise}\em}
\newcommand{\eDr}{\end{Exercise}}
\newcommand{\bpf}{\vspace{3pt}\begin{small}{\sc {Proof. }}}
\newcommand{\epf}{\end{small}\vspace{3pt}}
\newcommand{\st}{\;|\;}
\newcommand{\wha}[1]{{\widehat{#1}\,}}
\newcommand{\ul}[1]{{\underline{#1}}}
\newcommand{\dtac}[1]{\dt{\ac{#1}}{}}
\newcommand{\wti}[1]{{\widetilde{#1}}}
\newcommand{\vp}{\; \ba\cdot \;}
\newcommand{\vs}{\; \ba+ \;}
\newcommand{\co}[2]{_{#1}{}^{#2}}
\newcommand{\col}[3]{_{#1}{}^{#2}{}_{#3}}
\newcommand{\END}{{\,\text{\footnotesize\qedsymbol}}}
\DeclareMathOperator{\fib}{{{fib}}}

\DeclareMathOperator{\Map}{{{Map}}}
\DeclareMathOperator{\usec}{{\underline{sec}}}
\DeclareMathOperator{\tub}{{{tub}}}
\DeclareMathOperator{\utub}{{\underline{tub}}}

\DeclareMathOperator{\Fsec}{{{F--sec}}}

\DeclareMathOperator{\pro}{\mathrm{pro}}

\DeclareMathOperator{\aff}{\mathrm{aff}}

\DeclareMathOperator{\cns}{\mathrm{cns}}
\DeclareMathOperator{\lin}{\mathrm{lin}}

\DeclareMathOperator{\map}{{{map}}}

\DeclareMathOperator{\Sec}{{{Sec}}}
\DeclareMathOperator{\pln}{{{pol}}}
\DeclareMathOperator{\cpt}{{{cpt}}}
\DeclareMathOperator{\id}{{{id}}}

\newcommand{\ltag}[1]{\leqno{\quad\text{#1)}}}
\newcommand{\QED}{{\,\,\text{\rm{\footnotesize QED}}}}
\newcommand{\myskip}{\smallskip}
\newcommand{\Wob}{\Big{\{}}
\newcommand{\Wcb}{\Big{\}}}

\newcommand{\lnb}{\linebreak}

\newcommand{\f}[1]{{\boldsymbol{#1}}}

\newcommand{\ba}[1]{{{\bar{#1}}}}

\newcommand{\R}[1]{{{\rm{#1}}}}
\newcommand{\Rarr}{\qquad\Rightarrow\qquad}
\newcommand{\br}[1]{{\breve{#1}}{}}

\newcommand{\brac}[1]{\br{\ac{#1}}{}}

\newcommand{\sst}{\;\;|\;\;}

\newcommand{\dix}[2]%old[3]
{\frac{\der_{#2}(#1^{#2} \, \rtd g)}{\rtd g}}

\newcommand{\diG}[2]%old[3]
{\frac{\der_{#1}(G^{#1#2}_0 \, \rtd g)}{\rtd g}}
\newcommand{\diGl}[2]%old[3]
{\der_{#1}(G^{#1#2}_0 \, \rtd g) / \rtd g}

\newcommand{\ac}[1]{{\acute{#1}}{}}

\newcommand{\ucar}[1]{\underset{#1}{\times}}

\newcommand{\bEq}{\begin{equation}}
\newcommand{\eEq}{\end{equation}}
\newcommand{\beq}{\begin{equation*}}
\newcommand{\eeq}{\end{equation*}}
\newcommand{\bal}{\begin{align*}}
\newcommand{\bml}{\begin{multline*}}
\newcommand{\bgt}{\begin{gather*}}
\newcommand{\eqv}{\,\equiv\,}
\newcommand{\seq}{\,\simeq\,}

\newcommand{\wob}{\big{\{}}
\newcommand{\wcb}{\big{\}}}

\newcommand{\baf}[1]{{\bar{\boldsymbol{#1}}}}

\newcommand{\w}[1]{\big{#1}}
\newcommand{\W}[1]{\Big{#1}}

\newcommand{\der}{\partial}
\newcommand{\ten}{\otimes}
\newcommand{\oset}[1]{\overset{#1}}
\newcommand{\ol}[1]{{\overline{#1}}}
\newcommand{\mto}{\mapsto}
\newcommand{\car}{\times}
\newcommand{\rtd}[1]{\sqrt{|#1|}}
\newcommand{\wed}{\wedge}

\newcommand{\hto}{\hookrightarrow}

\newcommand{\whaE}[1]{{\widehat{\eurm{#1}}}}

\newcommand{\Rn}{\mathbb{R}}
\newcommand{\byd}{\,{\raisebox{.092ex}{\rm :}{\rm =}}\,}

\newcommand{\sub}{\subset}
\newcommand{\nab}{\nabla}
\newcommand{\ssep}[1]{{\qquad\text{\rm{#1}}\qquad}}

\newcommand{\dt}[1]{{\dot{#1}}}
\newcommand{\1}{\mathbf 1}
\newcommand{\com}{{}\circ{}}

\newcommand{\C}[1]{{\mathcal{#1}}}

\newcommand{\Upa}{^{\uparrow}{}}

\newcommand{\emp}{\emph}

\newcommand{\eps}{\epsilon}
\newcommand{\sep}[1]{{\quad\text{\rm{#1}}\quad}}
\newcommand{\acf}[1]{\acute{\f{#1}}{}}
\newcommand{\fz}{\footnotesize}

\newcommand{\baC}[1]{{\bar{{\mathcal{#1}}}}}
\newcommand{\chC}[1]{{\check{\mathcal{#1}}}}

\newcommand{\whaC}[1]{{\widehat{\mathcal{#1}}}}

\newcommand{\acC}[1]{{\acute{{\mathcal{#1}}}}}

\newcommand{\nin}{\notin}
\newcommand{\bat}{\begin{alignat*}}
\newcommand{\bAt}{\begin{alignat}}

\newcommand{\bsm}{\begin{quotation}\small}
\newcommand{\esm}{\end{quotation}\vspace{8pt}}
\newcommand{\alp}{\alpha}
\newcommand{\gam}{\gamma}

\newcommand{\sig}{\sigma}

\newcommand{\lam}{\lambda}
\newcommand{\Lam}{\Lambda}

\newcommand{\ome}{\omega}

\newcommand{\zet}{\zeta}
\newcommand{\bdg}{\beq\begin{diagram}}
\newcommand{\edg}{\end{diagram}\eeq}
\newcommand{\bDg}{\bEq\begin{diagram}}
\newcommand{\eDg}{\end{diagram}\eEq}
\newcommand{\black}{\color{black}}

\newcommand{\introduction}
{
  \chapter*{\introductionname}
  \addcontentsline{toc}{chapter}{\hspace*{1.25em}\introductionname}
}
\newcommand\introductionname{Introduction}
\newcommand{\thelistofsymbols}
               {\chapter*{\symbolsname}
                \addcontentsline{toc}{chapter}{\hspace*{1.25em}%
                  \symbolsname}
               }
\newcommand\symbolsname{List of Symbols}
\newcommand\bibliographyname{Bibliography}
\makeindex
\begin{document}
%--------------------------------------------------------------------%
\pagestyle{empty}

{\ } 
 
\vspace{4cm}
 
\begin{center}
{\bf{\Huge  Smooth and F--smooth systems }}

\vspace{2cm}

{\large{\bf Josef Jany\v{s}ka, Marco Modugno}}
\end{center}

%--------------------------------------------------------------------%

\newpage

{ Josef Jany\v{s}ka

{\fz Department of Mathematics and Statistics, Masaryk University}

{\fz Kotl\'a\v{r}sk\'a 2, 611 37 Brno, Czech Republic}

{\fz email: {\tt janyska@math.muni.cz}}

\vspace{4mm}

 Marco Modugno
 
{\fz Department of Mathematics and Informatics ``U. Dini", 
University of Florence}

{\fz Via S. Marta 3, 50139 Florence, Italy}

{\fz email: {\tt marco.modugno@unifi.it}}
}
%--------------------------------------------------------------------%
%--------------------------------------------------------------------%

\newpage
%\pagenumbering{epmty}
{ }
\vglue1cm

\centerline{\Large\bf PREFACE}
\vglue1cm

{
 We review the geometric theory of \emp{smooth systems of smooth maps}, of \emp{smooth systems of smooth sections} of a smooth double fibred manifold and of \emp{smooth systems of smooth connections} of a smooth fibred manifold.

 Moreover, after reviewing the concept of \emp{F--smooth space} due to A. Fr\"olicher, we discuss the \emp{F--smooth systems of smooth maps}, the \emp{F--smooth systems of fibrewisely smooth sections} of a smooth double fibred manifold, the \emp{F--smooth systems of fibrewisely smooth connections} of a smooth fibred manifold and of \emp{F--smooth connections} of an F--smooth system of fibrewisely smooth sections.
}

\myskip


\myskip


58D15 - Manifolds of mappings

53C05 - Connections general theory
%--------------------------------------------------------------------%
%\newpage
%{\ }

%--------------------------------------------------------------------%
\newpage
\pagestyle{headings}
\pagenumbering{arabic}
\setcounter{page}{4}

{\fz
\tableofcontents}
\markboth{\rm Contents}{\rm Contents}
%--------------------------------------------------------------------%
\newpage

\introduction

\markboth{\rm Introduction}{\rm Introduction}

 The present report is aimed at combining two independent geometric constructions: \emp{smooth systems of smooth sections of a smooth double fibred manifold} and \emp{F--smooth spaces}.
 
\myskip
 
 The notion of \emp{smooth systems of connections} goes back to the paper
\cite{Gar72},
where P.L. Garc\'ia studies the \emp{systems of principal connections} of a principal bundle and the associated ``\emp{universal connection}".
 Indeed, this is a fruitful geometric idea which exploits the properties of the Lie algebra associated with the principal bundle and deserves several applications.
 
 Later, it has been shown that the above geometric construction can be extended to any fibred manifold, at a more basic level, regardless of a possible symmetry group, so detaching the notion of ``system" from principal bundles and their structure group
(see, for instance,
\cite{Cab87,ManMod83d,Mic93}).
 In a few words, given a double fibred manifold
$\f G \to \f F \to \f B \,,$
a ``system of sections" is defined to be a 3--plet
$(\f S, \zet, \eps) \,,$
where
$\zet : \f S \to \f B$
is a fibred manifold and
$\eps : \f S \car_\f B \f F \to \f G$
a fibred morphism over
$\f F \,.$
 Thus, the fibred space
$\f S$
behaves as a space of ``parameters" and the ``evaluation map"
$\eps$
maps sections 
$s : \f B \to \f S$
of the fibred manifold
$\f S \to \f B$
to sections 
$\br s : \f F \to \f G$
of the fibred manifold
$\f G \to \f F \,.$
 Therefore, the choice of such a system of sections turns out to be just a smooth selection of a distinguished family of sections of the fibred manifold
$\f G \to \f F \,.$

 In particular, the above notion can be easily used to define a ``system of connections" of a fibred manifold
$\f F \to \f B \,,$
by setting
$\f G \byd T^*\f B \ten T\f F \,.$
 Actually, in this basic framework, we can recover the ``universal connection" of the system, along with its properties, without explicit reference to the structure group and its Lie algebra.
 In practice, the choice of the ``fibred manifold of parameters"
$\f S$
and of the ``evaluation map"
$\eps$
play the selective role that is played by the equivariance with respect to the structure group in the language of systems of principal connections.

\myskip

 The notion of \emp{F--smooth space} goes back to the paper
\cite{Fro82},
where A. Fr\"olicher introduces a notion of ``smoothness" which is alternative with respect to standard notion.
 Actually, in the present report, we use a definition of F--smoothness with a mild simplification with respect to the original one
(see
\cite{Slo86}).

 In a few words, the notion of \emp{standard smoothness} for a topological manifold
$\f M$
is based on the choice of a family
$\C A$
of local charts
$(x^i) : \f M \to \Rn^n \,,$
which fulfills a standard smoothness compatibility condition.
 On the other hand, the \emp{F--smoothness} for any set 
$\f S$
is based on the choice of a family
$\C C$
of curves
$c_\f I : \f I \to \f S \,,$
which fulfills a condition of smooth re--parametrisation, but no mutual compatibility condition.

 At a first insight, it might appear that the two notions of smoothness above be mutually dual and rather equivalent.
 But, it is not so.
 Indeed, the F--smoothness turns out to be a rather weak condition: it is even more feeble than continuity!
 Nevertheless, one can prove that in the framework of F--smooth spaces it is possible to achieve several advanced geometric constructions
(see, for instance,
\cite{CabKol95,CabKol98,Kol87,KolMod98}).
 It is also worth mentioning further deep investigation in this framework
(see, for instance,
\cite{KriMic97,Mic80}).

\myskip

 In the present report, we are mainly concerned with F--smooth spaces consisting of smooth maps between standard smooth manifolds.
 In our opinion, this is a quite fruitful application of the general notion of F--smoothness.
 
 Accordingly, we generalise the concept of system of sections of a double fibred manifold by replacing the smooth fibred manifold of parameters
$\zet : \f S \to \f B$
with an F--smooth fibred set
$\zet : \f S \to \f B \,.$
 In this way, we can achieve systems of sections, which, in a sense, are infinite dimensional, by skipping the hard methods of infinite dimensional manifolds.
 
\myskip
%\markright{Introduction\hfill \thepage}

 At a first insight, some constructions of the present report might appear very cumbersome.
 But, their core idea is quite simple and intuitive.
 Unfortunately, a detailed account requires odd subtleties; but the reader can grasp the basic simple ideas at a first reading and go throughout details in a second reading, when necessary.
 
\myskip

 As far as applications to mathematical physics are concerned, we have used, in the framework of Covariant Quantum Mechanics, the smooth systems of connections and their universal connection in order to define the ``upper quantum connection", the F--smooth systems of sections for the definition of "sectional quantum bundle" over time and the ``F--smooth connections" in order to regard the Schr\"odinger operator as a connection of the sectional quantum bundle
(see, for instance,
\cite{JanMod02c,JanMod20}).

\myskip

 We use the following symbols:
 
 - given two smooth manifolds
$\f M$
and
$\f N \,,$
we denote the set of \emp{global smooth maps} between the two manifolds by
$\Map(\f M, \f N) \,,$
 
 - given a smooth fibred manifold
$\f F \to \f B \,,$
we denote the set of \emp{global smooth sections}
$s : \f B \to \f F$
by
$\Sec(\f B, \f F) \,.$
 
 - given a smooth fibred manifold
$\f F \to \f B \,,$
we denote the sheaf of \emp{local smooth sections}
$s : \f B \to \f F$
by
$\sec(\f B, \f F) \,.$

%--------------------------------------------------------------------%
\chapter{Smooth manifolds and F--smooth spaces}
\label{Smooth manifolds and F--smooth spaces}
%--------------------------------------------------------------------%
\bsm
 We briefly recall the standard definition of \emp{smooth manifold} and discuss the notion of \emp{F--smooth space}.
 Moreover, we compare these concepts.

\myskip
 
 The reader can be interested to go back to the original literature concerning F--smooth spaces
(see 
\cite{Fro82}
and
\cite{CabKol95,CabKol98,Kol87,KolMod98,KriMic97,Slo86}).
\esm
%--------------------------------------------------------------------%
%\newpage
\markboth{\rm Chapter \thechapter. Smooth manifolds and F--smooth spaces}{\rm \thesection. Smooth manifolds}
\section{Smooth manifolds}
\label{Smooth manifolds}
\markboth{\rm Chapter \thechapter. Smooth manifolds and F--smooth spaces}{\rm \thesection. Smooth manifolds}
%--------------------------------------------------------------------%
\bsm
 We briefly recall the standard definition of \emp{smooth manifold} just in view of a comparison with the forthcoming definition of F--smooth space
(see next
Section \ref{F--smooth spaces}).
\esm

 We observe that the basic original concept of derivative and the consequent concept of smoothness can be achieved in the framework of affine spaces
(see, for instance,
\cite{Sch67}).
 However, for practical reasons, one usually replaces a generic modelling affine space with
$\Rn^n \,.$

 Accordingly, one can introduce the standard notion of smooth manifold as follows.
 
\bDf
\label{Definition: smooth manifolds}
\index{smooth systems!smooth manifolds}
 A \emp{smooth manifold} of dimension
$m$
is defined to be a pair
$(\f M, \C A) \,,$
where
$\f M$
is a topological manifold of dimension
$n$
and 
$\C A$
is a topological atlas whose local charts
$(x^i) : \f M \to \Rn^m \,,$
fulfill ``smooth transition rules" (in the sense of smoothness of affine spaces),  in the intersections of their domains, and a ``maximality condition".

 Moreover, a map
$f : \f M \to \f N$
between two smooth manifolds of dimension
$m$
and 
$n \,,$
respectively, is said to be \emp{smooth} if it yields locally ``local smooth maps" (in the sense of smoothness of affine spaces)
$\Rn^m \to \Rn^n \,.$\END
\eDf

 Thus, the standard smoothness of a manifold
$\f M$
involves a background structure of topological manifold and a  smoothness compatibility condition of local topological charts.
  Hence, the smoothness of manifolds turns out to be a global property achieved via a local condition of compatibility between local topological charts.
%--------------------------------------------------------------------%
%\newpage
\markboth{\rm Chapter \thechapter. Smooth manifolds and F--smooth spaces}{\rm \thesection. F--smooth spacess}
\section{F--smooth spaces}
\label{Section: F--smooth spaces}
\markboth{\rm Chapter \thechapter. Smooth manifolds and F--smooth spaces}{\rm \thesection. F--smooth spaces}
%--------------------------------------------------------------------%
\bsm
 We discuss the notion of \emp{F--smooth space}
$(\f S, \C C)$
and of \emp{F--smooth maps} 
$f : \f S \to \acf S$
between F--smooth spaces.
 
 Moreover, we compare the notions of smoothness and F--smoothness.
\esm
%--------------------------------------------------------------------%
\subsection{F--smooth spaces}
\label{F--smooth spaces}
%--------------------------------------------------------------------%
\bsm
 We introduce the notion of ``\emp{F--smooth space}", with minor changes, with respect to the original definition due to A. Fr\"olicher \cite{Fro82}
(see also 
\cite{CabKol98,KolMicSlo93,KolMod98,KriMic97,Slo86}).

 Roughly speaking, an ``\emp{F--smooth space}" is defined to be a set
$\f S$ 
equipped with a family
$\C C$
of curves
$c : \f I_c \to \f S$
which fulfills a certain feeble requirement of reparametrisation, without compatibility conditions.
 
 This apparently simple difference of the two concepts of smooth manifold and F--smooth space yields great difference in their consequences.
 
 Actually, the notion of ``F--smooth space" is weaker than that of ``smooth manifold".
 Nevertheless, it allows us to achieve several geometric constructions.
 
 On the other hand, any smooth manifold turns out to be an F--smooth space in a natural way
(see the forthcoming
Section \ref{Smooth manifolds as F--smooth spaces}).
 
 We can define the notion of ``\emp{F--smooth subspace}" of an      F--smooth space in a natural way. 
 Indeed, the F--smooth spaces fulfill a remarkable property concerning F--subspaces, which has no analogue for topological or smooth manifolds.

\myskip

 The notion of ``F--smooth space" provides the geometric context suitable for the further discussions on ``F--smooth systems of smooth maps" and ``F--smooth systems of smooth sections"
(see the forthcoming
Section \ref{F--smooth systems of smooth maps}
and
Section \ref{F--smooth systems of smooth sections}).
\esm

\bDf
\label{Definition: F--smooth spaces}
\index{F--smooth systems!F--smooth space}
 An \emp{F--smooth space} 
is defined to be a pair
$(\f S, \C C) \,,$
where 
$\f S$ 
is a non empty set and
$\C C$
is a family of curves, called \emp{basic curves},
\beq
\C C \eqv \{c : \f I_c \to \f S\} \,,
\eeq 
where
$\f I_c \sub \Rn$
is an open subset, such that:

 1) for each
$s \in \f S \,,$
there is at least one curve
$c : \f I_c \to \f S$ 
belonging to
$\C C$
and a
$\lam \in \f I_c \,,$
such that
\beq
c(\lam) = s \,,
\eeq

 2) if 
$c \in \C C$ 
and 
$\gam : \f I_\gam \to \f I_c$ 
is a smooth map defined on an open subset
$\f I_\gam \sub \Rn \,,$ 
then 
\beq
c \circ \gam \in \C C \,.\END
\eeq
\eDf

\myskip

 We have immediate consequences of the above Definition.
 
\myskip

\bPr
\label{Proposition: F--smooth spaces}
\index{F--smooth systems!F--smooth space}
 Let
$(\f S, \C C)$
be an F--smooth space.

 1) All constant curves belong to
$\C C \,.$

 2) If the curve
$c : \f I_c \to \f S$
belongs to
$\C C$,
then its restriction 
$c'$
to any open subset
$\f I_{c'} \sub \f I_c$
belongs to
$\C C$.
\ePr

\bpf
 1) Let
$s \in \f S \,.$
 Then, in virtue of condition 1) in
Definition \ref{Definition: F--smooth spaces}, 
there exists a curve 
$c : \f I_c \to \f S$
which belongs to
$\C C$
and such that, for a certain
$\lam \in \f I_c$,
we have
$c(\lam) = s \,.$

 Moreover, let us consider any open subset
$\f I_\gam \sub \Rn \,,$
the above
$\lam \in \f I_c$
and the constant map
$\gam : \f I_\gam \to \f I_c : \mu \mto \lam \,.$
 Then, according to condition 2) in
Definition \ref{Definition: F--smooth spaces},
the constant curve
$c \com \gam : \f I_\gam \to \f S$
belongs to
$\C C \,.$

 2) The curve
$c' : \f I_{c'} \to \f S$
can be regarded as the composition
$c' = c \com \gam$,
where
$\gam : \f I_{c'} \hto \f I_c$
is the smooth inclusion map.\QED
\epf

\myskip

 The families of basic curves of F--smooth spaces fulfill the following remarkable property, which has no analogue concerning the charts of smooth manifolds.
 
\myskip

\bPr
\label{Proposition: no compatibility conditions required by basic curves}
\index{F--smooth systems!F--smooth space}
 Let us consider a set
$\f S$
and two families
$\C C$
and
$\acC C$
of
$\f S \,,$
which fulfill the requirements of
Definition \ref{Definition: F--smooth spaces}.

 Then, the union
$\C C \cup \acC C$
of the two families of curves fulfills the requirements of
Definition \ref{Definition: F--smooth spaces}.\END
\ePr

\myskip

 Next, we compare different F--smooth structures of a set
$\f S \,.$

\myskip

\bNt
\label{Note: F--smooth space}
\index{F--smooth systems!F--smooth space}
 Any non empty set
$\f S$
admits at least one F--smooth structure.
 Even more, if
$\f S$
has infinitely many elements, then it admits infinitely many F--smooth structures.\END
\eNt

\bDf
\label{Definition: finer F--smooth structure}
\index{F--smooth systems!F--smooth space}
 Given two F--smooth structures
$\C C$
and
$\acC C$
of a set
$\f S \,,$
such that
$\C C \sub \acC C \,,$
then we say that
$\acC C$
is \emp{finer} than
$\C C \,.$\END
\eDf

\myskip

 It is worth comparing the notions of F--smooth space and smooth manifold.
 
\myskip

\bRm
\label{Remark: comparison of F--smooth spaces and smooth manifolds}
\index{F--smooth systems!F--smooth space} 
 In a sense,
  
a) the assignment of the set of curves
$\C C$
for an F--smooth space plays a role analogous to the assignment of a
maximal smooth atlas
$\C A$
for a smooth manifold,

b) the assignment of a set
$\ol{\C C}$
as in the forthcoming 
Example \ref{Example: F--smooth set}
plays a role analogous to the assignment of a smooth atlas
$\baC A$
for a smooth manifold,

 c) the conditions 1) and 2) in 
Definition \ref{Definition: F--smooth spaces} 
for an F--smooth space play a role analogous to the conditions of smooth transitions for the charts of a smooth atlas.\END
\eRm

\myskip

 Eventually, we consider a few trivial or exotic examples of F--smooth spaces, just to account for the range of the notion of F--smooth space.

\myskip

\bEx
\label{Example: F--smooth set}
\index{F--smooth systems!F--smooth space}
 Let us consider any set
$\f S \,.$

 Moreover, let us choose any set
$\ol{\C C}$
consisting of curves
$c : \f I_c \to \f S$
fulfilling condition 1) in
Definition \ref{Definition: F--smooth spaces}.

 Furthermore, let us define the set
$\C C$
consisting of all curves of the type
\beq
c \com \gam : \f I_\gam \to \f S \,,
\ssep{where}
c \in \ol{\C C} \,,
\eeq
and where
$\gam : \f I_\gam \to \f I_c$
is any smooth map defined in a open subset
$\f I_\gam \sub \Rn \,.$

 Then, the pair
$(\f S, \C C)$
turns out to be an F--smooth space.\END
\eEx

\bEx
\label{Example: exotic F--smooth space with constant curves}
\index{F--smooth systems!F--smooth space}
 Let us consider any set
$\f S$
and the set
$\C C$
consisting just of all constant curves
$c : \f I_c \to \f S \,.$
Then, the pair
$(\f S, \C C)$
turns out to be an F--smooth space.

 In this case, all other possible F--smooth structures are finer than the above structure.\END
\eEx

\bEx
\label{Example: exotic F--smooth space with strange curves}
\index{F--smooth systems!F--smooth space}
 Let us consider the set
$\f S \byd \Rn^2 \,,$
along with the natural smooth chart
$(x,y) \,,$
and the set
$\C C$
consisting of all smooth curves 
$c$
whose coordinate expression is of the type 
\beq
c^x(\lam) = a \,, 
\qquad
c^y(\lam) = \gam(\lam) \,,
\eeq
where
$a \in \Rn$
and
$\gam : \f I_\gam \to \Rn$
is any smooth curve.

 Then, the pair
$(\f S, \C C)$
turns out to be an F--smooth space.\END
\eEx

\bEx
\label{Example: exotic F--smooth space with other strange curves}
\index{F--smooth systems!F--smooth space}
 Let us consider the set
$\f S \byd \Rn$
and the set 
$\C C$
consisting of all smooth curves
$c : \f I_c \to \f S \,,$
such that
\beq
0 \leq |c(\lam_2) - c(\lam_1)| \leq 1 \,,
\ssep{for each}
\lam_1, \lam_2 \in \f I_c \,.
\eeq

 Then, the pair
$(\f S, \C C)$
turns out to be an F--smooth space.\END
\eEx
%--------------------------------------------------------------------%
\subsection{F--smooth subspaces}
\label{F--smooth subspaces}
%--------------------------------------------------------------------%
\bsm
 Further, we discuss the F--smooth subspaces of an F--smooth space and emphasise an interesting property, which has no analogue for topological or smooth manifolds.
\esm

\bDf
\label{Definition: F--smooth subspace}
\index{F--smooth systems!F-smooth subspace}
 Let 
$(\f S, \C C)$
be an F-smooth space.
 Then, an \emp{F--smooth subspace} is defined to be a pair
\beq
(\f S', \C C') \subseteq (\f S, \C C) \,,
\eeq
where
$\f S' \subseteq \f S$
is a non empty subset and
$\C C' \subseteq \C C$
is the subset consisting of all curves belonging to
$\C C$
and with values in
$\f S' \,.$

 For short, we say also that
$\f S' \subseteq \f S$
is an \emp{F--smooth subspace}.\END
\eDf

\bPr
\label{Proposition: F--smooth subspace}
\index{F--smooth systems!F-smooth subspace}
 Let us consider an F--smooth space
$(\f S, \C C)$
and \emp{any} non empty subset
$\f S' \subseteq \f S \,.$
 Then,
$\f S'$
turns out to be an F--smooth space in a natural way.

 In fact, let 
$\C C' \subseteq \C C$
be the subset consisting of all curves
$c \in \C C \,,$
whose image is contained in
$\f S' \,.$
 Then, the pair
$(\f S', \C C')$
turns out to be an F--smooth space.
\ePr

\bpf
 Let us check the two conditions of
Definiton \ref{Definition: F--smooth spaces}.

 1) The 1st condition is fulfilled because all constant curves of
$\f S'$
belong to
$\C C' \,.$

 2) The 2nd condition is also fulfilled.
 In fact, let
$c \in \C C'$
and let
$\gam : \f I_\gam \to \f I_c$
be a smooth map.
 Then, the composition
$c \com \gam$
belongs to
$\C C$
and has values in
$\f S' \,.$\QED
\epf

\bRm
 Let us consider an F--smooth space
$(\f S, \C C)$
and any non empty subset
$\f S' \subseteq \f S \,.$

 Then, besides the above structure of F--smooth subspace 
(see
Definition \ref{Definition: F--smooth subspace}),
the subset
$\f S'$
might be equipped with other F-smooth structures given by a further
subset
$\C C'' \subseteq \C C' \subseteq \C C \,.$

 Thus, the above F--smooth structure 
$\C C'$
is the finest F--smooth structure among those which are comparable with the F--smooth structure
$\C C$
of the environmental space
$\f S \,.$\END
\eRm

\myskip

 We exhibit an exotic example of F--smooth subspace.
 
\myskip

\bEx
\label{Example: F--smooth structure given by irrational numbers}
\index{F--smooth systems!F--smooth space}
 Let us consider the F--smooth space
$(\f S, \C C)$,
where
$\f S \byd \Rn$
and
$\C C$
is the set of all smooth curves of
$\Rn \,.$

 Moreover, consider the subset
$\f S' \sub \f S$
consisting of all irrational numbers.

 Then, the curves belonging to
$\C C$
and with values in
$\f S'$
are just the constant curves with irrational image.
 Let
$\C C'$
be the set of these curves.

 Indeed, the pair
$(\f S', \C C')$
turns out to be an F--smooth subspace of
$(\f S, \C C) \,.$\END
\eEx
%--------------------------------------------------------------------%
\subsection{F--smooth maps}
\label{F--smooth maps}
%--------------------------------------------------------------------%
\bsm
 We introduce the notion of ``\emp{F--smooth map}" between F--smooth
spaces in a natural way.
 According to this definition, F--smooth spaces along with global F--smooth maps constitute a category.

 Indeed, the restriction property of F--smooth maps is analogous to a property holding for smooth maps; on the other hand a gluing property of F--smooth maps does not hold for F--smooth spaces.
\esm

\bDf
\label{Definition: F--smooth maps}
\index{F--smooth systems!F-smooth map}
 Let us consider two F--smooth spaces
$(\f S, \C C)$ 
and 
$(\f S', \C C') \,.$

 Then, a (local) map 
\beq
f : \f S \to \f S'
\eeq
(defined on a subset
$\f U \sub \f S$)
is said to be \emp{F--smooth} if, for each 
$c \in \C C$
(with values in
$\f U$)
we have 
\beq
c' \eqv f \com c \in \C C' \,.\END
\eeq
\eDf

\bPr
\label{Proposition: composition of global F--smooth maps}
\index{F--smooth systems!F-smooth map}
 If
$\f S$
is an F-smooth space, then the global map
\beq
\id_\f S : \f S \to \f S
\eeq
is F--smooth.

 Mooreover, if
$\f S,\,$ 
$\f S', \,$
$\f S''$
are F--smooth spaces and
\beq
f : \f S \to \f S'
\ssep{and}
f' : \f S' \to \f S''
\eeq
are global F--smooth maps, then
$f' \com f : \f S \to \f S''$
turns out to be a global F--smooth map.\END
\ePr

\bRm
\label{Remark: two F--smooth structures}
\index{F--smooth systems!F-smooth map}
 Let us consider a set
$\f S$
equipped with two F--smooth structures
$\C C$
and
$\C C' \,,$
and suppose that
$\C C' \sub \C C \,.$

 Then,
\beq
\id_\f S : (\f S, \C C') \to (\f S, \C C)
\eeq
turns out to be F--smooth, but
\beq
\id_\f S : (\f S, \C C) \to (\f S, \C C')
\eeq
is not F--smooth.\END
\eRm

\myskip

 We have a restriction property of F--smooth maps analogous to a
property of smooth maps.

\bPr
\label{Proposition: restriction property of F--smooth maps}
\index{F--smooth systems!F-smooth map}
 Let us consider two F--smooth spaces
$(\f S, \C C)$ 
and 
$(\f S', \C C') \,.$

 If
$f : \f S \to \f S'$
is an F--smooth map, defined on a subset
$\f U \sub \f S$,
then the restriction
$\ba f$
of
$f$
to a further subset
$\baf U \sub \f U$
turns out to be F--smooth.
\ePr

\bpf
 If
$c' \eqv f \com c \in \C C' \,,$
for each
$c \in \C C$
with values in
$\f U \,,$
then
$c' \eqv f \com \ba c \in \C C' \,,$
for each
$\ba c \in \C C$
with values in
$\baf U \,.$\QED
\epf

\myskip

 The F--smooth maps do not have a gluing property
analogous to the gluing property of smooth maps.

 The reason for this difference of behaviour between smooth and F--smooth maps is due to the fact that F--smoothness of
maps is a \emp{global property}, while smoothness of maps is a \emp{local property}.

\myskip

\bRm
\label{Remark: gluing property of F--smooth maps}
\index{F--smooth systems!F-smooth map}
 Let us consider two F--smooth spaces
$(\f S, \C C)$ 
and 
$(\f S', \C C')$ 
and two maps
\beq
f_1 : \f S \to \f S'
\ssep{and}
f_2 : \f S \to \f S' \,,
\eeq
defined respectively on the subsets
$\f U_1 \sub \f S$
and
$\f U_2 \sub \f S \,.$

Clearly, if
$f_1 = f_2$
on
$\f U_1 \cap \f U_2$,
then 
$f_1$
and
$f_2$ 
yield a ``glued map"
$f : \f S \to \f S'$,
defined on
$\f U_1 \cup \f U_2 \,.$

 However, the hypothesis that
$f_1$
and
$f_2$
be F--smooth does not imply that
$f$
be F--smooth.

 In fact, let us provide an exotic countre--example.

 Let us consider the F--smooth spaces
$(\f S, \C C)$ 
and
$(\f S', \C C') \,,$
where
$\f S = \f S' = \Rn \,.$

 Suppose that
$\C C$
consists of all smooth curves of
$\f S$
and that
$\C C'$
consists of all smooth curves of
$\f S'$
considered in
Example \ref{Example: exotic F--smooth space with other strange curves}.

 Next consider the maps
\beq
f_1 = \id : \f U_1 = (0, 1) \to \f S'
\ssep{and}
f_2 = \id : \f U_2 = (1/2, 3/2) \to \f S' \,.
\eeq

 Clearly, these maps are F--smooth and coincide in the intersection of
their domains.

 However, the glued map
$f : \f S \to \f S' \,,$
defined in 
$(0, 3/2)$
is not F--smooth.\END
\eRm

\myskip

 Next, we discuss the behaviour of F--smooth maps with respect different smooth structures.
 
\myskip

\bPr
\label{Proposition: two different F--smooth structures}
 Let 
$\f S$
be a set equipped with two different F--smooth structures
$(\f S, \C C)$
and
$(\f S, \ol{\C C}) \,,$
with
$\C C \subseteq \ol{\C C}$
and let
$(\f S', \C C')$
be another F--smooth space.

 Let
\beq
\Map_{\C C \C C'}(\f S, \, \f S') \,,
\qquad
\Map_{\ol{\C C} \C C'}(\f S, \, \f S')
\sep{and}
\Map_{\C C' \C C}(\f S', \, \f S) \,,
\qquad
\Map_{\C C' \ol{\C C}}(\f S', \, \f S)
\eeq
be, respectively, the sets of maps
$\f S \to \f S'$
and
$\f S' \to \f S \,,$
which are F--smooth with respect to the two different F--smooth
structures of
$\f S \,.$

 Then, we have
\beq
\Map_{\ol{\C C} \C C'}(\f S, \, \f S') \subseteq 
\Map_{\C C \C C'}(\f S, \, \f S')
\sep{and}
\Map_{\C C' \C C}(\f S', \, \f S) \subseteq 
\Map_{\C C' \ol{\C C}}(\f S', \, \f S) \,.\END
\eeq
\ePr

\bEx
\label{Example: F--smooth space with constant maps}
 Let
$(\f S, \C C)$
be an F--smooth space and let
$\C C$
consist just of constant curves of
$\f S$
(see
Example \ref{Example: exotic F--smooth space with constant curves}).

 Moreover, let
$(\f S', \C C')$
be another F--smooth space and let 
$\C C'$
contain also non constant curves of
$\f S'$
passing through each point of
$\f S' \,.$

 Then, a map
$f : \f S' \to \f S$
is F--smooth if and only if it is constant.

 Moreover, all maps
$f : \f S \to \f S'$
are F--smooth.\END
\eEx

\myskip

 Next, we consider an exotic example of F--smooth maps, just to account for the range of the notion of F--smooth map.

\myskip

\bEx
\label{Example: two different smooth and F--smooth structures}
 Let us consider the set
$\f S \byd \Rn$
and equip it with two different smooth structures induced respectively by the two charts
(see also
Remark \ref{Remark: 1st tangent prolongation and tangent prolongation  of a system of maps need not to be injective})
\beq
x : \f S \to \Rn : s \mto s
\ssep{and}
\ac x : \f S \to \Rn : s \mto s^3 \,.
\eeq

 Clearly, the 1st smooth structure is the natural one, while the 2nd one is an exotic smooth structure.
 Indeed, these smooth structures are different because the transition rule
\beq
x = (\ac x)^{1/3}
\eeq
is not differentiable in
$s = 0 \in \f S \,.$
 
 Next, let us define two F--smooth structures
$(\f S, \C C)$
and
$(\f S, \C C')$
on
$\f S \,,$
where, respectively,
$\C C$
and
$\acC C$
consist of the curves
$c : \f I_c \to \f S$
and
$\ac c : \f I_{\ac c} \to \f S$
whose coordinate expressions
\beq
x \com c : \f I_c \to \Rn
\ssep{and}
\ac x \com \ac c : \f I_c \to \Rn \,,
\eeq
are smooth in the standard sense, with respect to the above smooth structures of
$\f S \,,$
respectively.

 Indeed, the above F--smooth structures are different.
 To prove this fact, the curves with coordinate expressions
\beq
x \com c : \f I_c \to \Rn : \lam \mto \lam
\ssep{and}
\ac x \com \ac c : \f I_{\ac c} \to \Rn : \lam \mto \lam
\eeq
belong, respectively, to
$\C C$
and
$\acC C \,.$
 
 Actually, we have
$\ac c \nin \C C \,,$
because the coordinate expression
\beq
x \com \ac c : \f I_{\ac c} \to \Rn : \lam \mto \lam^{1/3}
\eeq
is not differentiable in
$s = 0 \in \f S \,.$

 On the other hand, we have
$\C C \sub \acC C \,,$
because, for each
$c \in \C C \,,$
the coordinate expression
\beq
\ac x \com c : \f I_c \to \Rn : \lam \mto c(\lam)^3
\eeq
is smooth in the standard sense.\END
\eEx
%--------------------------------------------------------------------%
\subsection{Cartesian product of F--smooth spaces}
\label{Cartesian product of F--smooth spaces}
%--------------------------------------------------------------------%
\bsm
 Eventually, we analyse the F--smooth cartesian product of F--smooth spaces.
\esm

\bPr
\label{Proposition: F--smooth cartesian product}
 Let
$(\f S', \C C')$ 
and 
$(\f S'', \C C'')$ 
be F-smooth spaces and set
\beq
\C C \byd 
\Wob
(c',c'') \in \C C' \car \C C'' 
\;|\;
\f I_{c'} = \f I_{c''}
\Wcb \,.
\eeq

 Then, the pair
$(\f S' \car \f S'', \, \C C)$ 
turns out to be an F--smooth space.

 Moreover, the natural projections
\beq
\pro' : \f S' \car \f S'' \to \f S'
\ssep{and}
\pro'' : \f S' \car \f S'' \to \f S''
\eeq
turn out to be F-smooth.
\ePr

\bpf
 1) The pair
$(\f S' \car \f S'', \, \C C)$
is an F--smooth space.

 In fact, we can easily see that the set
$\C C$
fulfills the properties of
Definition \ref{Definition: F--smooth spaces}.

\myskip

 2) The maps
\beq
\pro' : \f S' \car \f S'' \to \f S'
\ssep{and}
\pro'' : \f S' \car \f S'' \to \f S''
\eeq
are F--smooth.

 In fact, for each F--smooth curve
$c \byd (c', c'') : \f I_c \to \f S' \car \f S'' \,,$
the curves
\beq
c' = \pro' \com c
\ssep{and}
c'' = \pro'' \com c
\eeq
are F--smooth by hypothesis.\QED
\epf
%--------------------------------------------------------------------%
%\newpage
\markboth{\rm Chapter \thechapter. Smooth manifolds and F--smooth spaces}{\rm \thesection. Smooth manifolds as F--smooth spaces}
\section{Smooth manifolds as F--smooth spaces}\black
\label{Smooth manifolds as F--smooth spaces}
\markboth{\rm Chapter \thechapter. Smooth manifolds and F--smooth spaces}{\rm \thesection. Smooth manifolds as F--smooth spaces}
%--------------------------------------------------------------------%
\bsm
 Each smooth manifold turns out to be an F--smooth space in a natural way.

 On the other hand, each set
$\f S$
might be equipped with many smooth structures; hence, these smooth manifolds turn out to be F--smooth spaces in many ways.

 Indeed, we can achieve a ``universal" F--smooth structure of a set
$\f S$
by taking the union of all families of basic curves of all F--smooth structures.
 
 Moreover, an F-smooth space
$\f S$
might be a smooth manifold in many ways.

\myskip

 Here, we do not address the above relation between smooth structures and F--smooth structures in detail and full generality.
 But, we just clarify this question by an example.
\esm

\bTh
\label{Theorem: relation between smoothness and F--smoothness}
 Each smooth manifold
$\f M \,,$
along with the set
$\C C$
consisting of all smooth curves
$c : \f I_c \to \f M \,,$
turns out to be an F--smooth space.

 Moreover, a map between smooth manifolds turns out to be F--smooth (in the sense of the above natural F--smooth structures) if and
only if it is smooth.
\eTh

\bpf
 The non trivial proof of this result is due to 
A. Fr\"olicher \cite{Fro82}
and is based on the 
Boman's theorem \cite{Bom67}.QED
\epf

\bDf
\label{Definition: natural F--smooth structure of smooth manifolds}
\index{F--smooth systems!natural F--smooth structure}
 According to the above
Theorem \ref{Theorem: relation between smoothness and F--smoothness},
for each smooth manifold
$\f M \,,$
we define its \emp{natural F--smooth structure} to be provided by the set 
$\C C$
consisting of all smooth curves
$c : \f I_c \to \f M \,.$\END
\eDf

 Indeed, from now on, for each smooth manifold 
$\f M \,,$
we shall refer to its natural F--smooth structure.

\bCr
\label{Corollary: F--smooth maps}
 Let
$\f M$
be a smooth manifold.
 Then, a function
$f : \f M \to \Rn$
is F--smooth if and only if the map
$f \com c : \f I_c \to \Rn$
is smooth for each smooth curve
$c : \f I_c \to \f M \,.$\END
\eCr

\bCr
\label{Corollary: basic curves and F--smooth curves}
 Let
$(\f S, \C C)$
be an F--smooth space and
$c : \f I_c \to \f S$
any curve.

 Then,
$c$ 
is F--smooth if and only if
$c \in \C C \,.$
\eCr

\bpf
 Let us consider
$\f I_c$
as a smooth manifold equipped with its natural F--smooth structure.

 1) If
$c \in \C C \,,$
then, in virtue of 
Definition \ref{Definition: F--smooth spaces}, 
for each smooth curve
$\gam : \f I_\gam \to \f I_c \,,$
we have
$c \com \gam \in \C C \,.$

 Hence, in virtue of
Definition \ref{Definition: F--smooth maps},
$c$
is F--smooth.

\myskip

 2) If
$c$
is F--smooth, then, in virtue of
Definition \ref{Definition: F--smooth maps},
for each smooth curve
$\gam : \f I_\gam \to \f I_c \,,$
we have
$c \com \gam \in \C C \,.$

 Hence, in particular, we have
$c = c \com \id_{\f I_c} \in \C C \,.$\QED
\epf

\myskip

 Eventually, we provide an example suitable to compare different smooth and F-smooth structures.
 
\myskip

\bEx
\label{Example: several smooth structures underlying an F--smooth space}
 Let us consider the set
$\f S \byd \Rn \,.$

 On this set
$\f S$
we have a ``natural" smooth structure provided by the natural global chart
$x \byd \id : \f S \to \Rn \,.$

 Moreover, we can equip the set
$\f S$
with the further ``exotic" smooth structure provided by the global chart
$\ac x \byd x^3 \,.$

 If we regard
$\f S$
as an F--smooth space according to its natural smooth structure
$x \,,$
then the family of F--smooth curves
$\C C_x$
consists of all curves
$c : \f I_c \to \f S \,,$
such that the function
$x \com c : \f I_c \to \Rn$
are smooth.

 If we regard
$\f S$
as an F--smooth space according to its exotic smooth structure
$\ac x \,,$
then the family of F--smooth curves
$\C C_{\ac x}$
consists of all curves
$c : \f I_c \to \f S \,,$
such that the function
$\ac x \com c : \f I_c \to \Rn$
are smooth.

 Clearly, we have
\beq
\C C_x \sub \C C_{\ac x} \,.
\eeq

\myskip

 In an analogous way, we can equip
$\f S$
with other ``exotic" smooth structures in infinitely many ways and regard
$\f S$
as an F--smooth space in infinitely many ways according to each smooth structure.
 On the other hand, the union of all families of F--smooth curves obtained in this way equips
$\f S$
with a ``universal" F--smooth structure.

 Thus, conversely, the above construction provides an example of an F--smooth space, which can be regarded as a smooth manifold in infinitely many ways.\END
\eEx
%--------------------------------------------------------------------%
\chapter{Systems of maps}
\label{Systems of maps}
%--------------------------------------------------------------------%
\bsm
 First, we discuss the \emp{smooth systems} 
$(\f S, \eps)$
of \emp{smooth maps}
$f : \f M \to \f N \,.$
 Here, the ``space of parameters"
$\f S$
is a \emp{smooth manifold} and the ``evaluation map"
$\eps : \f S \car \f M \to \f N$
a \emp{smooth map}.

\myskip

 Then, we discuss the \emp{F--smooth systems}
$(\f S, \eps)$
of \emp{smooth maps}
$f : \f M \to \f N \,.$
 Here, the ``space of parameters"
$\f S$
is an \emp{F--smooth space} and the ``evaluation map"
$\eps : \f S \car \f M \to \f N$
an \emp{F--smooth map}.

\myskip

 The above notions are intended as an introduction to the more sophisticated notions of systems of sections discussed in the next 
Chapter \S\ref{Systems of sections}.

 The reader can find further discussions concerning the present subject in
\cite{ManMod83d}.
\esm
%--------------------------------------------------------------------%
%\newpage
\markboth{\rm Chapter \thechapter. Systems of maps}{\rm \thesection. Smooth systems of smooth maps}
\section{Smooth systems of smooth maps}
\label{Section: Smooth systems of smooth maps}
\markboth{\rm Chapter \thechapter. Systems of maps}{\rm \thesection. Smooth systems of smooth maps}
%--------------------------------------------------------------------%
\bsm
 We discuss the notion of \emp{smooth system of smooth maps} and the tangent prolongations of a smooth system of smooth maps as an introduction to the subsequent notions of smooth systems of smooth sections and smooth systems of smooth connections.
\esm
%--------------------------------------------------------------------%
\subsection{Smooth systems of smooth maps}
\label{Smooth systems of smooth maps}
%--------------------------------------------------------------------%
\bsm
 We start by defining the concept of \emp{smooth system of smooth maps} between two smooth manifolds.
 
 In simple words, a ``smooth system of smooth maps" is defined to be a family
$\{f\}$
of global smooth maps
$f : \f M \to \f N$
between two smooth manifolds
$\f M$
and
$\f N \,,$
which is smoothly parametrised by the points
$s \in \f S$
of a smooth manifold
$\f S \,.$
 
\myskip
 
 This simple concept will serve as an introduction to the further notions of ``smooth system of smooth sections" of a smooth double fibred manifold and ``smooth system of smooth connections" of a smooth fibred manifold
(see, later, 
Definition \ref{Definition: systems of sections}
and
Definition \ref{Definition: systems of connections}).
 
 Later, this concept will serve also as an introduction to the further more sophisticated notions of ``F--smooth system of smooth maps" between smooth manifolds, of ``F--smooth system of fibrewisely smooth sections" of a smooth double fibred manifold and of
``F--smooth system of fibrewisely F--smooth connections" of a smooth double fibred manifold
(see, later,
Definition \ref{Definition: F--smooth systems of smooth maps}
and
Definition \ref{Definition: F--smooth systems of smooth sections}).
\esm

 Let us consider two smooth manifolds
$\f M$
and
$\f N \,,$
and denote the set of global smooth maps between the two manifolds by
$\Map(\f M, \f N) \byd \wob f : \f M \to \f N \wcb \,.$

\myskip

\bDf
\label{Definition: smooth systems of smooth maps}
\index{smooth systems!system of maps}
\index{smooth systems!evaluation map}
\index{smooth systems!injective systems of maps}
 We define a \emp{smooth system of smooth maps},
between the smooth manifolds
$\f M$
and
$\f N \,,$
to be a pair
$(\f S, \eps) \,,$
where 

 1) $\f S$ 
is a smooth manifold,

 2) $\eps$
is a global smooth map, called \emp{evaluation map},
\beq
\eps : \f S \car \f M \to \f N \,.
\eeq

 Thus, the evaluation map
$\eps$
yields the map
\beq
\eps_\f S : \f S \to \Map(\f M, \f N) : 
s \mto \br s \,,
\eeq
where, for each
$s \in \f S \,,$ 
the global smooth map
$\br s$
is defined by
\beq
\br s : \f M \to \f N : m \mto \eps (s, m) \,.
\eeq

 Therefore, the map
$\eps_\f S :\f S \to \Map(\f M, \f N)$
provides a \emp{selection} of the global smooth maps
$\f M \to \f N \,,$
given by the subset
\beq
\Map_\f S(\f M, \f N) \byd \eps_\f S(\f S) \sub \Map(\f M, \f N) \,.
\eeq

\myskip

 The smooth system of smooth maps
$(\f S, \eps)$
is said to be \emp{injective} if the map
\beq
\eps_\f S : \f S \to \Map(\f M, \f N)
\eeq
is injective, i.e. if, for each
$s \,, \ac s \in \f S \,,$
\beq
\eps_\f S (s) = \eps_\f S(\ac s)
\Rarr
s = s' \,.
\eeq

 If the system is injective, then we obtain the bijection
\beq
\eps_\f S : \f S \to \Map_\f S(\f M, \f N) : s \mto \br s \,,
\eeq
whose inverse is denoted by
\beq
(\eps_\f S)^{-1} : 
\Map_\f S(\f M, \f N) \to \f S : f \mto \wha f \,.\END
\eeq
\eDf

\myskip

 Indeed, we are essentially interested in injective F--smooth systems of smooth maps.

\myskip

 We denote the local smooth charts of
$\f M \,,$
$\f N \,,$
$\f S$
respectively, by
\beq
(y^i) : \f M \to \Rn^{d_\f M} \,,
\qquad
(z^a) : \f N \to \Rn^{d_\f N} \,,
\qquad
(w^A) : \f S \to \Rn^{d_\f S} \,.
\eeq

\myskip

 We have the following elementary examples of smooth systems of smooth maps.

\myskip

\bEx
\label{Example: system of linear maps}
\index{smooth systems!system of linear maps}
 If
$\f M$
and
$\f N$
are vector spaces, then we obtain the injective smooth system of linear maps by setting
\beq
\f S \byd \lin(\f M, \f N) \,.
\eeq
 
 In this case, the smooth manifold
$\f S$
turns out to be a vector space of dimension
\beq
d_\f S = d_\f M \cdot d_\f N \,.
\eeq
 
 The choice of a basis of
$\f M$
and a basis of
$\f N$
yields the distinguished linear (global) charts
\beq
(y^i) : \f M \to \Rn^{d_\f M} \,,
\qquad
(z^a) : \f N \to \Rn^{d_\f N} \,,
\qquad
(w^a_i) : \f S \to \Rn^{d_\f M \cdot d_\f M} \,.
\eeq

 Then, the coordinate expression of
$\eps$
becomes
\beq
\eps^a = w^a_i \, y^i \,.\END
\eeq
\eEx

\bEx
\label{Example: system of affine maps}
\index{smooth systems!system of affine maps}
 If
$\f M$
and
$\f N$
are affine spaces, then we obtain the injective smooth system of affine maps by setting
\beq
\f S \byd \aff(\f M, \f N) \,.
\eeq
 
 In this case, the smooth manifold
$\f S$
turns out to be an affine space, of dimension
\beq
d_\f S = d_\f M \cdot d_\f N + d_\f N \,,
\eeq
which is associated with the vector space
$\aff (\f M, \baf N) \,,$
where
$\baf N$
is the vector space associated with the affine space
$\f N \,.$
  
 The choice of a basis and an origin of
$\f M$
and a basis and an origin of
$\f N$
yields the distinguished affine (global) charts
\beq
(y^i) : \f M \to \Rn^{d_\f M} \,,
\qquad
(z^a) : \f N \to \Rn^{d_\f N} \,,
\qquad
(w^a_i, w^a) : \f S \to 
\Rn^{d_\f M \cdot d_\f M} \car \Rn^{d_\f N} \,.
\eeq

 Accordingly, the coordinate expression of
$\eps$
becomes
\beq
\eps^a = w^a_i \, y^i + w^a \,.\END
\eeq
\eEx

\bEx
\label{Example: system of polynomial maps}
\index{smooth systems!system of polynomial maps}
 If
$\f M$
and
$\f N$
are affine spaces, then, analogously to the above 
Example \ref{Example: system of affine maps}, 
we can define 

 - the injective smooth system of polynomial maps of a given degree
$r \,,$
with
$0 \leq r \,,$

 - the injective smooth system of polynomial maps of any degree
$r \,,$
with
$0 \leq r \leq k \,,$
where
$k$
is a given positive integer.\END
\eEx

\myskip

 All examples above deal with finite dimensional smooth systems of maps, as it is implicitly requested in
Definition \ref{Definition: smooth systems of smooth maps}.

 On the other hand, we can easily extend the concept of smooth system of smooth maps between two smooth manifolds, by considering an infinite dimensional system, which is the direct limit of finite dimensional systems, according to the following
Example \ref{Example: infinite dimensional system of polynomial maps of any degree}.
 In a sense, this example is intermediate between the finite dimensional case and the infinite dimensional case.

\myskip

\bEx
\label{Example: infinite dimensional system of polynomial maps of any degree}
\index{smooth systems!system of polynomial maps}
 If
$\f M$
and
$\f N$
are affine spaces, then we obtain the injective (infinite dimensional) smooth system of all polynomial maps by considering the set
\beq
\f S \byd \pln(\f M, \f N) \,,
\eeq
consisting of all polynomials
$\f M \to \f N$
of any degree
$r \,,$
with
$0 \leq r < \infty \,.$

 Later, we shall see that such a smooth system has a natural F--smooth structure
(see, later,
Definition \ref{Definition: F--smooth spaces}).\END
\eEx

\myskip

 Eventually we compare three simple examples in order to emphasise further features of smooth systems of global smooth maps.
 
\myskip

\bEx
\label{Example: smooth systems of maps of Rn}
 Let us consider the smooth manifolds
\beq
\f M \byd \Rn \,,
\qquad
\f N \byd \Rn \,,
\qquad
\f S \byd \Rn \,,
\eeq
equipped with their standard charts
\beq
(y) : \f M \to \Rn \,,
\qquad
(z) : \f N \to \Rn \,,
\qquad
(w) : \f S \to \Rn \,.
\eeq

 Moreover, let us consider three smooth systems of smooth maps defined, respectively, by the following evaluation maps
\beq
\eps_1 : \f S \car \f M \to \f N \,,
\qquad
\eps_2 : \f S \car \f M \to \f N \,,
\qquad
\eps_3 : \f S \car \f M \to \f N \,,
\eeq
with coordinate expressions
\beq
z \com \eps_1 = w \, y \,,
\qquad
z \com \eps_2 = w^3 \, y \,,
\qquad
z \com \eps_3 = w^2 \, y \,.
\eeq

 Then, 
 
 - $(\f S, \eps_1)$ 
is the system of linear maps
$f : \f M \to \f N \,,$

 - $(\f S, \eps_2)$ 
is the system of linear maps
$f : \f M \to \f N \,,$

 - $(\f S, \eps_3)$ 
is the system of linear maps
$f : \f M \to \f N$
with coefficients
$w^2 \geq 0 \,.$

 Thus, the systems
$(\f S, \eps_1)$
and
$(\f S, \eps_2)$
select the same smooth maps, but with different smooth parametrisations.
 In these cases we have assumed the same smooth structure of
$\f S \,,$
but the transition between the two parametrisations is non smooth in one way.

 Moreover, the systems
$(\f S, \eps_1)$
and
$(\f S, \eps_2)$
are injective, while the system
$(\f S, \eps_3)$
is non injective.\END
\eEx

\myskip

 It is useful to compare the above systems
$(\f S, \eps_1)$
and
$(\f S, \eps_2)$
with the system
$(\f S, \eps)$
discussed in the forthcoming
Example \ref{Example: different smooth structures of a system of maps}.

\myskip

 Eventually, we exhibit a further weird example of F--smooth systems of smooth maps.
\myskip

\bEx
\label{Example: m + s u}
 Let us consider two vector spaces
$\f M$
and
$\f N \,,$
along with the smooth manifold
$\f S \byd \lin(\f M, \f N)$
and the smooth map
\beq
\eps : \f S \car \f M \to \f N : (s, m) \mto n + s (X) \,,
\eeq
where
$n \in \f N$
and
$X \in \f M$
is a non vanishing given element.

 Then, the pair
$(\f S, \eps)$
turns out to be an injective smooth system of smooth maps.\END
\eEx
%--------------------------------------------------------------------%
\subsection{Smooth structure of $(\f S,\eps)$}
\label{Smooth structure of (S,eps)}
%--------------------------------------------------------------------%
\bsm
 In our definition of ``smooth system of smooth maps" between two smooth manifolds
(see
Definition \ref{Definition: smooth systems of smooth maps})
we have required a priori that the set of parameters
$\f S$
be finite dimensional and smooth and that the evaluation map
$\eps$
be smooth as well.

 On the other hand, we might ask whether the fact that the manifolds
$\f M$
and
$\f N$
are smooth and that the maps
$f : \f M \to \f N$
selected by the system are smooth allows us to recover uniquely the smooth structure of
$\f S \,.$

 The answer to the above question is negative.
 Here, we do not fully address this question.
 But, we present a simple example
(see
Example \ref{Example: different smooth structures of a system of maps}),
where we show that, if
$\f S$
admits a finite dimensional smooth structure compatible with
$\eps \,,$
then this structure needs not to be unique.

 Moreover, we observe that, if we do not assume a priori a finite dimensional smooth structure on
$\f S \,,$
then it might be that no finite dimensional smooth structure at all could be recovered on
$\f S \,.$
 To prove this, just consider the system
$(\f S, \eps) \,,$
where
$\f S$
is the set of \emp{all} smooth maps
$f : \f M \to \f N$
(see, later,
\S\ref{F--smooth systems of smooth maps}).

\myskip

 The present question might arise also in comparison with a result which will be achieved later, in the next Section, in the context of ``F--smooth systems of smooth maps", where we do not assume a priori any smooth structure of
$\f S \,,$
but we recover uniquely its F--smooth structure
(see
Definition \ref{Definition: F--smooth systems of smooth maps}
and
Theorem \ref{Theorem: F--smooth structure of F--smooth systems of smooth maps}).
\esm

\bEx
\label{Example: different smooth structures of a system of maps}
 Let us consider the following smooth system
$(\f S, \eps)$
of smooth maps between smooth manifolds.
 
 Let us consider the manifolds
\beq
\f M \byd \Rn \,,
\qquad
\f N \byd \Rn \,,
\qquad
\f S \byd \Rn \,.
\eeq

 We consider the ``natural" smooth structures of
$\f M$
and
$\f N$
induced by their ``natural" charts
\beq
y : \f M \to \Rn
\ssep{and}
z : \f N \to \Rn \,.
\eeq

 On the other hand, we consider the ``natural" smooth structure and the ``exotic" smooth structure of
$\f S$
given respectively by the natural and the exotic charts
\beq
w : \f S \to \Rn
\ssep{and}
\ac w \byd w^3 : \f S \to \Rn \,.
\eeq

 Moreover, let us consider the evaluation map
$\eps : \f S \car \f M \to \f N \,,$
whose coordinate expression in the above charts reads, respectively, as
\beq
z \com \eps = w^3 \, y
\ssep{and}
z \com \eps = \ac w \, y \,.
\eeq

 Indeed, the evaluation map
$\eps$
turns out to be the same and smooth with respect to the above different smooth structures of
$\f S \,.$\END
\eEx

\myskip

 It is useful to compare the previous systems
$(\f S, \eps_1)$
and
$(\f S, \eps_2) \,,$
discuss in 
Example \ref{Example: smooth systems of maps of Rn},
with the above system
$(\f S, \eps) \,.$

 Indeed, in the first two systems above we deal with the same smooth structure of
$\f S$
and with different parametrisations of the selected global smooth maps.

 Conversely, in the third system we deal with a different smooth structure of
$\f S \,.$

 On the other hand, these three systems select the same global smooth maps in different ways.
%--------------------------------------------------------------------%
\subsection{Smooth tangent prolongation of $(\f S, \eps)$}
\label{Smooth tangent prolongation of (S, eps)}
%--------------------------------------------------------------------%
\bsm
 Given a smooth system of smooth maps
$(\f S, \eps)$
between two smooth manifolds
$\f M$
and
$\f N \,,$
we discuss 
 
 - the smooth system
$(T\f S, T\eps)$
of smooth maps between the smooth manifolds
$T\f M$
and
$T\f N \,,$
which is achieved via the tangent prolongation
$T\eps$
of
$\eps$
with respect to the both factors
$\f S$
and
$\f M \,,$
 
 - the smooth system
$(T\f S, T_1\eps)$
of smooth maps between the smooth manifolds
$\f M$
and
$T\f N \,,$
which is achieved via the tangent prolongation
$T_1\eps$
of
$\eps$
with respect to the 1st factor
$\f S \,,$ 

 - the smooth system
$(\f S, T_2\eps)$
of smooth maps between the smooth manifolds
$T\f M$
and
$T\f N \,,$
which is achieved via the tangent prolongation
$T_2\eps$
of
$\eps$
with respect to the 2nd factor
$\f M \,.$

 We stress that, if the system
$(\f S, \eps)$
is injective, then its tangent prolongations
$(T\f S, T\eps)$
and
$(T\f S, T_1\eps)$
need not to be injective.

\myskip

 Later
(see
\S\ref{F--smooth systems of smooth maps}), 
we shall be involved with a generalised concept of system, where
$\f S$
is no longer a finite dimensional smooth manifold, hence we cannot avail of the standard approach to define its tangent space
$T\f S \,.$
 Actually, we will achieve the tangent space of
$\f S$
by an indirect procedure, via smooth maps between smooth manifolds.
 Indeed, the tangent prolongations of the system 
$(\f S, \eps)$
discussed in the present Section will be a hint for the more sophisticated cases discussed in the next Chapter
(see
Section \ref{F--smooth tangent prolongation of (S, eps)}
and
Section \ref{F--smooth tangent prolongation of (S, zet, eps)}).
\esm

 Let us consider two smooth manifolds
$\f M$
and
$\f N \,,$
and a smooth system of smooth maps
$(\f S, \eps)$
between these manifolds
(see 
Definition \ref{Definition: smooth systems of smooth maps}).

 We denote the charts of
$\f M \,,$
$\f N \,,$
$\f S \,,$
$T\f M \,,$
$T\f N \,,$
$T\f S$
respectively, by
\bgt
(y^i) : \f M \to \Rn^{d_\f M} \,,
\qquad
(y^i, \dt y^i) : T\f M \to \Rn^{2 d_\f M} \,,
\\
(z^a) : \f N \to \Rn^{d_\f N} \,,
\qquad
(z^a, \dt z^a) : T\f N \to \Rn^{2 d_\f N} \,,
\\
(w^A) : \f S \to \Rn^{d_\f S} \,,
\qquad
(w^A, \dt w^A) : T\f S \to \Rn^{2 d_\f S} \,.
\end{gather*}

\myskip

\bPr
\label{Proposition: tangent prolongations of systems of maps}
\index{smooth systems!1st partial tangent prolongation of systems of maps}
\index{smooth systems!2nd partial tangent prolongation of systems of maps}
\index{smooth systems!tangent prolongation of systems of maps}
 We consider the following \emp{tangent prolongations} of the system.

\myskip

 1) We have the \emp{total tangent prolongation} of
$\eps \,,$
which yields the smooth map
\beq
T\eps : T\f S \car T\f M \to T\f N \,,
\eeq
with coordinate expression
\beq
(z^a \com T\eps) = 
\eps^a
\ssep{and}
(\dt z^a \com T\eps) = 
\der_A\eps^a \, \dt w^A + \der_i \eps^a \, \dt y^i \,.
\eeq

\myskip

 2) We have the \emp{partial tangent prolongation} of
$\eps \,,$
\emp{with respect to the 1st factor}
which yield the smooth map
\beq
T_1\eps : T\f S \car \f M \to T\f N \,,
\eeq
with coordinate expression
\beq
(z^a \com T_1\eps) = 
\eps^a
\ssep{and}
(\dt z^a \com T_1\eps) = 
\der_A\eps^a \, \dt w^A \,.
\eeq 

 3) We have the \emp{partial tangent prolongation} of
$\eps \,,$
\emp{with respect to the 2nd factor}, which yields the smooth map
\beq
T_2\eps : \f S \car T\f M \to T\f N \,,
\eeq
with coordinate expression
\beq
(z^a \com T_2\eps) = 
\eps^a
\ssep{and}
(\dt z^a \com T_2\eps) = 
\der_i \eps^a \, \dt y^i \,.
\eeq

\myskip

 Indeed, the following diagrams commute
\newdiagramgrid{Square 3}
{1, 1, 1, 1, 1, 1, 1, 1, 1, 1, .4}
{.7, .7}
\bdg[grid=Square 3]
T\f S \car T\f M &\rTo^{T\eps} &T\f N
&&T\f S \car \f M &\rTo^{T_1\eps} &T\f N 
&&\f S \car T\f M &\rTo^{T_2\eps} &T\f N
\\
\dTo^{} &&\dTo_{}
&&\dTo^{} &&\dTo_{}
&&\dTo^{} &&\dTo_{}
\\
\f S \car \f M &\rTo^{\eps} &\f N
&&\f S \car \f M &\rTo^{\eps} &\f N &
&\f S \car \f M &\rTo^{\eps} &\f N &.
\edg

\myskip

 Thus, according to the above commutative diagrams,

 - the pair
$(T\f S, T\eps)$
turns out to be a smooth system of smooth maps between the smooth manifolds
$T\f M$
and
$T\f N \,,$

 - the pair
$(T\f S, T_1\eps)$
turns out to be a smooth system of smooth maps between the smooth manifolds
$\f M$
and
$T\f N \,,$

 - the pair
$(\f S, T_2\eps)$
turns out to be a smooth system of smooth maps between the smooth manifolds
$T\f M$
and
$T\f N \,.$

 All above prolonged smooth systems project over the smooth system
$(\f S, \eps) \,.$ 

 Indeed, the smooth system
$(T\f S, T\eps)$
is characterised by the smooth system
$(T\f S, T_2\eps)$
and, conversely, the smooth system
$(T\f S, T_2\eps)$
is characterised by the smooth system
$(T\f S, T\eps)$
in virtue of the equalities
\beq
T\eps = T_1\eps + T_2\eps
\ssep{and}
T_2\eps = T\eps - T_1\eps \,.\END
\eeq
\ePr

\bCr
\label{Corollary: injectivity of prolonged system}
\index{smooth systems!injective system of maps}
 If the smooth system
$(\f S, \eps)$
is injective, then also the smooth system 
$(\f S, T_2\eps)$
turns out to be injective, in virtue of the projectability on
$(\f S, \eps) \,.$\END
\eCr

\bPr
\label{Proposition: injectivity of the tangent prolongations of a system of maps}
\index{smooth systems!injective system of maps}
 The following implications hold:
 
 - if the smooth system
$(T\f S, T\eps)$
is injective then the smooth system
$(\f S, \eps)$
is injective as well,

 - if the smooth system
$(T\f S, T_1\eps)$
is injective then the smooth system
$(\f S, \eps)$
is injective as well.
\ePr

\bpf
 In fact, the smooth system
$(\f S, \eps)$
is obtained by projection of its tangent prolongations, according to the commutative diagrams in
Proposition \ref{Proposition: tangent prolongations of systems of maps}.\QED
\epf

\bRm
\label{Remark: 1st tangent prolongation and tangent prolongation of a system of maps need not to be injective}
\index{smooth systems!injective system of maps}
 If the smooth system
$(\f S, \eps)$
is injective, then its tangent prolongations
$(T\f S, T\eps)$
and
$(T\f S, T_1\eps)$
need not to be injective.\END
\eRm

\myskip

 Now, let us examine the tangent prolongations of two elementary examples of smooth systems of smooth maps.
 
\myskip

\bEx
\label{Example: tangent prolongations of systems of linear maps}
\index{smooth systems!system of linear maps}
 Let us refer to the smooth system
$(\f S, \eps)$
of linear maps discussed in
Example \ref{Example: system of linear maps}.

 Then, the coordinate expressions of the smooth systems
\beq
T\eps : T\f S \car T\f M \to T\f N
\ssep{and}
T_1\eps : T\f S \car \f M \to T\f N
\eeq
are
\bat{4}
z^a \com T\eps
&=
w^a_i \, y^i
&&\ssep{and}
\dt z^a \com T\eps 
&&= 
\dt w^a_i \, y^i + w^a_i \, \dt y^i \,,
\\
z^a \com T_1\eps
&=
w^a_i \, y^i
&&\ssep{and}
\dt z^a \com T_1\eps 
&&= 
\dt w^a_i \, y^i \,.
\end{alignat*}

 Indeed, the above coordinate expressions show that the two tangent prolongations are injective.\END
\eEx

\bEx
\label{Example: tangent prolongations of systems of affine maps}
\index{smooth systems!system of affine maps}
 Let us refer to the smooth system
$(\f S, \eps)$
of affine maps discussed in
Example \ref{Example: system of affine maps}.

 Then, the coordinate expressions of the smooth systems
\beq
T\eps : T\f S \car T\f M \to T\f N
\ssep{and}
T_1\eps : T\f S \car \f M \to T\f N
\eeq
are
\bat{4}
z^a \com T\eps
&=
w^a_i \, y^i + w^a
&&\ssep{and}
\dt z^a \com T\eps 
&&= 
\dt w^a_i \, y^i + w^a_i \, \dt y^i + \dt w^a \,,
\\
z^a \com T_1\eps
&=
w^a_i \, y^i + w^a
&&\ssep{and}
\dt z^a \com T_1\eps 
&&= 
\dt w^a_i \, y^i + \dt w^a \,.
\end{alignat*}

 Indeed, the above coordinate expressions show that the two tangent prolongations are injective.\END
\eEx

\myskip

 In a similar way we can define 
(see
Example \ref{Example: system of polynomial maps})
 
 - the tangent prolongations of the smooth system of polynomial maps of a given degree
$r \,,$
with
$0 \leq r \,,$

 - the tangent prolongations of the system of polynomial maps of any degree
$r \,,$
with
$0 \leq r \leq k \,,$

 - the tangent prolongations of the smooth system of polynomial maps of any degree
$r \,,$
with
$0 \leq r \leq \infty \,.$
%--------------------------------------------------------------------%
%\newpage
\markboth{\rm Chapter \thechapter. Systems of maps}{\rm \thesection. F--smooth systems of smooth maps}
\section{F--smooth systems of smooth maps}
\label{Section: F--smooth systems of smooth maps}
\markboth{\rm Chapter \thechapter. Systems of maps}{\rm \thesection. F--smooth systems of smooth maps}
%--------------------------------------------------------------------%
\bsm
 We discuss the \emp{F--smooth systems} 
$(\f S, \eps)$
of smooth maps
$f : \f M \to \f N$
between smooth manifolds and its F--smooth tangent prolongation
$(\E T\f S, \E T\eps) \,.$

 Moreover, we compare the F--smooth systems of smooth maps with the smooth systems of smooth maps.
\esm
%--------------------------------------------------------------------%
\subsection{F--smooth systems of smooth maps}
\label{F--smooth systems of smooth maps}
%--------------------------------------------------------------------%
\bsm
 In the previous 
Section \ref{Section: Smooth systems of smooth maps},
we have discussed the concept of a ``smooth system of smooth maps" 
$(\f S, \eps)$
between two smooth manifolds
$\f M$
and
$\f N$
(see
Section \ref{Systems of maps}).

 Here, we introduce the generalised notion of ``\emp{F--smooth system of smooth maps}"
$(\f S, \eps)$
between two smooth manifolds
$\f M$
and
$\f N \,,$
by releasing the hypotheses of smoothness and finite dimension of
$\f S \,.$
 
\myskip

 The general concept of F--smooth space, which has been introduced in the above Chapter
(see
Section \ref{Section: F--smooth spaces}) 
can be applied to a large spectrum of contexts.
 
 On the other hand, for our purposes, the most interesting examples of F--smooth spaces are the F--smooth spaces of smooth maps between two smooth manifolds.
 Even more, this kind of examples provide the true reason of our interest on F--smooth spaces in the present report.
 
  To be more precise, our main interest deals with the particular case  of F--smooth systems of fibrewisely smooth sections, which will be introduced in
Section \ref{F--smooth systems of smooth sections}.
 Thus, here the concept of ``F--smooth system of smooth maps" is intended as an introduction to the more sophisticated concept of 
``F--smooth system of fibrewisely smooth sections".

\myskip
 
 Indeed, a geometric approach to the space of smooth maps between smooth manifolds, in terms of the standard differential geometry, would possibly involve subtle and hard problems concerning infinite dimensional smooth manifolds.
 Conversely, the spaces of smooth maps between two smooth manifolds, regarded as F--smooth spaces, allow us to achieve several geometric results, even if the structure of F--smooth space is weaker than that of smooth manifold (see, for instance \cite{CabKol95,CabKol98,KolMod98}).
 
  In the next
Section \ref{F--smooth tangent prolongation of (S, eps)}, 
as a starting example of such geometric constructions, we sketch the tangent space of F--smooth spaces of global smooth maps.
\esm

 Let us consider two smooth manifolds
$\f M$
and
$\f N \,.$

\myskip

 The following Definition provides a generalisation of the concept of smooth system of smooth maps
(see
Definition \ref{Definition: smooth systems of smooth maps}),
as here we do not require that
$\f S$
be a finite dimensional smooth manifold
(hence, that the map
$\eps$
be smooth).

 Even more, we do not assume a priori any kind of smoothness on the set
$\f S \,,$
but later we will uniquely recover its F--smooth structure by
Theorem \ref{Theorem: F--smooth structure of F--smooth systems of smooth maps}.
 Indeed, the specification ``\emp{F--smooth}" system, which is anticipated in the following
Definition \ref{Definition: F--smooth systems of smooth maps},
will be justified later by this 
Theorem.

\myskip

\bDf
\label{Definition: F--smooth systems of smooth maps}
\index{F--smooth systems!F--smooth system of smooth maps} 
\index{F--smooth systems!injective F--smooth system of smooth maps} 
 We define an \emp{F--smooth system of smooth maps}
between the smooth manifolds
$\f M$
and
$\f N$
to be a pair
$(\f S, \eps) \,,$
where 

 1) $\f S$
is a set,

 2) $\eps : \f S \car \f M \to \f N$
is a map, called \emp{evaluation map}, such that, for each
$s \in \f S \,,$
the induced map
\beq
\eps_s : \f M \to \f N : m \mto \eps(s,m)
\eeq
is \emp{globally defined and smooth}.

 Thus, the evaluation map
$\eps$
yields the map
\beq
\eps_\f S : \f S \to \Map(\f M, \f N) : 
s \mto \br s \,,
\eeq
where, for each
$s \in \f S \,,$ 
the global smooth map
$\br s$
is defined by
\beq
\br s : \f M \to \f N : m \mto \eps (s, m) \,.
\eeq

 Therefore, the map
\beq
\eps_\f S :\f S \to \Map(\f M, \f N)
\eeq
provides a \emp{selection} of the global smooth maps
$\f M \to \f N \,,$
given by the subset
\beq
\Map_\f S(\f M, \f N) \byd \eps_\f S(\f S) \sub \Map(\f M, \f N) \,.
\eeq

\myskip

 The F--smooth system of smooth maps
$(\f S, \eps)$
is said to be \emp{injective} if the map
\beq
\eps_\f S : \f S \to \Map(\f M, \f N)
\eeq
is injective, i.e. if, for each
$s \,, \ac s \in \f S \,,$
\beq
\eps_\f S (s) = \eps_\f S(\ac s)
\Rarr
s = s' \,.
\eeq

 If the system is injective, then we obtain the bijection
\beq
\eps_\f S : \f S \to \Map_\f S(\f M, \f N) : s \mto \br s \,,
\eeq
whose inverse is denoted by
\beq
(\eps_\f S)^{-1} : 
\Map_\f S(\f M, \f N) \to \f S : f \mto \wha f \,.\END
\eeq
\eDf

\myskip

 Indeed, we are essentially interested in injective F--smooth systems of smooth maps.

\myskip

 Next, we prove that any system of smooth maps has a natural F--smooth structure.
 Actually, this property justifies the fact that we have anticipated the F--smoothness in
Definition \ref{Definition: F--smooth systems of smooth maps}.

\myskip

\bTh
\label{Theorem: F--smooth structure of F--smooth systems of smooth maps}
\index{F--smooth systems!F--smooth system of maps}
 Let us consider an F--smooth system 
$(\f S, \eps)$
of global smooth maps and let
\beq
\C C \byd \{c : \f I_c \to \f S\} \,,
\eeq
be the set consisting of \emp{all} curves, such that the induced map
\beq
c^*(\eps) : \f I_c \car \f M \to \f N :
(\lam, m) \mto \eps\w(c(\lam),m\w)
\eeq
be smooth.

 Then, the pair
$(\f S, \C C)$
turns out to be an F--smooth space.
\eTh

\bpf
 Let us check that
$\C C$
fulfills the two requirements of 
Definition \ref{Definition: F--smooth spaces}.

 1) For each
$s \in \f S \,,$
the constant curve
$c : \Rn \to \f S : \lam \mto s$
yields the map
\beq
c^*(\eps) : \f I_c \car \f M \to \f N :
(\lam, m) \mto \eps\w(c(\lam),m\w) = \eps(s,m) \,,
\eeq
which is smooth in virtue of condition 2) in
Definition \ref{Definition: F--smooth systems of smooth maps}.
 Hence,
$c$
turns out to be an element of
$\C C$
passing through
$s \,.$

 2) Let us consider a curve
$c \in \C C$
and a smooth curve
$\gam : \f I_\gam \to \f I_c \,.$

 Then, the induced map
\beq
(c \com \gam)^* (\eps) : \f I_\gam \car \f M \to \f N :
(\lam, m) \mto \eps\W(c\w(\gam(\lam)\w), \, m\W)
\eeq
is the composition of two smooth maps
\newdiagramgrid{1}
{1.5, 1.5, 1.2, 1.2, .4}
{1}
\bdg[grid=1]
\f I_\gam \car \f M &\rTo^{(\gam, \id_\f M)} 
&\f I_c \car \f M
&\rTo^{c^*(\eps)} &\f N &,
\edg
hence it turns out to be smooth.\QED
\epf

\bCr
 Let us consider an F--smooth system 
$(\f S, \eps)$
of smooth maps.

 Then, the map
$\eps : \f S \car \f M \to \f N$
turns out to be F-smooth, with reference to the F-smooth structure of
$\f S$
and the natural F-smooth structures of
$\f M$
and
$\f N$
(see the above
Theorem \ref{Theorem: F--smooth structure of F--smooth systems of smooth maps},
Theorem \ref{Theorem: relation between smoothness and F--smoothness},
Definition \ref{Definition: natural F--smooth structure of smooth manifolds}
and
Proposition \ref{Proposition: F--smooth cartesian product}).
\eCr

\bpf
 By definition
(see
Definition \ref{Definition: F--smooth maps}),
the map
$\eps : \f S \car \f M \to \f N$
is F--smooth if, for each F--smooth curve
$c : \f I_c \to \f S \car \f M \,,$
the curve
$\eps \com c : \f I_c \to \f N \,,$
given by the commutative diagram
\newdiagramgrid{2x2}
{1.5, 1.5, .6}
{.7, .7}
\bdg[grid=2x2]
\f S \car \f M &\rTo^{\eps} &\f N
\\
\uTo^{c} &&\dTo_{\id_\f N}
\\
\f I_c &\rTo^{\eps \com c} &\f N &,
\edg
is smooth.

 Thus, let us consider any F--smooth curve
$c : \f I_c \to \f S \car \f M \,.$

 In virtue of
Proposition \ref{Proposition: F--smooth cartesian product}
and
Theorem \ref{Theorem: relation between smoothness and F--smoothness},
we can write
\beq
c \eqv (c_\f S, c_\f M) : \f I_c \to \f S \car \f M \,,
\eeq
where
$c_\f S : \f I_c \to \f S$
is an F--smooth curve and
$c_\f M : \f I_c \to \f M$
is a smooth curve.

 Let us define the smooth map
\beq
c^*_\f S(\eps) : \f I_c \car \f M \to \f N : 
(\lam, m) \mto \eps \w(c_\f S(\lam), m\w) \,.
\eeq

 Then, the curve
\beq
\eps \com c = \eps \com (c_\f S, c_\f M) : \f I_c \to \f N :
\lam \mto \eps\w(c_\f S(\lam), \, c_\f M(\lam)\w) =
c^*_\f S (\eps) (\lam, m) \,,
\eeq
which is a composition of smooth maps according to the following digram
\newdiagramgrid{2x2}
{1.3, 1.7, .7}
{.7, .7}
\bdg[grid=2x2]
\f I_c &\rTo^{(\id, c_\f M)} &\f I_c \car \f M
\\
\uTo^{\id} &&\dTo_{c^*_\f S(\eps)}
\\
\f I_c &\rTo^{\eps \com c} &\;\f N &,
\edg
turns out to be smooth hence F--smooth
(see
Theorem \ref{Theorem: F--smooth structure of F--smooth systems of smooth maps}).\QED
\epf

\myskip

 Now, we provide examples of F--smooth systems of smooth maps, by starting with ``infinite dimensional" examples of systems of smooth maps. 
 
\myskip

\bEx
\label{Example: F--smooth system of all global smooth maps}
\index{F--smooth systems!F--smooth system of maps}
 Let us consider two smooth manifolds
$\f M$
and
$\f N \,.$

 Then, the set
\beq
\f S \byd \Map(\f M, \f N) \byd \wob f : \f M \to \f N \wcb
\eeq
consisting of \emp{all global smooth maps} 
$f : \f M \to \f N$
yields an F--smooth system of smooth maps.\END
\eEx

\bEx
\label{Example: F--smooth system of all smooth maps with compact support}
\index{F--smooth systems!F--smooth system of maps}
 Let us consider a smooth manifold
$\f M$
and a vector space
$\f N \,.$

 Then, the subset
\beq
\f S \byd \Map_{\cpt}(\f M, \f N) \sub \Map(\f M, \f N) \byd \wob f : \f M \to \f N \wcb
\eeq
consisting of \emp{all global smooth maps} 
$f : \f M \to \f N$
\emp{with compact support} yields an F--smooth system of smooth maps.

 Indeed, this system is a subsystem of the system of all global smooth maps
$\f M \to \f N$
considered in the above 
Example \ref{Example: F--smooth system of all smooth maps with compact support}.\END
\eEx

\myskip

 We can reconsider, in the present context of F--smooth spaces,
``finite dimensional" examples of systems of global smooth maps already discussed in
Section \ref{Smooth systems of smooth maps}.
 
\myskip

\bEx
\label{Example: F--smooth system of linear maps}
\index{F--smooth systems!F--smooth system of maps}
 According to
Example \ref{Example: system of linear maps},
if
$\f M$
and
$\f N$
are vector spaces, then we obtain the F--smooth system of linear maps by setting
\beq
\f S \byd \lin(\f M, \f N) \,.
\eeq

 Indeed, the F--smooth structure of
$\f S$
turns out to be just the natural F--smooth structure underlying the smooth structure of
$\f S \,,$
according to
Theorem \ref{Theorem: F--smooth structure of F--smooth systems of smooth maps}
and
Example \ref{Example: system of linear maps}.\END
\eEx

\bEx
\label{Example: F--smooth system of affine maps}
\index{F--smooth systems!F--smooth system of maps}
 According to
Example \ref{Example: system of affine maps},
if
$\f M$
and
$\f N$
are affine spaces, then we obtain the F--smooth system of affine maps by setting
\beq
\f S \byd \aff(\f M, \f N) \,.
\eeq

 Indeed, the F--smooth structure of
$\f S$
turns out to be just the natural F--smooth structure underlying the smooth structure of
$\f S \,,$
according to
Theorem \ref{Theorem: F--smooth structure of F--smooth systems of smooth maps}
and
Example \ref{Example: system of affine maps}.\END
\eEx

\bEx
\label{Example: F--smooth system of polynomial maps}
\index{F--smooth systems!F--smooth system of maps}
 According to
Example \ref{Example: system of polynomial maps},
if
$\f M$
and
$\f N$
are affine spaces, then we obtain the F--smooth systems
$\f S$
of polynomial maps of a given degree
$r$
and the F--smooth system
$\f S$
of polynomial maps of all degrees less than a given integer
$k \,.$

 Indeed, the F--smooth structures of the above systems of smooth maps
turn out to be just the natural F--smooth structures underlying the smooth structures of
$\f S \,,$
according to
Theorem \ref{Theorem: F--smooth structure of F--smooth systems of smooth maps}
and
Example \ref{Example: system of polynomial maps}.\END
\eEx

\myskip

 Eventually, we revisit in the present context of F--smooth spaces the ``infinite dimensional" example of polynomial maps of any degree. 
 
\myskip

\bEx
\label{Example: F--smooth system of polynomial maps of any degree}
\index{F--smooth systems!F--smooth system of maps}
 According to
Example \ref{Example: system of polynomial maps},
if
$\f M$
and
$\f N$
are affine spaces, then we obtain the F--smooth system
$\f S$
of polynomial maps of any degree.\END
\eEx

\myskip

 Further, we consider two examples which play an intermediate role between the F--smooth systems of smooth maps and the F--smooth systems of fibrewisely smooth sections, which will be discussed later
(see
\S\ref{F--smooth systems of smooth sections}).

\myskip

\bEx
\label{Example: F--smooth system of global smooth sections}
\index{F--smooth systems!F--smooth system of maps}
 Let us suppose that the smooth manifold
$\f N$
be a smooth bundle
$q : \f N \to \f M$
over the smooth manifold
$\f M \,.$

 Then,  the subset
\beq
\f S \byd \Sec(\f M, \f N) \sub \Map(\f M, \f N)
\eeq
consisting of \emp{all global smooth sections} 
$\phi : \f M \to \f N$
is an F--smooth system of smooth maps.

 Indeed, this system is a subsystem of the system of all smooth maps
$\f M \to \f N$
considered in the 
Example \ref{Example: F--smooth system of all global smooth maps}.\END
\eEx

\bEx
\label{Example: F--smooth system of global smooth sections with compact support}
\index{F--smooth systems!F--smooth system of maps}
 Let us suppose that the smooth manifold
$\f N$
be a vector space and a smooth bundle
$q : \f N \to \f M$
over the smooth manifold
$\f M \,.$

 Then,  the subset
\beq
\f S \sub \Sec(\f M, \f N) \sub \Map(\f M, \f N)
\eeq
consisting of \emp{all global smooth sections with compact support} 
$\phi : \f M \to \f N$
is an F--smooth system of smooth maps.

 Indeed, this system is a subsystem of the above system of all global smooth sections considered in the above 
Example \ref{Example: F--smooth system of global smooth sections}.\END
\eEx
%--------------------------------------------------------------------%
\subsection{F--smooth tangent prolongation of $(\f S, \eps)$}
\label{F--smooth tangent prolongation of (S, eps)}
%--------------------------------------------------------------------%
\bsm
 Now, we consider an F--smooth system
$(\f S, \eps)$
of smooth maps between two smooth manifolds
$\f M$
and
$\f N$
(see
Definition \ref{Definition: F--smooth systems of smooth maps})
and introduce the concept of ``\emp{F--smooth tangent space}"
$\E T\f S$
of the F--smooth space of parameters
$\f S \,.$

 Our formal construction of the tangent space
$\E T\f S$
reflects the intuitive idea, by which, for each
$s \in \f S \,,$
a tangent vector
$\E X_s \in \E T_s\f S$
is to be an ``infinitesimal variation" of the global smooth map
$\eps_s : \f M \to \f N \,.$
 It is remarkable the fact that this construction involves only global smooth maps between smooth manifolds, by exploiting the smooth structures of the manifolds
$\f M$
and
$\f N$
and the smoothness of the selected maps
$f : \f M \to \f N$
of the system.

\myskip

 Actually, we define a tangent vector
$\E X_s$
of
$\f S \,,$
at
$s \in \f S \,,$
to be an equivalence class of F--smooth curves
$\wha c : \f I_{\wha c} \to \f S \,,$
such that the induced smooth maps between smooth manifolds
$\wha c^*(\eps) : \f I_c \car \f M \to \f N$
have a 1st order contact in
$s \,.$

 The definition of 1st order contact shows that every
$\E X_s \in \E T_s\f S$
turns out to be represented by suitable smooth map
$\Xi_s : \f M \to T\f N \,,$
or that it can, equivalently, be represented by a suitable smooth map
$\ol\Xi_s : T\f M \to T\f N \,.$
 Both
$\Xi_s$
and
$\ol\Xi_s$ 
project over the smooth map
$\eps_s : \f M \to \f N \,.$

 Moreover, we show that the tangent space
$\E T\f S$
is equipped with a natural F--smooth structure and that the natural maps
\beq
\tau_\f S : \E T\f S \to \f S
\ssep{and}
\E T_1\eps : \E T\f S \car \f M \to T\f N
\eeq
turn out to be F--smooth.

 Indeed, the fibres of
$\tau_\f S : \E T\f S \to \f S$
are naturally equipped with a vector structure.
 
\myskip

 This discussion on the F--smooth tangent space
$\E T\f S$
of an F--smooth system of smooth maps between smooth manifolds is intended as an introduction to the more sophisticated case of the tangent space
$\E T\f S$
of an F--smooth system of fibrewisely smooth sections
(see
Section \ref{F--smooth tangent prolongation of (S, zet, eps)}).
\esm

 Thus, let us consider two smooth manifolds
$\f M$
and
$\f N$
and  denote the smooth charts of
$\f M \,,$
$T\f N \,,$
$\f N \,,$
$T\f N \,,$
respectively, by
\bgt
(y^i) : \f M \to \Rn^{d_\f M} \,,
\qquad
(y^i, \dt y^i) : T\f M \to \Rn^{2 d_\f M} \,,
\\
(z^a) : \f N \to \Rn^{d_\f N} \,,
\qquad
(z^a, \dt z^a) : T\f N \to \Rn^{2 d_\f N} \,.
\end{gather*}

 Moreover, let us consider an F--smooth system
$(\f S, \eps)$
of smooth maps between the smooth manifolds
$\f M$
and
$\f N$
(see
Definition \ref{Definition: F--smooth systems of smooth maps}),
along with the set of F--smooth curves
$\C C$
of
$\f S$
defined in 
Theorem \ref{Theorem: F--smooth structure of F--smooth systems of smooth maps}.

\myskip

 Then, we define the F--smooth tangent space
$\E T\f S$
of the F--smooth space of parameters
$\f S \,,$
via equivalence classes of suitable smooth maps between smooth manifolds, in the following way.

\myskip

\bDf
\label{Definition: 1st order contact of curves}
\index{F--smooth systems!1st order contact of curves}
 We say that two F--smooth curves
\beq
\wha c_1 : \f I_1 \to \f S
\ssep{and} 
\wha c_2 : \f I_2 \to \f S
\eeq
have a \emp{1st order contact} in
$(\lam_1, \lam_2) \in \f I_1 \car \f I_2$
if the induced smooth maps
(see
Theorem \ref{Theorem: F--smooth structure of F--smooth systems of smooth maps})
\beq
\wha c_1^*(\eps) : \f I_1 \car \f M \to \f N
\ssep{and}
\wha c_2^*(\eps) : \f I_2 \car \f M \to \f N
\eeq
fulfill the following condition
\beq
\w(T_1 \wha c_1^*(\eps)\w)|_{(\lam_1, 1)} = 
\w(T_1 \wha c_2^*(\eps)\w)|_{(\lam_2, 1)}
: \f M \to T\f N \,,\ltag{*}
\eeq
i.e. if, in coordinates,
\bal
\wha c^*_1(\eps)^a|_{\lam_1} 
&= 
\wha c^*_2(\eps)^a|_{\lam_2} \,,
\\
\der_0 \wha c^*_1(\eps)^a|_{\lam_1} 
&= 
\der_0 \wha c^*_2(\eps)^a|_{\lam_2} \,,
\end{align*}
where
$\der_0$
denotes the partial derivative with respect to the parameter.

 Clearly, the above 1st order contact yields an equivalence relation
$\sim$
in the set of pairs
\beq
\w(\wha c : \f I_{\wha c} \to \f S, \;\; \lam \in \f I_{\wha c}\w) \,,
\eeq
where
$\wha c : \f I_{\wha c} \to \f S$
are F--smooth curves of
$\f S \,.$\END
\eDf

\bDf
\label{Definition: F--smooth tangent space of the F--smooth space of smooth maps}
\index{F--smooth systems!tangent vectors}
\index{F--smooth systems!F--smooth tangent space}
 We define a \emp{tangent vector}
$\E X_s$
of
$\f S \,,$
at
$s \in \f S \,,$
to be an equivalence class
(see the above
Definition \ref{Definition: 1st order contact of curves})
\beq
\E X_s \byd \w[(\wha c_s, \lam)\w]_\sim
\eeq
where 
$\wha c_s : I_{\wha c_s} \to \f S$
are F--smooth curves, such that
$\wha c_s(\lam) = s \,.$

 Then, we define:
 
 1) the \emp{F--smooth tangent space} of
$\f S \,,$
at
$s \in \f S \,,$
to be the set of tangent vectors of
$\f S \,,$
at
$s \,,$
\beq
\E T_s \f S \byd \wob \E X_s \wcb
\eeq

 2) the \emp{F--smooth tangent space} of
$\f S$
to be the disjoint union
\beq
\E T\f S \byd \bigsqcup_{s \in \f S}{\E T_s}\f S \,.\END
\eeq
\eDf

\myskip

 Thus, in virtue of the above 
Definition \ref{Definition: 1st order contact of curves}
and
Definition \ref{Definition: F--smooth tangent space of the F--smooth space of smooth maps}, 
the tangent vectors
$\E X_s \in \E T_s\f S$
can be represented through suitable smooth maps
$\Xi : \f M \to T\f N \,,$
as follows.

\myskip

\bTh
\label{Theorem: representation of tangent space of space of maps}
 Let
$(\f S, \eps)$
be an F--smooth system of smooth maps.
 
 Then, in virtue of
Definition \ref{Definition: 1st order contact of curves}
and
Definition \ref{Definition: F--smooth tangent space of the F--smooth space of smooth maps},
every tangent vector
\beq
\E X_s \byd [(\wha c, \lam)]_\sim \in \E T_s\f S \,,
\eeq
can be regarded as the smooth map
\beq
\Xi_s \byd \w(T_1(\wha c^*(\eps))\w)_{|(\lam, 1)} : \f M \to T\f N
\,,
\eeq
which projects on
$\eps_s \,,$
according to the following commutative diagram
\newdiagramgrid{2x2}
{1.2, 1.2, .7}
{.7, .7}
\bdg[grid=2x2]
\f M &\rTo^{\Xi_s} &T\f N
\\
\dTo^{\id_\f M} &&\dTo_{\tau_\f N}
\\
\f M &\rTo^{\eps_s} &\;\f N &.
\edg

 We have the coordinate expression
\beq
(z^a, \dt z^a) \com \Xi_s = (\eps_s^a, \; \Xi^a_s)
\,,
\ssep{with}
\Xi^a_s = \der_0 \w(c^*(\eps)^a\w)_{|(\lam, 1)} \in \map(\f M, \, \Rn) \,.
\eeq

 Thus, in this way, we obtain a natural injective map
\beq
\E r_s \, : \, \E X_s \in \E T_s\f S \, \mto \, \Xi_s : \f M \to T\f N \,,
\eeq
which yields a representation of the tangent vectors
$\E X_s \in \E T_s\f S \,,$
through suitable smooth maps
$\Xi_s : \f M \to T\f N \,,$
which project on
$\eps_s : \f M \to \f N \,.$\END
\eTh

\myskip

 The F--smooth tangent space
$\E T\f S$
can be also represented in another equivalent way through smooth maps between smooth manifolds, as follows.

\myskip

\bCr
\label{Corollary: equivalent representation of the tangent space of space of maps}
 The tangent vectors
\beq
\Xi_s \seq \E X_s \byd \w[(\wha c, \lam)\w]_\sim \in \E T_s\f S
\eeq
can be regarded as affine fibred morphisms over the smooth map
$\eps_s : \f M \to \f N$
\beq
\ol\Xi_s : T\f M \to T\f N \,,
\eeq
whose fibre derivative is equal to the smooth tangent prolongation of
$\eps_s : \f M \to \f N$
\beq
D\ol\Xi_s = T(\eps_s) : T\f M \to T\f N \,,
\eeq
where we have considered the natural identifications
\beq
VT\f M \seq T\f M \ucar{\f M} T\f M
\ssep{and}
VT\f N \seq T\f N \ucar{\f M} T\f N \,.
\eeq

 Their coordinate expressions are of the type
\beq
(z^a, \dt z^\R a) \com \ol\Xi_s = 
(\eps_s^a, \; \Xi^a_s + \der_i \eps_s^a \, \dt y^i) \,,
\ssep{with}
\Xi^a_s \in \map(\f M, \, \Rn) \,.
\eeq
\eCr

\bpf
 Each smooth map
$\Xi_s : \f M \to T\f N \,,$
which projects over the smooth map
$\eps_s : \f M \to \f N \,,$
yields the affine fibred morphism over
$\eps_s : \f M \to \f N$
\beq
\ol\Xi_s \byd \Xi_s \com \tau_\f M + T(\eps_s) \,,
\eeq
with coordinate expression
\beq
(z^a, \dt z^a) \com \ol\Xi_s = 
(\eps^a_s, \; \Xi^a_s + \der_i \eps^a_s \, \dt y^i) \,,
\ssep{where}
\Xi^a_s \in \map(\f M, \, \Rn) \,.
\eeq

 Conversely, each smooth affine fibred morphism over
$\eps_s : \f M \to \f N \,,$
\beq
\ol\Xi_s : T\f M \to T\f N \,,
\eeq
whose fibre derivative is
$D\ol\Xi_s = T(\eps_s) \,,$
yields the smooth map
\beq
\Xi_s \byd \ol\Xi_s \com 0_\f M \,,
\ssep{where}
0_\f M : \f M \to T\f M \,,
\eeq
which projects over
$\eps_s : \f M \to \f N \,.$

 Indeed, the above correspondence
\beq
\Xi_s \mto \ol\Xi_s \mto \Xi_s
\eeq
is a bijection.\QED
\epf

\bRm
\label{Remark: higher order representation of tangent vectors}
 In the above 
Corollary \ref{Corollary: equivalent representation of the tangent space of space of maps},
we have exploited the fact that, given
$s \in \f S \,,$
we know the smooth map
$\eps_s : \f M \to \f N \,,$
along with its smooth tangent prolongation
$T(\eps_s) : T\f M \to T\f N \,.$\END
\eRm

\bNt
 By the way, by means of a reasoning analogous to that of the above
Corollary \ref{Corollary: equivalent representation of the tangent space of space of maps},
a representation of
$\E T\f S$
can be obtained by considering smooth maps of the type
\beq
\ol\Xi^k_s : T^k\f M \to T^k\f N \,,
\eeq
for any order
$k$
of tangent prolongation.\END
\eNt

\myskip
 
 In particular, the above results can be applied to the F--smooth system of compact support global smooth maps
(see
Example \ref{Example: F--smooth system of all smooth maps with compact support}).
 
\myskip
 
 Next, we introduce two natural maps associated with
$\E T\f S$
and show the natural vector structure of the fibres of
$\E T\f S \,.$

\myskip

\bLm
\label{Lemma: natural maps tau and T eps}
 We have the \emp{natural projection}
\beq
\tau_\f S : \E T\f S \to \f S : \Xi_s \mto s
\eeq
and the \emp{natural evaluation map}
\beq
\E T_1\eps : \E T\f S \car \f M \to T\f N :
(\Xi_s, m) \mto \Xi_s(m) \,.\END
\eeq
\eLm

\bPr
\label{Proposition: TS is a vector fibred space}
 The fibres of 
$\tau_\f S : \E T\f S \to \f S$ 
are naturally equipped with a vector structure induced by the maps
\bal
(\Xi_s + \ac\Xi_s) 
&: 
\f M \to T\f N : m \mto \Xi_s (m) + \ac\Xi_s (m) \in 
T_{\eps(s,m)}\f N \,,
\\
(k \, \Xi_s)
&:
\f M \to T\f N : m \mto k \, \Xi_s (m) \in 
T_{\eps(s,m)}\f N \,.\END
\end{align*}
\ePr

\myskip

 We have not assumed a priori any kind of smoothness on the set
$\E T\f S \,,$
but we can uniquely recover its natural F--smooth structure by the following
Theorem \ref{Theorem: the tangent space of the F--smooth system of smooth maps is F--smooth}.

 Indeed, the specification ``\emp{F--smooth}" tangent space, that we have anticipated in Definition \ref{Definition: F--smooth tangent space of the F--smooth space of smooth maps},
is justified by this Theorem.

\bTh
\label{Theorem: the tangent space of the F--smooth system of smooth maps is F--smooth}
  Let us define the set
\beq
\E T\C C \byd \wob \whaE c : \f I_{\whaE c} \to \E T\f S \wcb
\eeq
consisting of all curves of 
$\E T\f S$
such that 

 1) the base curve
\beq
\wha c_\f S \byd \tau_\f S \com \whaE c : \f I_{\whaE c} \to \f S
\eeq
be F--smooth,

 2) the induced map between smooth manifolds
\beq
\Xi_{\whaE c} : \f I_{\whaE c} \car \f M \to T\f N :
(\lam, m) \mto \Xi_{\wha c_\f S(\lam)}(m)
\eeq
be smooth.

 Then, the pair
$(\E T\f S, \E T\C C)$
turns out to be an F--smooth space and the maps
\beq
\tau_\f S : \E T\f S \to \f S
\ssep{and}
T_1\eps : \E T\f S \car \f M \to T\f N
\eeq
turn out to be F--smooth.
\eTh

\bpf 
 a) First of all, let us check that the set
$\E T\C C$
fulfills the two requirements of 
Definition \ref{Definition: F--smooth spaces}.

\myskip

 1) For each
$\E X \in \E T\f S \,,$
the constant curve
$\whaE c : \Rn \to \E T\f S : \lam \mto \E X$
turns out to be an element of
$\E T\C C$
passing through
$\E X \,.$

\myskip

 2) Let us consider a curve
$\whaE c \in \E T\C C$
and a smooth curve
$\gam : \f I_\gam \to \f I_{\whaE c} \,.$

 Then, the induced map
\beq
\Xi_{{\whaE c}} \com \gam_\f M : 
\f I_\gam \car \f M \to T\f N :
(\lam, m) \mto \eps_{(\whaE c \com \gam)(\lam)} (m)
\eeq
is the composition of two smooth maps
\newdiagramgrid{1}
{1.5, 3, 1.5, 1.5, .7}
{1}
\bdg[grid=1]
\f I_\gam \car \f M &\rTo^{\gam_\f M \byd (\gam, \id_\f M)} 
&\f I_{\whaE c} \car \f M
&\rTo^{\Xi_{\whaE c}} &T\f N &,
\edg
hence it turns out to be smooth.

\myskip

 b) In order to prove that
$\tau_\f S : \E T\f S \to \f S$
be F--smooth, we have to show that, for each F--smooth curve
$\whaE c : \f I_{\whaE c} \to \E T\f S \,,$
the induced curve
$\wha c_\f S \byd \tau_\f S \com \whaE c : \f I_{\whaE c} \to \f S$
be F--smooth.

 In fact, let us consider any F--smooth curve and the associated smooth map
\beq
\whaE c : \f I_{\whaE c} \to \E T\f S
\ssep{and}
\Xi_{\whaE c} : \f I_{\whaE c} \car \f M \to T\f N \,.
\eeq

 Indeed, the induced curve
\beq
\wha c_\f S \byd \tau_\f S \com \whaE c : \f I_{\whaE c} \to \f S
\eeq
is characterised by the composition of smooth maps
\beq
\tau_\f N \com \Xi_{\whaE c} : \f I_{\whaE c} \car \f M \to \f N \,,
\eeq
hence it turns out to be smooth.

\myskip

 c) In order to prove that
$\E T_1\eps : \E T\f S \car \f M \to T\f N$
be F--smooth, we have to show that, for each F--smooth curve
$\whaE c = 
(\whaE c_{\E T\f S}, c_\f M) : 
\f I_{\whaE c} \to \E T\f S \car \f M \,,$
the induced curve
$\E T_1\eps \com \whaE c : \f I_\whaE c \to T\f N$
be F--smooth.

 In fact, let us consider any F--smooth curve and the associated smooth map
\beq
\whaE c_{\E T\f S} : \f I_\whaE c \to \E T\f S
\ssep{and}
\Xi_c : \f I_\whaE c \car \f M \to T\f N \,.
\eeq

 Then, the curve
\beq
\E T_1\eps \com (\whaE c_{T\f S}, c_\f M) : \f I_\whaE c \to T\f N :
\lam \mto \E T_1\eps \w(\whaE c_{\E T\f S}(\lam), \, c_\f M(\lam)\w)
\eeq
which is a composition of smooth maps according to the following diagram
\newdiagramgrid{2x2}
{2, 2, .6}
{.7, .7}
\bdg[grid=2x2]
\f I_\whaE c &\rTo^{(\id, c_\f M)} &\f I_\whaE c \car \f M
\\
\uTo^{\id_{I_\whaE c}} &&\dTo_{\Xi_\whaE c}
\\
\f I_\whaE c &\rTo^{\E T_1\eps \com (\whaE c_{\E T\f S}, c_\f M)} 
&T\f N &,
\edg
turns out to be smooth hence F--smooth
(see
Theorem \ref{Theorem: F--smooth structure of F--smooth systems of smooth maps}).\QED
\epf

\myskip

 We say that
$\tau_\f S : \E T\f S \to \f S$
is an \emp{F--smooth fibred set}.

\myskip

 Eventually, we discuss the tangent prolongation of F--smooth curves of
$\f S \,.$

\myskip

\bDf
\label{Definition: tangent prolongation of an F--smooth curve}
 Let
$\wha c : \f I_\wha c \to \f S$
be an F--smooth curve.
 Then, we define its \emp{tangent prolongation} to be the curve
\beq
d \, \wha c : \f I_\wha c \to \E T\f S : 
\lam \mto d \wha c(\lam) \byd [(\wha c, \lam)]_\sim
\eeq
defined in
Definition \ref{Definition: F--smooth tangent space of the F--smooth space of smooth maps}.\END
\eDf

\bPr
\label{Proposition: the tangent prolongation of an F--smooth curve is F--smooth}
 If
$\wha c : \f I_\wha c \to \f S$
is an F--smooth curve, then its tangent prolongation
$d \, \wha c : \f I_\wha c \to \E T\f S$
turns out to be F--smooth.
\ePr

\bpf
 Let
$\gam : \f I_\gam \to \f I_\wha c$
be a smooth curve.
 Then, the induced map
\beq
\w(d \, (\wha c \com \gam)\w)^*(T_1\eps) : 
\f I_\gam \car \f M \to T\f N
\eeq
turns out to be the composition of smooth maps.
 Hence,
$d \, \wha c$
is F--smooth.\QED
\epf
%--------------------------------------------------------------------%
%\newpage
\markboth{\rm Chapter \thechapter. Systems of maps}{\rm \thesection. F--smooth vs smooth systems of maps}
\section{F--smooth vs smooth systems of maps}
\label{F--smooth vs smooth systems of maps}
\markboth{\rm Chapter \thechapter. Systems of maps}{\rm \thesection. F--smooth vs smooth systems of maps}
%--------------------------------------------------------------------%
\bsm
 In the above 
Sections \ref{Smooth tangent prolongation of (S, eps)}
and
\ref{F--smooth tangent prolongation of (S, eps)}
we have studied the smooth tangent space
$T\f S$
and the F--smooth tangent space
$\E T\f S$
with reference to a smooth system and to an F--smooth system
$(\f S, \eps)$
of smooth maps between two smooth manifolds, respectively.
 
 Clearly, a smooth system is a particular case of an F--smooth system, because a smooth system can be regarded as an F--smooth system along with the additional assumptions on the smoothness of
$\f S$
and
$\eps \,.$

 Then, the need of a comparison between the smooth approach to
$T\f S$
and the F--smooth approach to
$\E T\f S$
arises naturally, having in mind
Theorem \ref{Theorem: relation between smoothness and F--smoothness}.

 Actually, by regarding a smooth system
$(\f S, \eps)$
of smooth maps as a particular F--smooth system of smooth maps, we show a natural bijection
\beq
\imath : T\f S \to \E T\f S : X \mto \E X \,,
\eeq
which, in terms of the representation of
$\E X$
via the smooth map
\beq
\Xi \byd \E r(\E X) : \f M \to T\f N
\eeq
reads as
\beq
\Xi = (T_1\eps)_X \,.
\eeq
\esm

\bNt
 Summing up previous results, let us compare two types of systems of smooth maps between smooth manifolds.

\myskip

 1) In the case of an F--smooth system
$(\f S, \eps)$
of smooth maps
$f : \f M \to \f N \,,$
we make no assumptions on any kind of smoothness of the set
$\f S$
(see
Definition \ref{Definition: F--smooth systems of smooth maps}).

 Hence, we cannot avail of a smooth structure of
$\f S$
and, in order to achieve the F--smooth tangent space
$\E T\f S \,,$
we need to follow an indirect abstract procedure.

 In fact, we have defined the tangent vectors
$\E X_s \in \E T_s\f S \,,$
with
$s \in \f S \,,$
as equivalence classes
$\E X_s \byd [(\wha c, \lam)]_\sim$
of pairs
$(\wha c, \lam) \,,$
where
$\wha c : \f I_\wha c \to \f S$
are F--smooth curves and
$\lam \in \f I_\wha c \,,$
which have a 1st order contact in
$s$
(see
Definition \ref{Definition: F--smooth tangent space of the F--smooth space of smooth maps}).

 Here, the 1st order contact is defined through the smooth map
(see
Definition \ref{Definition: 1st order contact of curves})
\beq
\w(T_1 \wha c^*(\eps)\w)|_{(\lam, 1)} : \f M \to T\f N \,.
\eeq

 Then, we obtain the representation of tangent vectors
$\E X_s \byd [(\wha c, \lam)]_\sim\in \E T_s\f S \,,$
through the smooth maps
\beq
\Xi_s \byd \w(T_1(\wha c^*_s(\eps))\w)|_{(\lam, 1)} : \f M \to T\f N \,,
\eeq
which project on
$\eps_s : \f M \to \f N$
(see
Theorem \ref{Theorem: representation of tangent space of space of maps}).\END

\myskip

 2) In the case of a finite dimensional smooth system
$(\f S, \eps)$
of smooth maps
$f : \f M \to \f N \,,$
we assume that
$\f S$
is a smooth manifold and
$\eps$
a smooth map.

 Let us denote the typical smooth chart of
$\f S$
by
$(w^A) \,.$

 Hence, we can avail of the assumed smooth structure of
$\f S$
to achieve a direct approach of the tangent space
$T\f S \,.$ 

 In fact, according to a standard definition in Differential Geometry, we define the tangent vectors
$X_s \in T_s\f S \,,$
with
$s \in \f S \,,$
as equivalence classes
$X_s \byd [(\wha c, \lam)]_\sim$
of pairs
$(\wha c, \lam) \,,$
where
$\wha c : \f I_\wha c \to \f S$
are smooth curves and
$\lam \in \f I_\wha c \,,$
which have a 1st order contact in
$s \,.$

 Here, the 1st order contact is defined directly through the smooth structure of 
$\f S$
(without the need to consider the smooth map
$\wha c^*(\eps) : \f I_\wha c \car \f M \to \f N$).

 Then, the smooth tangent map
(see
Section \ref{F--smooth tangent prolongation of (S, eps)})
\beq
T_1\eps : T\f S \car \f M \to T\f N \,,
\eeq
with coordinate expression
\beq
z^a \com T_1\eps = \eps^a
\ssep{and}
\dt z^a \com T_1\eps = 
\der_A\eps^a \, \dt w^A \,,
\eeq
provides, for each
$X_s \in T_s\f S \,,$ 
the smooth map
\beq
(T_1\eps)_{X_s} : \f M \to T\f N \,.\END
\eeq
\eNt

\bPr
\label{Proposition: bijection TS to TS}
 Let us consider a smooth system
$(\f S, \eps)$
of smooth maps
$f : \f M \to \f N$
and regard it as a particular F--smooth system.
 
 Then, in virtue of the standard definition of the smooth tangent space
$T\f S$
and of the definition of the F--smooth tangent space
$\E T\f S$
(see
Definition \ref{Definition: F--smooth tangent space of the F--smooth space of smooth maps}),
we obtain a natural F--smooth map
\beq
\imath : T _s\f S \to \E T_s\f S : X_s \mto \E X_s \,,
\ssep{for each}
s \in \f S \,.
\eeq

 Indeed, in terms of the smooth representation of
$\E X_s \in \E T_s\f S \,,$
the above map turns out to be given by the smooth map
(see
Theorem \ref{Theorem: representation of tangent space of space of maps})
\beq
\Xi_s \byd \E r (\E X_s) = (T_1\eps)_{|X_s} : \f M \to T\f N \,,
\ssep{for each}
s \in \f S \,,
\eeq
i.e., in coordinates,
\beq
z^a \com \Xi_s = \eps^a_s
\ssep{and}
\dt z^a \com \Xi_s = (\der_A\eps^a)_s \, X^A_s \,.
\eeq
\ePr

\bpf
 Each smooth curve
$\wha c : \f I_\wha c \to \f S$
makes the map
\beq
\wha c^*(\eps) : \f I_\wha c \car \f M \to \f N
\eeq
smooth, hence
$\wha c : \f I_\wha c \to \f S$ 
turns out to be also F--smooth
(see
Theorem \ref{Theorem: F--smooth structure of F--smooth systems of smooth maps}).

 Moreover, we can easily see that, if two smooth pairs
$(\wha c_1, \lam_1)$
and
$(\wha c_2, \lam_2)$
have a 1st order contact in
$s \in \f S \,,$
in the sense of smooth manifolds, then they have also a 1st order contact in the sense of F--smooth spaces
(see
Definition \ref{Definition: 1st order contact of curves}).\QED
\epf

\myskip

 In the above 
Proposition \ref{Proposition: bijection TS to TS}, 
we have assumed a given smooth structure of
$\f S \,.$
 Now, we discuss what happens if we change this smooth structure.
 
\myskip

\bRm
\label{Remark: infinitely many smooth structures}
 We have already observed
(see
Example \ref{Example: several smooth structures underlying an F--smooth space})
that, given a smooth system
$(\f S, \eps)$
of smooth maps, there might exist infinitely many possible smooth structures of
$\f S$
which make the map
$\eps : \f S \car \f M \to \f N$
smooth, hence, which yield the same family of selected smooth maps
$f : \f M \to \f N \,.$

 Indeed, different smooth structures of
$\f S$
yield different F--smooth fibred morphisms
\beq
\imath : T \f S \to \E T\f S : X_s \mto \E X_s \,.
\eeq

 Therefore, the F--smooth tangent space
$\E T_s\f S$
contains the images
$\imath (X_s)$
of the tangent vectors
$X_s \in T\f S \,,$
for all possible smooth structures of
$\f S \,.$\END
\eRm

\bRm
\label{Remark: non injectivity of the map iota}
 In general, the map
$\imath : T \f S \to \E T\f S$
is not injective.
 We prove this fact by means of a counter--example.
 
 Let us consider the smooth manifolds 
$\f S \byd \Rn \,,$
$\f M \byd \Rn$
and
$\f N \byd \Rn \,,$ 
along with their natural smooth charts 
$w : \f S \to \Rn \,,$
$y : \f M \to \Rn$ 
and 
$z : \f N \to \Rn \,.$

 Moreover, let us consider the non--injective smooth system of smooth maps given by the evaluation map
$\eps : \f S \car \f M \to \f N \,,$ 
with coordinate expression
\beq
z \com \eps = w^2 \, y \,.
\eeq

 Let us consider the element 
$s = 1 \in \f S$ 
and the two smooth curves
\beq
\wha c_1 : \Rn \to \f S : \lam \mto \lam
\ssep{and}
\wha c_2 : \Rn \to \f S : \lam \mto - \lam \,.
\eeq
 
 Then, we obtain
\beq
\wha c_1 (1) = 1 = \wha c_2(-1)
\eeq
and two different tangent vectors
\beq
X_1 = T\wha c_1(1,0) = (1,1) \in T_1\f S
\ssep{and}
X_2 = T\wha c_2(-1,0) = (1,-1) \in T_1\f S \,.
\eeq

 On the other hand, the pairs 
$(\wha c_1,1)$ 
and 
$(\wha c_2,-1)$ 
are equivalent in the sense of
Definition \ref{Definition: 1st order contact of curves}.

 In fact, we have
\beq
w \com \w(\wha c^*_1(\eps)\w) (1) = 
y = 
w \com \w(\wha c^*_2(\eps)\w) (-1)
\eeq
and 
\beq
\w(\der_0 (\wha c^*_1(\eps))\w)|_1 = 
2 \, y =
\w(\der_0 (\wha c^*_2(\eps))\w)|_{-1} \,.
\eeq

 Then, we have
\beq
\iota(X_1) = \iota(X_2) \,.\END
\eeq
\eRm
%--------------------------------------------------------------------%
\chapter{Systems of sections}
\label{Systems of sections}
%--------------------------------------------------------------------%
\bsm
 First, we discuss the \emp{smooth systems} 
$(\f S, \zet, \eps)$
of \emp{smooth sections}
$\phi : \f F \to \f G$
of a smooth double fibred manifold
$\f G \to \f F \to \f B \,.$
 Here, the ``fibred space of parameters"
$\zet : \f S \to \f B$
is a \emp{smooth fibred manifold} and the ``evaluation map"
$\eps : \f S \ucar{\f B} \f F \to \f G$
a \emp{smooth fibred morphism}
over
$\f B \,.$

\myskip

 Then, we discuss the \emp{F--smooth systems}
$(\f S, \zet, \eps)$
of \emp{fibrewisely smooth sections}
$\phi : \f F \to \f G$
of a smooth double fibred manifold
$\f G \to \f F \to \f B \,.$
 Here, we have a weaker smoothness requirement, as the ``fibred space of parameters"
$\zet : \f S \to \f B$
turns out to be an \emp{F--smooth fibred} space and the ``evaluation map"
$\eps : \f S \ucar{\f B} \f F \to \f G$
an \emp{F--smooth fibred morphism}
over 
$\f B \,.$

\myskip

 The above notions are intended as an introduction to the particular cases of systems of connections discussed in the next 
Chapter \S\ref{Systems of connections}.

 The reader can find further discussions concerning the present subject in
\cite{Cab87,CabKol95,CabKol98,DodMod86,Gar72,JanMod02c,
Kol87,KolMicSlo93,KolMod98,ManMod83d,MarMod91,Mod91}.
\esm
%--------------------------------------------------------------------%
%\newpage
\markboth{\rm Chapter \thechapter. Systems of sections}{\rm \thesection. Smooth systems of smooth sections}
\section{Smooth systems of smooth sections}
\label{Section: Smooth systems of smooth sections}
\markboth{\rm Chapter \thechapter. Systems of sections}{\rm \thesection. Smooth systems of smooth sections}
%--------------------------------------------------------------------%
\bsm
 We discuss the notion of \emp{smooth systems of smooth sections} of a smooth double fibred manifold.
\esm
%--------------------------------------------------------------------%
\subsection{Smooth systems of smooth sections}
\label{Smooth systems of smooth sections}
%--------------------------------------------------------------------%
\bsm
 First of all, given a smooth double fibred manifold
$\f G \oset{q}{\to} \f F \oset{p}{\to} \f B \,,$
we define the \emp{``tubelike" smooth sections}
$\phi : \f F \to \f G \,.$

 Next, we discuss the notion of \emp{``smooth system of smooth sections"} of the smooth double fibred manifold.
 
 In simple words, a ``smooth system of smooth sections" is defined to be a selected family
$\{\phi\}$
of smooth sections
$\phi : \f F \to \f G$
of the smooth fibred manifold
$q : \f G \to \f F \,,$
which is parametrised by the smooth sections
$\sig : \f B \to \f S$
of a smooth fibred manifold
$\zet : \f S \to \f B \,.$
 
\myskip
 
 The notion of smooth system of smooth sections will be used as an introduction to the concept of ``smooth system of smooth connections", which is developed in the forthcoming 
Section \ref{Smooth systems of smooth connections}.

\myskip

 Later, in the next
Section \ref{F--smooth systems of smooth sections}, 
we shall revisit the notion of system of sections in a larger context, detached from the hypothesis of smoothness and finite dimension and approached by means of the concept of F--smoothness
(see
Section \ref{F--smooth spaces}).
\esm

 Let us consider a smooth double fibred manifold
$\f G \oset{q}{\to} \f F \oset{p}{\to} \f B \,.$

\myskip

 We start by introducing the preliminary notions of ``tubelike subset" and ``tubelike section".

\myskip

\bDf
\label{Definition: tubelike subsets}
\index{smooth systems!tubelike subset}
 We define the \emp{tubelike}
subsets of
$\f F$
and of
$\f G$ 
to be the open subsets of the type
\beq
p^{-1} (\f U) \sub \f F
\ssep{and}
(p \com q)^{-1} (\f U) \sub \f G \,,
\eeq
where
$\f U \sub \f B$
is an open subset, according to the following commutative diagrams
\newdiagramgrid{2x2-2x2}
{1.5, 1.5, .4, 2.5, 1.5, 1.5, .6}
{.7, .7}
\bdg[grid=2x2-2x2]
p^{-1}(\f U) &\rInto^{\sub} &\f F 
&&(p \com q)^{-1}(\f U) &\rInto^{\sub} &\f G
\\
\dTo^{p} &&\dTo_{p}
&&\dTo^{p \com q} &&\dTo_{p \com q}
\\
\f U &\rInto^{\sub} &\f B &
&\f U &\rInto^{\sub} &\f B &.\END
\edg
\eDf

\myskip

\bNt
\label{Note: tubelike topology}
\index{smooth systems!tubelike topology}
 The tubelike open subsets
$p^{-1} (\f U) \sub \f F$
and
$(p \com q)^{-1} (\f U) \sub \f G$
yield a topology on
$\f F$
and
$\f G \,.$

 In the following, we shall usually refer to this topology.\END
\eNt

\bDf
\label{Definition: tubelike sections}
\index{smooth systems!tubelike section}
 A smooth section
$\phi \in \sec(\f F, \f G)$
is said to be \emp{tubelike}
if it is globally defined on a tubelike open subset
$p^{-1} (\f U) \sub \f F \,.$

 We denote the subsheaf (with respect to the tubelike topology) of \emp{tubelike smooth sections}
$\phi \in \sec(\f F, \f G)$
by
\beq
\tub(\f F, \f G) \sub \sec(\f F, \f G) \,.\END
\eeq
\eDf

\myskip

 Then, we introduce the notion of ``smooth system of smooth sections".
 
\myskip

\bDf
\label{Definition: systems of sections}
\index{smooth systems!system of sections} 
\index{smooth systems!injective systems of sections}
 We define a \emp{smooth system of smooth sections}
of the smooth fibred manifold
$q : \f G \to \f F$
to be a 3--plet
$(\f S, \zet, \eps) \,,$
where 

1) $\zet : \f S \to \f B$ 
is a smooth fibred manifold,

2) $\eps : \f S \ucar{\f B} \f F \to \f G$
is a smooth fibred map over 
$\f F \,,$
according to the following commutative diagram
\newdiagramgrid{2x2}
{1.2, 1.2, .6}
{.7, .7}
\bdg[grid=2x2]
\f S\ucar{\f B} \f F &\rTo^{\eps} &\f G
\\
\dTo^{\pro_2} &&\dTo_{q}
\\
\f F &\rTo^{\id_\f F} &\f F &.
\edg

 We call
$\eps$
the \emp{evaluation map} of the system.

 Thus, the evaluation map
$\eps$
yields the sheaf morphism
\beq
\eps_\f S : \sec(\f B, \f S) \to \tub(\f F, \f G) : 
\sig \mto \br\sig \,,
\eeq
where, for each
$\sig \in \sec(\f B, \f S) \,,$ 
the tubelike smooth section
$\br\sig$
is defined by
\beq
\br\sig : \f F \to \f G : f_b \mto \eps (\sig(b), f_b) \,,
\ssep{for each}
b \in \f B \,.
\eeq

 Therefore, the map
$\eps_\f S : \sec(\f B, \f S) \to \tub(\f F, \f G)$
provides a \emp{selection} of the tubelike smooth sections
$\phi : \f F \to \f G \,,$
given by the subset
\beq
\tub_\f S(\f F, \f G) \byd \eps_\f S\w(\sec(\f B, \f S)\w) \sub \tub(\f F, \f G) \,.
\eeq

\myskip

 The smooth system of smooth sections
$(\f S, \zet, \eps)$
is said to be \emp{injective} if the map
$\eps_\f S : \sec(\f B, \f S) \to \tub(\f F, \f G)$
is injective, i.e. if, for each
$\sig \,, \ac\sig \in \sec(\f B, \f S) \,,$
\beq
\br\sig \eqv \eps_\f S (\sig) = \brac\sig \eqv \eps_\f S(\ac\sig)
\Rarr
\sig = \ac\sig \,.
\eeq

 If the system is injective, then we obtain the bijection
\beq
\eps_\f S : 
\sec(\f B, \f S) \to \tub_\f S(\f F, \f G) :
\sig \mto \br\sig \,,
\eeq
whose inverse is denoted by
\beq
(\eps_\f S)^{-1} : 
\tub_\f S(\f F, \f G) \to \sec(\f B, \f S) :
\phi \mto \wha\phi \,.\END
\eeq
\eDf

\myskip

 Indeed, we are essentially interested in injective smooth systems of smooth sections.

\myskip

 We shall denote the local fibred smooth charts of
$\f B \,,$
$\f S \,,$
$\f F \,,$
$\f G \,,$
respectively, by
\bgt
(x^\lam) : \f B \to \Rn^{d_\f B} \,,
\qquad
(x^\lam, w^A) : \f S \to \Rn^{d_\f S} \,,
\\
(x^\lam, y^i) : \f F \to \Rn^{d_\f F} \,,
\qquad
(x^\lam, y^i, z^a) : \f G \to \Rn^{d_\f G} \,.
\end{gather*}

\myskip

 We discuss the following elementary examples of smooth systems of smooth sections.

\myskip

\bEx
\label{Example: system of linear sections}
\index{smooth systems!system of linear sections}
 If
$p : \f F \to \f B$
and
$p \com q : \f G \to \f B$
are vector bundles, then the family of tubelike sections
$\br\sig : \f F \to \f G \,,$ 
which are linear fibred morphisms over
$\f B \,,$
yields the injective smooth system of smooth sections
\beq
\f S \byd \lin_\f B(\f F, \f G) \,.
\eeq
 
 With reference to a linear fibred chart of the smooth double fibred manifold, the coordinate expression of the tubelike sections
$\br\sig : \f F \to \f G$
of this system is of the type
\beq
x^\lam \com \br\sig = x^\lam \,,
\qquad
y^i \com \br\sig = y^i \,,
\qquad
z^a \com \br\sig = K^a_i \, y^i \,,
\ssep{with}
K^a_i \in \map(\f B, \Rn) \,.
\eeq

 Hence, the above fibred chart of the smooth double fibred manifold yields a distinguished fibred chart
$(x^\lam, w^a_i)$
of
$\f S$
and the coordinate expression of
$\eps$
becomes
\beq
\eps^a = w^a_j \, y^j \,.\END
\eeq
\eEx

\bEx
\label{Example: system of affine sections}
\index{smooth systems!system of affine sections}
 If
$p : \f F \to \f B$
and
$p \com q : \f G \to \f B$
are affine bundles, then the family of tubelike sections
$\br\sig : \f F \to \f G \,,$ 
which are affine fibred morphisms over
$\f B \,,$
yields the injective smooth system of smooth sections
\beq
\f S \byd \aff_\f B(\f F, \f G) \,.
\eeq
 
 With reference to an affine fibred chart of the smooth double fibred manifold, the coordinate expression of the tubelike sections
$\br\sig : \f F \to \f G$
of this system is of the type
\beq
x^\lam \com \br\sig = x^\lam \,,
\quad
y^i \com \br\sig = y^i \,,
\quad
z^a \com \br\sig = K^a_i \, y^i + K^a \,,
\sep{with}
K^a_i, \, K^a \in \map(\f B, \Rn) \,.
\eeq

 Hence, the above fibred chart of the smooth double fibred manifold yields the distinguished fibred chart
$(x^\lam, w^a_i, w^a)$
of
$\f S$
and the coordinate expression of
$\eps$
becomes
\beq
\eps^a = w^a_j \, y^j + w^a \,.\END
\eeq
\eEx

\bEx
\label{Example: system of polynomial sections}
\index{smooth systems!system of polynomial sections}
 If
$p : \f F \to \f B$
and
$p \com q : \f G \to \f B$
are affine bundles, then analogously to the above 
Example \ref{Example: system of affine sections}, 
we can define 

 - the injective smooth system of polynomial sections of degree
$r \,,$
with
$0 \leq r \,,$

 - the injective smooth system of polynomial sections of any degree 
$0 \leq r \leq k \,,$
where
$k$
is a given positive integer.\END
\eEx

\myskip

 All examples above deal with finite dimensional smooth systems of smooth sections, as it is implicitly requested by
Definition \ref{Definition: systems of sections}.

 On the other hand, we can easily extend the concept of smooth system of smooth sections, by considering an infinite dimensional smooth system, which is the direct limit of finite dimensional smooth systems, according to the following
Example \ref{Example: infinite dimensional system of polynomial sections of any degree}.

\myskip

\bEx
\label{Example: infinite dimensional system of polynomial sections of any degree}
\index{smooth systems!system of polynomial sections}
 If
$p : \f F \to \f B$
and
$p \com q : \f G \to \f B$
are affine bundles, then we obtain the smooth system of polynomial sections by considering the family of polynomial sections
$s : \f F \to \f G$
of any degree
$r \,,$
with
$0 \leq r \leq \infty \,.$

 Later, we shall see that such a smooth system has a natural F--smooth structure
(see, later,
Definition \ref{Definition: F--smooth spaces})\END
\eEx
%--------------------------------------------------------------------%
\subsection{Smooth structure of $(\f S,\zet,\eps)$}
\label{Smooth structure of (S,zet,eps)}
%--------------------------------------------------------------------%
\bsm
 In our definition of ``smooth system of smooth sections"
(see
Definition \ref{Definition: systems of sections})
we have required a priori that the fibred set of parameters
$\zet : \f S \to \f B$
be smooth and that the evaluation map
$\eps$
be smooth as well.

 On the other hand, we might ask whether the fact that the manifold
$\f B \,,$
the fibred manifolds
$p : \f F \to \f B$
and
$q : \f G \to \f F$
are smooth and that the sections
$\phi : \f F \to \f G$
selected by the system are smooth allows us to recover uniquely the smooth structure of
$\zet : \f S \to \f B \,.$

 The answer to the above question is negative.
 Here we do not fully address this problem, which is too far from the true scope of the present report.
  But, we present a simple example
(see
Example \ref{Example: different smooth structures of a system of sections}),
where we show that, if
$\f S$
admits a finite dimensional smooth structure compatible with
$\eps \,,$
then this structure needs not to be unique.

 Moreover, we observe that if we do not assume a priori a finite dimensional smooth structure on
$\zet : \f S \to \f B \,,$
then it might be that no finite dimensional smooth structure at all could be recovered on
$\zet : \f S \to \f B \,.$
 To prove this, just consider the system
$(\f S, \zet, \eps)$
of \emp{all} smooth sections
$\phi : \f F \to \f G$
(see, later,
Section \ref{F--smooth systems of smooth sections}).

\myskip

 The above question might arise also in comparison with a result which will be achieved later, in the next Section, in the context of ``F--smooth systems of smooth sections", where we do not assume a priori any smooth structure of the fibred set
$\zet : \f S \to \f B \,,$
but we recover uniquely its F--smooth structure
(see
Definition \ref{Definition: F--smooth systems of smooth sections}
and
Theorem \ref{Theorem: F--smooth structure of F--smooth systems of smooth sections}).
\esm

 Analogously to
Example \ref{Example: different smooth structures of a system of maps},
we can exhibit a system of smooth sections
$(\f S, \zet, \eps) \,,$
where
$\f S$
is equipped with different smooth structures.

\myskip

\bEx
\label{Example: different smooth structures of a system of sections}
 Let us consider the following system
$(\f S, \zet, \eps)$
of smooth sections:
 
 1) we consider the manifolds
\bgt
\f B \byd \Rn \,,
\qquad
\f F \byd \Rn \car \Rn \,,
\qquad
\f G \byd \Rn \car \Rn \car \Rn \,,
\\
\f S \byd \Rn \car \Rn \,,
\end{gather*}
equipped with their ``natural" smooth structures, and denote their ``natural" charts by
\bgt
x : \f B \to \Rn \,,
\qquad
(x,y) : \f F \to \Rn \car \Rn \,,
\qquad
(x,y, z) : \f G \to \Rn \car \Rn \car \Rn \,,
\\
(x, w) : \f S \to \Rn \car \Rn \,.
\end{gather*}

 Moreover, let us consider the smooth evaluation map
\beq
\eps : \f S \ucar{\f B} \f F \to \f G \,,
\eeq
defined by the coordinate expression
\beq
(x, y; \, z) \com \eps = (x, y; \, w \, y) \,.\ltag{1}
\eeq
 
 Besides the above ``natural" smooth structure of
$\f S \,,$
we can consider several further smooth structures of
$\f S \,,$
which make the \emp{same} evaluation map
$\eps$
smooth, hence which essentially yield the ``same" system of sections.

 For instance, let us consider the ``exotic" smooth structure of
$\f S$
induced by the fibred bijection
\beq
(x, \ac w) \byd (x, w^3) : \f S \to \Rn \car \Rn \,.
\eeq

 Then, the evaluation map
$\eps$
reads, in the above exotic fibred chart as
\beq
(x, y; \, z) \com \eps = (x, y; \, w^3 \, y) \,.\ltag{2}
\eeq

 Clearly, the equalities
$1)$
and
$2)$
define the same evaluation map
$\eps \,.$
 Moreover, 
$\eps$
turns out to be smooth with respect to both smooth structures of
$\f S \,.$

 By the way, the system
$(\f S, \eps)$
turns out to be injective in both cases.\END
\eEx

\myskip

 We might even ask whether we can characterise the ``natural" smooth structure of
$\f S$
via a suitable feature of the system; this question might arise later if we compare ``smooth systems of smooth maps" and ``F--smooth systems of smooth maps" 
(see, later, in the next
Section \ref{F--smooth systems of smooth maps}).

 Here, we skip a general answer to this question, which is too far from the true scope of the present report.
 
 On the other hand, we observe, as a hint, that the set of curves of
$\f S$
which are smooth with respect to its natural smooth structure is smaller than the set of curves of
$\f S$
which are smooth with respect to its exotic smooth structure.
 The converse occurs for the sets of smooth functions.
%--------------------------------------------------------------------%
\subsection{Smooth lifted fibred manifold}
\label{Smooth lifted fibred manifold}
%--------------------------------------------------------------------%
\bsm
 Given a smooth double fibred manifold
$\f G \oset{q}{\to} \f F \oset{p}{\to} \f B$
and a smooth system of smooth sections
$(\f S, \zet, \eps) \,,$
it is useful to define the ``\emp{lifted smooth fibred manifold}"
\beq
p\Upa : \f F\Upa \byd \f S \ucar{\f B} \f F \to \f S
\eeq
of the fibred manifold
$p : \f F \to \f B \,.$ 
\esm

 Thus, let us consider a smooth double fibred manifold
$\f G \oset{q}{\to} \f F \oset{p}{\to} \f B$
and a smooth system 
$(\f S, \zet, \eps)$
of smooth sections,
where
$\eps: \f S \ucar{\f B} \f F \to \f G$
(see
Definition \ref{Definition: systems of sections}).

\myskip

\bDf
\label{Definition: lift fibred manifold}
\index{smooth systems!lifted fibred manifold}
 We define the \emp{lifted smooth fibred manifold} of the smooth fibred manifold
$p : \f F \to \f B$
to be the fibred product over
$\f B$
\beq 
\f F\Upa \byd \wob(s_b, f_b) \in \f S \ucar{\f B} \f F \st 
s_b \in \f S_b, \; f_b \in \f F_b, \; b \in \f B\wcb 
= \f S \ucar{\f B} \f F \,,
\eeq
which can be regarded as the pullback of the smooth fibred manifold
$p : \f F \to \f B \,,$
with respect to the smooth projection
$\zet : \f S \to \f B \,.$

 The natural map
\beq
p\Upa : \f F\Upa \to \f S : (s_b, f_b) \to s_b
\eeq
makes
$\f F\Upa$
a smooth fibred manifold over
$\f S \,.$

 In simple words, the smooth fibred manifold
$p\Upa : \f F\Upa \to \f S$
can be regarded as the ``extension" of the smooth fibred manifold
$p : \f F \to \f B$
obtained by ``extending" the base space
$\f B$
to
$\f S \,.$

 Then,
$\f F\Upa$
turns out to be also a smooth double fibred manifold
\newdiagramgrid{1}
{.9, .9, .8, .8, .4}
{1}
\bdg[grid=1]
\f F\Upa 
&\rTo^{\pro_2} &\f F
&\rTo^{p} &\f B &, 
\edg
according to the following commutative diagram
\newdiagramgrid{2x2}
{1.3, 1.3, .6}
{.7, .7}
\bdg[grid=2x2]
\f S \ucar{\f B} \f F &\rTo^{\pro_2} &\f F
\\
\dTo^{\pro_1} &&\dTo_{p}
\\
\f S &\rTo^{\zet} &\f B &.
\edg

 Hence, we can regard the smooth evaluation map
$\eps : \f S \ucar{\f B}\f F \to \f G$
as a smooth fibred morphism over
$\f F$
\beq
\eps : \f F\Upa \to \f G \,.
\eeq

 The induced fibred chart of
$p\Upa : \f F\Upa \to \f S$
is
$(x^\lam, w^A, y^i) \,.$\END
\eDf
%--------------------------------------------------------------------%
%\newpage
\markboth{\rm Chapter \thechapter. Systems of sections}{\rm \thesection. F--smooth systems of smooth sections}
\section{F--smooth systems of smooth sections}
\label{Section: F--smooth systems of smooth sections}
\markboth{\rm Chapter \thechapter. Systems of sections}{\rm \thesection. F--smooth systems of smooth sections}
%--------------------------------------------------------------------%
\bsm
 We discuss the notions of \emp{F--smooth systems of fibrewisely smooth sections}
$(\f S, \zet, \eps)$
and of its \emp{F--smooth tangent prolongation}
$(\E T\f S, \E T\zet, \E T\eps) \,.$

 Moreover, we compare the F--smooth and smooth structure of 
$\f S$
for the particular case of a smooth system of smooth sections.

 Eventually, we discuss the F--smooth differential operators.
\esm
%--------------------------------------------------------------------%
\subsection{F--smooth systems of smooth sections}
\label{F--smooth systems of smooth sections}
%--------------------------------------------------------------------%
\bsm
 In the previous 
Section \ref{Smooth systems of smooth sections}, 
we have studied the smooth systems
$(\f S, \zet, \eps)$
of smooth sections of a smooth double fibred manifold
$\f G \to \f F \to \f B \,.$

 Now, we analyse the generalised notion of \emp{F--smooth system} 
$(\f S, \zet, \eps)$
of \emp{fibrewisely smooth sections} of a smooth double fibred manifold
$\f G \to \f F \to \f B \,,$
by releasing the hypotheses of smoothness and finite dimension of
$\f S \,.$
\esm

 Let us consider a smooth double fibred manifold
$\f G \oset{q}{\to} \f F \oset{p}{\to} \f B \,.$

\bDf
\label{Definition: fibrewisely smooth tubelike sections}
 We denote by
(see
Definition \ref{Definition: tubelike sections})
\beq
\utub(\f F, \f G) \sub \wob s : \f F \to \f G \wcb
\eeq
the subsheaf consisting of tubelike sections
$s : \f F \to \f G \,,$
which are \emp{fibrewiselyly smooth}, i.e. which fulfill the following condition, \emp{without any further local smoothness requirement},

 - $s_b : \f F_b \to \f G_b$
is global and smooth, for each 
$b \in \f B \,.$

 Thus, the sheaf of tubelike smooth sections turns out to be a subsheaf
(see
Definition \ref{Definition: tubelike sections})
\beq
\tub(\f F, \f G) \sub \utub(\f F, \f G) \,.\END
\eeq
\eDf

\myskip

 The following Definition is a generalisation of
Definition \ref{Definition: systems of sections},
as here we do not require that
$\f S$
be a finite dimensional smooth manifold
(hence, that the maps
$\zet$
and
$\eps$
be smooth).

\myskip

\bDf
\label{Definition: F--smooth systems of smooth sections}
\index{F--smooth systems!F--smooth system of smooth sections} 
\index{F--smooth systems!injective F--smooth system of smooth sections}
 We define an \emp{F--smooth system of fibrewisely smooth sections}
of the smooth double fibred manifold
$\f G \to \f F \to \f B$
to be a 3--plet
$(\f S, \zet, \eps) \,,$
where 

 1) $\f S$
is a set,

 2) $\zet : \f S \to \f B$ 
is a fibred set (i.e. 
$\f S$
is a set and
$\zet$
a surjective map, without any smoothness requirements),

 3) $\eps : \f S \ucar{\f B} \f F \to \f G$
is a fibred map over 
$\f F \,,$
according to the following commutative diagram
\newdiagramgrid{2x2}
{1.2, 1.2, .6}
{.7, .7}
\bdg[grid=2x2]
\f S \ucar{\f B} \f F &\rTo^{\eps} &\f G
\\
\dTo^{\pro_2} &&\dTo_{q}
\\
\f F &\rTo^{\id_\f F} &\;\f F &,
\edg
which fulfills the following condition:

 *) for each
$s \in \f S_b \,,$
with
$b \in \f B \,,$
the induced section
\beq
\eps_s : \f F_b \to \f G_b
\eeq
of the restricted smooth fibred manifold
$q_b : \f G_b \to \f F_b$
is \emp{smooth and globally defined on}
$\f F_b \,.$

 The map
$\eps : \f S \ucar{\f B} \f F \to \f G$
is called the \emp{evaluation map} of the system.

 We denote by
\beq
\usec(\f B, \f S) \sub \wob s : \f B \to \f S \wcb
\eeq
the subsheaf consisting of \emp{local} sections
$s : \f B \to \f S \,,$
\emp{without any smoothness requirement}.

 Thus, the evaluation map
$\eps$
yields the sheaf morphism
\beq
\eps_\f S : \usec(\f B, \f S) \to \utub(\f F, \f G) : 
\sig \mto \br\sig \,,
\eeq
where, for each
$\sig \in \usec(\f B, \f S) \,,$ 
the tubelike fibrewisely smooth section
$\br\sig$
is defined by
\beq
\br\sig : \f F_b \to \f G_b : f_b \mto \eps \w(\sig(b), f_b\w) \,,
\ssep{for each}
b \in \f S \,.
\eeq

 Therefore, the map
$\eps_\f S : \usec(\f B, \f S) \to \utub(\f F, \f G)$
provides a \emp{selection} of the tubelike fibrewisely smooth sections
$\phi : \f F \to \f G \,,$
given by the subset
\beq
\utub_\f S(\f F, \f G) \byd \eps_\f S\w(\usec(\f B, \f S)\w) \sub \utub(\f F, \f G) \,.
\eeq

\myskip

 The smooth system of fibrewisely smooth sections
$(\f S, \zet, \eps)$
is said to be \emp{injective} if the map
$\eps_\f S : \usec(\f B, \f S) \to \utub(\f F, \f G)$
is injective, i.e. if, for each
$\sig \,, \ac\sig \in \usec(\f B, \f S) \,,$
\beq
\br\sig \eqv \eps_\f S (\sig) = \brac\sig \eqv \eps_\f S(\ac\sig)
\Rarr
\sig = \ac\sig \,.
\eeq

 If the system is injective, then we obtain the bijection
\beq
\eps_\f S : 
\usec(\f B, \f S) \to \utub_\f S(\f F, \f G) :
\sig \mto \br\sig \,,
\eeq
whose inverse is denoted by
\beq
(\eps_\f S)^{-1} : 
\utub_\f S(\f F, \f G) \to \usec(\f B, \f S) :
\phi \mto \wha\phi \,.\END
\eeq
\eDf

\myskip

 Indeed, we are essentially interested in injective F--smooth systems of tubelike fibrewisely smooth sections.

\myskip

\bRm
 We stress that we have not assumed a priori any smooth or F--smooth structure on
$\f S \,.$
 On the other hand, the specification ``\emp{F--smooth}" system in the above
Definition \ref{Definition: F--smooth systems of smooth sections}
will be justified later by
Theorem \ref{Theorem: F--smooth structure of F--smooth systems of 
smooth sections}.

 Moreover, we stress that the local sections
$\sig : \f B \to \f S \,,$
without any F--smoothness requirement, yield tubelike sections
$\br \sig : \f F \to \f G \,,$
which need not to be smooth, even if their restrictions
$\br\sig_b : \f F_b \to \f G_b \,,$ 
are smooth, for each
$b \in \f B \,,$
(according to condition *) in
Definition \ref{Definition: F--smooth systems of smooth sections}).

 Later, we shall see that the set
$\f S$
is an F--smooth space in a natural way
(see
Theorem \ref{Theorem: F--smooth structure of F--smooth systems of 
smooth sections})
and that the selected fibrewisely smooth sections
$\br\sig \in \utub_\f S(\f F, \f G)$
turn out to be smooth if and only if their source sections
$\sig \in \usec(\f B, \f S)$
are F--smooth
(see
Theorem \ref{Theorem: sections br sig are smooth if and only if sig are F--smooth}).

 This result further justifies the name ``\emp{F--smooth systems of fibrewisely smooth sections}" in
Definition \ref{Definition: F--smooth systems of smooth sections}.\END
\eRm

\myskip

 Let us examine some examples of F--smooth systems of fibrewisely smooth sections.

 We start by considering ``infinite dimensional" examples, by following the analogous thread of smooth systems of global smooth maps
(see
\S\ref{Systems of sections}). 

\myskip

\bEx
\label{Example: F--smooth system of all tubelike smooth sections}
 For each smooth double fibred manifold
$\f G \to \f F \to \f B \,,$
the set
(see
Definition \ref{Definition: fibrewisely smooth tubelike sections})
\beq
\utub_\f S(\f F, \f G) \byd \utub(\f F, \f G)
\eeq
consisting of \emp{all} 
$\br\phi \in \utub(\f F, \f G)$
yields in a natural way an injective F--smooth system of fibrewisely smooth sections.\END
\eEx

\bEx
\label{Example: F--smooth system of all tubelike smooth sections with compact support}
 For each smooth double fibred manifold
$\f G \to \f F \to \f B \,,$
where
$\f G \to \f F$
is a \emp{vector bundle}, the subset
\beq
\utub_\f S(\f F, \f G) \sub \utub(\f F, \f G)
\eeq
consisting of \emp{all}
$\br\phi \in \utub(\f F, \f G) \,,$
whose fibrewisely restrictions
$\br\sig_b : \f F_b \to \f G_b$
have \emp{compact support} for each
$b \in \f B \,,$
yields in a natural way an injective F--smooth system of tubelike fibrewisely smooth sections.

 Indeed, this system is a subsystem of the system of sections
$\phi : \f F \to \f G$
considered in the above 
Example \ref{Example: F--smooth system of all tubelike smooth sections}.\END
\eEx

\myskip

 We can reconsider ``finite dimensional" examples of smooth systems of tubelike fibrewisely smooth sections in the present context of F--smooth spaces
(see
Section \ref{Smooth systems of smooth sections}). 
 
\myskip

\bEx
\label{Example: F--smooth system of linear sections}
 For each smooth double fibred manifold
$\f G \to \f F \to \f B \,,$
where
$\f F \to \f B$
and
$\f G \to \f B$
are \emp{vector bundles}, 
the subset
\beq
\utub_\f S(\f F, \f G) \sub \utub(\f F, \f G)
\eeq
consisting of \emp{all}
$\br\phi \in \utub(\f F, \f G) \,,$
whose fibrewise restrictions
$\br\sig_b : \f F_b \to \f G_b$
are \emp{linear maps} for each
$b \in \f B \,,$
yields in a natural way an injective F--smooth system of fibrewisely smooth sections.

 Indeed, the F--smooth structure of
$\f S$
turns out to be just the natural F--smooth structure underlying the natural smooth structure of
$\f S$
(see, later,
Theorem \ref{Theorem: F--smooth structure of F--smooth systems of smooth sections}).\END
\eEx

\bEx
\label{Example: F--smooth system of affine sections}
 For each smooth double fibred manifold
$\f G \to \f F \to \f B \,,$
where
$\f F \to \f B$
and
$\f G \to \f B$
are \emp{affine bundles}, 
the subset
\beq
\utub_\f S(\f F, \f G) \sub \utub(\f F, \f G)
\eeq
consisting of \emp{all}
$\br\phi \in \utub(\f F, \f G) \,,$
whose fibrewise restrictions
$\br s_b : \f F_b \to \f G_b \,,$
are \emp{affine maps} for each
$b \in \f B \,,$
yields in a natural way an injective F--smooth system of fibrewisely smooth sections.

 Indeed, the F--smooth structure of
$\f S$
turns out to be just the natural F--smooth structure underlying the natural smooth structure of
$\f S$
(see, later,
Theorem \ref{Theorem: F--smooth structure of F--smooth systems of smooth sections}).\END
\eEx

\bEx
\label{Example: F--smooth system of polynomial sections of a given degree}
 For each smooth double fibred manifold
$\f G \to \f F \to \f B \,,$
where
$\f F \to \f B$
and
$\f G \to \f B$
are \emp{affine bundles}, 
the subset
\beq
\utub_\f S(\f F, \f G) \sub \utub(\f F, \f G)
\eeq
consisting of \emp{all}
$\br\phi \in \utub(\f F, \f G) \,,$
whose fibrewise restrictions
$\br\phi_b : \f F_b \to \f G_b$
are \emp{polynomial maps} of a \emp{given degree} for each
$b \in \f B \,,$
yields in a natural way an injective F--smooth system of fibrewisely smooth sections.

 Indeed, the F--smooth structure of
$\f S$
turns out to be just the natural F--smooth structure underlying the natural smooth structure of
$\f S$
(see, later,
Theorem \ref{Theorem: F--smooth structure of F--smooth systems of smooth sections}).\END
\eEx

\myskip

 We can also reconsider the ``infinite dimensional" example of F--smooth systems of polynomial sections of any degree in the present context of F--smooth spaces
(see
\S\ref{Systems of sections}). 
 
\myskip

\bEx
\label{Example: F--smooth system of polynomial sections of any degree}
 For each smooth double fibred manifold
$\f G \to \f F \to \f B \,,$
where
$\f F \to \f B$
and
$\f G \to \f B$
are \emp{affine bundles}, 
the subset
\beq
\utub_\f S(\f F, \f G) \sub \utub(\f F, \f G)
\eeq
consisting of \emp{all}
$\br\sig \in \utub(\f F, \f G) \,,$
whose fibrewise restrictions
$\br\sig_b : \f F_b \to \f G_b \,,$
are \emp{polynomial maps} of \emp{any degree} for each
$b \in \f B \,,$
yields in a natural way an injective F--smooth system of fibrewisely smooth sections.

 Indeed, the F--smooth structure of
$\f S$
turns out to be just the natural F--smooth structure underlying the natural infinite dimensional smooth structure of
$\f S$
(see, later,
Theorem \ref{Theorem: F--smooth structure of F--smooth systems of smooth sections}).\END
\eEx
%--------------------------------------------------------------------%
\subsection{F--smooth structure of $\f S$}
\label{F--smooth structure of S}
%--------------------------------------------------------------------%
\bsm
 We show that an F--smooth system 
$(\f S, \zet, \eps)$
of fibrewisely smooth sections turns out to have a natural F--smooth structure.

 Namely, the set
$\f S$
has a natural F--smooth structure and the maps
\beq
\zet : \f S \to \f B
\ssep{and}
\eps : \f S \ucar{\f B} \f F \to \f G
\eeq
turn out to be F--smooth.

 Furthermore, we exhibit a bijection between F--smooth sections
$\sig : \f B \to \f S$
and tubelike smooth sections
$\br\sig : \f F \to \f G \,.$
\esm

 Let us consider an F--smooth system of fibrewisely smooth sections
$(\f S, \zet, \, \eps) \,.$

\myskip

 Let us start by exhibiting a natural F--smooth structure of the set
$\f S \,.$
 Preliminarily, we need a few technical Lemmas, which provide some pullback objects.
 
\myskip

\bLm
\label{Lemma: c*(F)}
 If
$c : \f I_c \to \f B$
is a smooth curve, then we obtain the smooth submanifold
\beq
c^* (\f F) \byd
\{(\lam, f) \in \f I_c \car \f F \sst c(\lam) = p(f)\} \sub
\f I_c \car \f F \,,
\eeq
along with the smooth projections
\beq
c^*(p) : c^*(\f F) \to \f I_c : (\lam, f) \mto \lam
\ssep{and}
c^*_\f F : c^*(\f F) \to \f F : (\lam, f) \mto f \,.\END
\eeq
\eLm

\bLm
\label{Lemma: gam*}
 If
$c : \f I_c \to \f B$
and
$\gam : \f I_{\gam} \to \f I_c$
are smooth maps, then we obtain the smooth map
\beq
\gam^* : (c \com \gam)^* (\f F) \to c^* (\f F) :
(\lam, f) \mto \w(\gam(\lam), f\w) \,,
\eeq
which provides just a smooth reparametrisation of
$c^* (\f F) \,.$\END
\eLm

\bLm
\label{Lemma: c* : c*(F) to G}
 If
$\wha c : \f I_{\wha c} \to \f S$
is a curve, which projects on a smooth curve
$c \byd \zet \com \wha c : \f I_{\wha c} \to \f B \,,$
then we obtain the map
\beq
\wha c^*(\eps) : c^*(\f F) \to \f G :
(\lam, f) \mto \eps\w(\wha c(\lam), f\w) \,.\END
\eeq
\eLm

\myskip

\bTh
\label{Theorem: F--smooth structure of F--smooth systems of smooth sections}
 Let us consider the set
\beq
\C C \byd \{\wha c : \f I_\wha c \to \f S\}
\eeq
consisting of all curves
$\wha c : \f I_{\wha c} \to \f S \,,$
such that the following induced maps between smooth manifolds be smooth
(see the above
Lemma \ref{Lemma: c*(F)}
and
Lemma \ref{Lemma: c* : c*(F) to G})
\bat{2}
(a)
&\qquad\qquad\quad\;\;\;
c
&&: 
\f I_\wha c \to \f B : \lam \mto \zet \w(\wha c(\lam)\w) \,,
\\
(b)
&\qquad\qquad\;
\wha c^*(\eps)
&&: 
c^*(\f F) \to \f G : 
(\lam, f) \mto \eps\w(\wha c(\lam), f\w) \,.
\end{alignat*}

 Then, the pair
$(\f S, \C C)$
turns out to be an F--smooth space.
\eTh

\bpf
 Let us prove that the pair
$(\f S, \C C)$
be an F--smooth space, by showing that it fulfills the two conditions 1) and 2) in 
Definition \ref{Definition: F--smooth spaces}.

\myskip

 1) For every element
$s \in \f S \,,$
let us consider the constant curve
$\wha c : \f I_{\wha c} \to \f S : \lam \mto s \,,$
which clearly passes through
$s \,.$
 Then, the following facts hold.
 
 a) The induced curve 
$c \byd \zet \com \wha c : \f I_{\wha c} \to \f B$
turns out to be constant as well, hence smooth.

 b) The map
\beq
\wha c^*(\eps) : c^*(\f F) \to \f G : 
(\lam, f) \mto \eps\w(\wha c(\lam), f\w)
\eeq
can be regarded as the map
\beq
\eps_s : \f F_{\zet(s)} \to \f G_{\zet(s)} \,,
\eeq
hence, in virtue of condition *) in
Definition \ref{Definition: F--smooth systems of smooth sections}, 
it turns out to be smooth.
 
 Therefore, the constant curves
$\wha c : \f I_\wha c \to \f S$
belong to
$\C C \,.$

\myskip

 2) If
$\gam : \f I_\gam \to \f I_\wha c$
is any smooth curve, then 

 a) the map
\beq
c \com \gam : \f I_\wha c \to \f B
\eeq
is a composition of smooth maps between smooth manifolds, hence it turns out to be smooth,

 b) in virtue of the above
Lemma \ref{Lemma: gam*}, 
the map
\beq
(\wha c \com \gam)^*(\eps) : 
(c \com \gam)^* (\f F) \to \f G : 
(\lam, f) \mto \eps \W(\wha c \w((\gam(\lam)\w), \, f\W)
\eeq
is given by the composition of smooth maps between smooth manifolds
\newdiagramgrid{2x3}
{2, 1.5, 2, 2, 1.7}
{.9, .9}
\bdg[grid=2x3]
(c \com \gam)^* (\f F) &\rTo^{\gam^*} &c^*(\f F) 
&\rTo^{\wha c^*(\eps)} &\f G &:
\\
\W(\lam, f_{c\w(\gam(\lam)\w)}\W)
&\rTo^{\gam^*} 
&\W(\gam(\lam), f_{c\w(\gam(\lam)\w)}\W)
&\rTo^{\wha c^*(\eps)} 
&\eps\W(\wha c\w(\gam(\lam)\w), f_{c\w(\gam(\lam)\w)}\W) &,
\edg
hence it turns out to be smooth.
 
 Therefore, according to the above conditions (a) and (b), the curves
$\wha c \com \gam : \f I_\wha c \to \f S$
belong to
$\C C \,.$\QED
\epf

\myskip

 Then, we show that the maps
$\zet$
and 
$\eps$
are F--smooth.
 Preliminarily, we need a technical Lemma, which provides some pullback objects.

\bLm
\label{Lemma: factorisation of (c, CF}
 Let us consider an F--smooth curve
\beq
(\wha c, c_\f F) : \f I_{\wha c} \to \f S \ucar{\f B} \f F \,,
\eeq
where
\beq
\wha c : \f I_{\wha c} \to \f S
\ssep{and}
c_\f F : \f I_{\wha c} \to \f F
\eeq
are, respectively, an F--smooth curve and a smooth curve, which project on the same smooth base curve
$c : \f I_{\wha c} \to \f B \,.$

 Then, the F--smooth curve 
$(\wha c, c_\f F) : \f I_{\wha c} \to \f S \car_{\f B} \f F$
factorises through a smooth curve
(see
Lemma \ref{Lemma: c*(F)})
\beq
(c, c_\f F) : \f I_{\wha c} \to  c^*(\f F) \,,
\eeq
according to the following commutative diagram
\newdiagramgrid{2x2}
{1.3, 1.3, .7}
{.7, .7}
\bdg[grid=2x2]
\f I_{\wha c} &\rTo^{(\wha c, c_\f F)} &\f S \ucar{\f B} \f F
\\
\dTo^{\id_{\f I_{\wha c}}} &&\uTo_{\wha c \car \id_\f F}
\\
\f I_{\wha c} &\rTo^{(\id_{\f I_{\wha c}}, c_\f F)} &c^*(\f F) &.
\edg
\eLm

\bpf
 In fact, for each
$\lam \in \f I_{\wha c} \,,$
the following diagram commutes
\newdiagramgrid{2x2}
{1.3, 1.5, 1.5}
{.7, .7}
\bdg[grid=2x2]
\lam &\rTo^{(\wha c, c_\f F)} &\w(\wha c(\lam), c_\f F(\lam)\w)
\\
\dTo^{\id} &&\uTo_{\wha c \car \id_\f F}
\\
\lam &\rTo &\w(\lam, c_\f F(\lam)\w) &.\QED
\edg
\epf

\bPr
\label{Proposition: zet and eps are F--smooth}
 The maps
\beq
\zet : \f S \to \f B
\ssep{and}
\eps : \f S \ucar{\f B} \f F \to \f G
\eeq
turn out to be F--smooth.
\ePr

\bpf
 According to
Definition \ref{Definition: F--smooth maps},
we have to prove that
$\zet$
and
$\eps$
map F--smooth curves of the source space into F--smooth curves of the target space.

\myskip

 1) For each F--smooth curve
$\wha c : \f I_{\wha c} \to \f S \,,$
the composed map
$\zet \com \wha c : \f I_{\wha c} \to \f B$
is a smooth curve, by assumption, hence it is an F--smooth curve
(see
Definition \ref{Definition: natural F--smooth structure of smooth manifolds}).

\myskip

 2) Let
$(\wha c, c_\f F) : \f I_{\wha c} \to \f S \ucar{\f B} \f F$
be an F--smooth curve, where
\beq
\wha c : \f I_{\wha c} \to \f S
\ssep{and}
c_\f F : \f I_{\wha c} \to \f F
\eeq
are, respectively  an F--smooth curve and a smooth curve, which project on the same smooth base curve
$c : \f I_{\wha c} \to \f B \,.$

 Then, in virtue of the above
Lemma \ref{Lemma: factorisation of (c, CF},
the curve
\beq
\eps \com (\wha c, c^*_\f F) : \f I_{\wha c} \to \f G
\eeq
turns out to be the composition of smooth curves
\beq
\eps \com (\wha c, c^*_\f F) = 
c^*(\eps) \com (\id_{\f I_{\wha c}}, c_\f F) \,,
\eeq
according to the following commutative diagram
\newdiagramgrid{2x3}
{1.3, 1.3, 1.3, 1.3, .7}
{.7, .7}
\bdg[grid=2x3]
\f I_{\wha c} &\rTo^{(\wha c, c_\f F)} &\f S \ucar{\f B} \f F &\rTo^{\eps} &\f G
\\
\dTo^{\id_{\f I_{\wha c}}} &&\uTo_{\wha c \car \id_\f F} &&\uTo_{\id_\f G}
\\
\f I_{\wha c} &\rTo^{(\id, c_\f F)} &c^*(\f F) &\rTo^{c^*(\eps)} &\f G &,
\edg
hence it is smooth.\QED
\epf

\myskip

 Next, we show a natural bijection between local F--smooth sections
$\sig : \f B \to \f S$
and tubelike smooth sections
$\br\sig : \f F \to \f G \,.$
 Preliminarily, we need a technical Lemma.

\myskip

\bLm
\label{Lemma: pullback section induced by a curve}
 Let us consider a section
$\sig \in \usec(\f B, \f S) \,,$
the induced tubelike fibrewisely smooth section
$\br\sig \in \utub_\f S(\f F, \f G)$
and a smooth curve
$c : \f I_c \to \f B$
(see
Definition \ref{Definition: F--smooth systems of smooth sections}).

 Then, we obtain, by pullback, the fibred morphism over
$\f F$
\beq
c^*(\br\sig) : c^*(\f F) \to \f G \,,
\eeq
given by the composition
(see
Lemma \ref{Lemma: c*(F)})
\newdiagramgrid{1x5}
{1.5, 1.5, 1.3, 1.3, 1.3, 1.3, 1.2, .7, .6}
{.7}
\bdg[grid=1x5]
c^*(\f F) &\rTo^{(c^*(p), c^*_\f F)} &\f I_c \car \f F
&\rTo^{c \car \id_\f F} &\f B \car \f F
&\rTo^{\sig \car \id_\f F} &\f S \ucar{\f B} \f F
&\rTo^{\eps} &\f G &.\END
\edg
\eLm

%\myskip

 We denote by
\beq
\Fsec(\f B, \f S) \sub \usec(\f B, \f S)
\eeq
the subsheaf consisting of \emp{F--smooth} local sections of the fibred set
$\zet : \f S \to \f B \,.$

 Moreover, we denote by
\beq
\tub_\f S(\f F, \f G) \byd 
\utub_\f S(\f F, \f G) \cap \tub(\f F, \f G)
\eeq
the subsheaf consisting of tubelike \emp{smooth} sections of the smooth fibred manifold
$\f G \to \f F \,,$
which are selected by the F--smooth system
$(\f S, \zet, \eps) \,.$

\myskip

\bTh 
\label{Theorem: sections br sig are smooth if and only if sig are F--smooth}
 Let
$\sig \in \usec(\f B, \f S)$
be a local section
and
\beq
\br\sig \byd \eps_\f S(\sig) \in \utub_\f S(\f F, \f G)
\eeq
the induced tubelike section
(see
Definition \ref{Definition: F--smooth systems of smooth sections}).

 Then,
$\br\sig$
is smooth if and only if
$\sig$
is F--smooth.

 In other words, the sheaf morphism
\beq
\eps_\f S : \usec(\f B, \f S) \to \utub_\f S(\f F, \f G)
\eeq
restricts to a sheaf \emp{isomorphism} (denoted by the same symbol)
\beq
\eps_\f S : \Fsec(\f B, \f S) \to \tub_\f S(\f F, \f G) \,,
\eeq
according to the following diagram commutes
\newdiagramgrid{2x2}
{1.7, 1.7, 1.3, .6}
{.7, .7}
\bdg[grid=2x2]
\usec(\f B, \f S) &\rTo^{\eps_\f S} &\utub_\f S(\f F, \f G)
\\
\uInto^{\cup} &&\uInto_{\cup}
\\
\Fsec(\f B, \f S) &\rTo^{\eps_\f S} &\;\tub_\f S(\f F, \f G) &.
\edg
\eTh

\bpf
 1) Let us prove that, if
$\sig$
is F--smooth, then
$\br\sig$
is smooth.

 In fact, the map between smooth manifolds
$\br\sig : \f F \to \f G$
is given by a composition of F--smooth maps, according to the following commutative diagram
(see
Definition \ref{Definition: natural F--smooth structure of smooth manifolds}
and
Proposition \ref{Proposition: zet and eps are F--smooth})
\newdiagramgrid{2x2}
{1.5, 1.3, 1.5, 1.5, 1.5, .4, .6}
{.7, .7}
\bdg[grid=2x2]
\f F &\rTo^{(p, \id_\f F)} &\f B \car \f F
&\rTo^{\sig \car \id_\f F} &\f S \ucar{\f B} \f F
&\rTo^{\eps} &\f G
\\
\uTo^{\id_\f F} &&&&&&\dTo_{\id_\f G}
\\
\f F &&\rTo^{\br\sig} &&&&\;\f G &.
\edg
 
 Hence, according to
Proposition \ref{Proposition: composition of global F--smooth maps},
$\br\sig : \f F \to \f G$
is F-smooth.
 On the other hand, in virtue of
Definition \ref{Definition: natural F--smooth structure of smooth manifolds}
and
Theorem \ref{Theorem: relation between smoothness and F--smoothness},
it means that the map
$\br\sig : \f F \to \f G$
is smooth.

\myskip

 2) Let us prove that, if
$\br\sig$
is smooth, then
$\sig$
is F--smooth.

 By definition of F--smooth map between F--smooth spaces
(see
Definition \ref{Definition: F--smooth maps}), 
we have to prove that, for smooth curve
$c : \f I_{\wha c} \to \f B \,,$
the composed curve
\beq
\wha c_\sig \byd \sig \com c : \f I_{\wha c} \to \f S
\eeq
be F--smooth.

 Hence, according to
Theorem \ref{Theorem: F--smooth structure of F--smooth systems of 
smooth sections},
we have to prove that the induced maps between smooth manifolds
\beq
c_\sig \byd \zet \com \wha c_\sig : \f I_{c_\sig} \to \f B
\ssep{and}
(\wha c_\sig)^* (\eps) : (c_\sig)^* (\f F) \to \f G \,,
\eeq
be smooth.

 a) The curve
$c_\sig \byd \zet \com \wha c_\sig : \f I_{c_\sig} \to \f B$
is F--smooth because it is just the smooth curve
$c_\sig : \f I_{\wha c} \to \f B \,.$

 b) The map
$(\wha c_\sig)^*(\eps) : (c_\sig)^* (\f F) \to \f G$
turns out to be smooth, because it is the restriction to the smooth sub manifold over
$c : \f I_{\wha c} \to \f B$
\beq
(c_\sig)^* (\f F) \sub \f F \,,
\eeq
according to the following commutative diagram
\newdiagramgrid{2x2}
{1.5, 1.5, 1.3, .7}
{.7, .7}
\bdg[grid=2x2]
(c_\sig)^* (\f F) &\rTo^{(\wha c_\sig)^*(\eps)} &\f G
\\
\dTo^{(c_\sig)^*_\f F} &&\dTo_{\id_\f G}
\\
\f F &\rTo^{\br\sig} &\f G &.\QED
\edg
\epf

\bRm
\label{Remark: meaning of the Theorem usec utub}
 The above
Theorem \ref{Theorem: sections br sig are smooth if and only if sig are F--smooth}
clarifies the name ``F--smooth system of fibrewisely \emp{smooth} sections" in
Definition \ref{Definition: F--smooth systems of smooth sections}.

 Roughly speaking, we can interpret the above
Theorem \ref{Theorem: sections br sig are smooth if and only if sig are F--smooth}
in the following way.

 Given a section
$\sig \in \usec(\f B, \f S) \,,$
the induced section
$\br\sig \in \utub(\f F, \f G)$
is smooth along the fibres of the fibred manifold
$p : \f F \to \f B$
in virtue of condition *) in
Definition \ref{Definition: F--smooth systems of smooth sections}.

 Then, in order to check whether
$\br\sig \in \utub(\f F, \f G)$
is fully smooth, we have to show that it is smooth ``transversally", i.e. along smooth curves
$c : \f I_c \to \f B \,.$

 On the other hand, in virtue of the definition of F-smooth map
(see
Definition \ref{Definition: F--smooth maps}),
the section
$\sig \in \usec(\f B, \f S)$
is F--smooth if and only if its composition with any smooth curve
$c : \f I_{\wha c} \to \f B$
is smooth.

 Thus, to check the smoothness of both 
$\br\sig \in \utub(\f F, \f G)$
and the F--smoothness of\lnb
$\sig \in \usec(\f B, \f S)$
is subject to analysing the behaviour of both sections along smooth curves
$c : \f I_{\wha c} \to \f B \,.$

\myskip

 In the particular case when
$\f B = \Rn \,,$
the above 
Theorem \ref{Theorem: sections br sig are smooth if and only if sig are F--smooth}
reduces to a tautology.

 In fact, a local section
$\sig : \f B \to \f S$
turns out to be a local curve
$\wha c \byd \sig : \Rn \to \f S \,.$
 Moreover, in this case we have locally
$c^*(\f F) = \f F \,.$

 So, checking that
$\wha\sig$
maps F--smooth curves of
$\f B$
into F--smooth curves of
$\f S$
reduces to check that the fibrewisely smooth section
$\sig : \f F \to \f G$
be smooth.
 Indeed, no other independent check is required.
 In other words, in the particular case when
$\f B = \Rn \,,$
the above Theorem reduces to say that
$\sec(\f B, \f S)$
is just, by definition, the subsheaf
$\sec(\f B, \f S) \sub \usec(\f B, \f S)$
which yields
$\tub(\f F, \f G) \,.$\END
\eRm

\myskip

\bRm
\label{Remark: smooth structure and F--smooth structure of systems of sections}
 In our definition of ``smooth system of smooth sections"
(see
Definition \ref{Definition: systems of sections}),
we have assumed a priori a finite dimensional smooth structure of
$\f S \,,$
while, for our concept of ``F--smooth system of fibrewisely smooth sections"
(see
Definition \ref{Definition: F--smooth systems of smooth sections}),
we have not assumed a priori any smooth structure of
$\f S \,,$
but have recovered an F--smooth structure in a unique way
(see
Theorem \ref{Theorem: F--smooth structure of F--smooth systems of 
smooth sections}).

\myskip

 Now, the procedure we have used for F--smooth systems of fibrewisely smooth sections to recover the F--smooth structure of
$\f S$
can be applied to smooth systems of smooth sections as well and a natural question arises: ``is the assumed smooth structure of
$\f S$
compatible with the recovered F--smooth structure?"
 
 In general, the answer is negative.
 In fact, in
Example \ref{Example: different smooth structures of a system of sections},
we have shown that we can assume several smooth structures on
$\f S$
providing the same system of smooth sections.\END
\eRm

\myskip

 Next, we show that, in the case when the smooth fibred manifold
$q : \f G \to \f F$
is an affine (vector) bundle, the F--smooth fibred space
$\zet : \f S \to \f B$
of an injective system
$(\f S, \zet, \eps)$
of smooth sections
$\phi : \f F \to \f G$
inherits in a natural way an affine (vector) structure.

\myskip

\bPr
\label{Proposition: vector structure of S in the case when G to F is a vector bundle}
 Let us suppose that the fibred manifold
$q : \f G \to \f F$
be a vector bundle and consider an \emp{injective} F--smooth system
$(\f S, \zet, \eps)$
of fibrewisely smooth sections
$\phi : \utub(\f F, \f G)$
(see
Example \ref{Example: F--smooth system of all tubelike smooth sections}).

 Then, the fibres of the F--smooth fibred space
$\zet : \f S \to \f B$
inherit in a natural way a vector structure given, for each
$k \in \Rn$
and
$s, \, \ac s \in \f S_b \,,$
with
$b \in \f B \,,$ 
by
\beq
k \, s \byd \wha{k \, \br s}
\ssep{and}
s + \ac s \byd \wha{\br s + \brac s} \,.\END
\eeq 
\ePr

\bPr
\label{Proposition: affine structure of S in the case when G to F is an affine bundle}
 Let us suppose that the fibred manifold
$q : \f G \to \f F$
be an affine bundle associated with the vector bundle
$\ba q : \baf G \to \f F \,,$
and consider an \emp{injective} F--smooth system
$(\f S, \zet, \eps)$
of fibrewisely smooth sections
$\phi : \utub(\f F, \f G)$
and the associated injective F--smooth system
$(\baf S, \ba\zet, \ba\eps)$
of fibrewisely smooth sections
$\ba\phi : \utub(\f F, \baf G) \,.$

 Then, the fibres of the F--smooth fibred space
$\zet : \f S \to \f B$
inherit in a natural way an affine structure given, for each
$s \in \f S_b$
and
$\ba s \in \baf S_b$
with
$b \in \f B \,,$ 
by
\beq
s + \ba s \byd \wha{\br s + \br{\ba s}} \,.\END
\eeq 
\ePr
%--------------------------------------------------------------------%
\subsection{F--smooth tangent prolongation of $(\f S, \zet, \eps)$}
\label{F--smooth tangent prolongation of (S, zet, eps)}
%--------------------------------------------------------------------%
\bsm
 Now, we consider an F--smooth system
$(\f S, \zet, \eps)$
of fibrewisely smooth sections of a smooth double fibred manifold
$\f G \oset{q}{\to} \f F \oset{p}{\to} \f B$
(see
Definition \ref{Definition: F--smooth systems of smooth sections})
and introduce the concept of ``F--smooth tangent space"
$\E T\f S$
of the F--smooth space of parameters
$\f S \,,$
which is equipped with its family
$\C C$
of F--smooth curves
$\wha c : \f I_{\wha c} \to \f S$
(see
Theorem \ref{Theorem: F--smooth structure of F--smooth systems of smooth sections}).

 Our formal construction of
$\E T\f S$
reflects the intuitive idea, by which, for each
$s \in \f S_b \,,$
with
$b \in \f B \,,$
a tangent vector
$\E X_s \in \E T_s\f S$
is to be an ``infinitesimal variation" of the global smooth map
$\eps_s : \f F_b \to \f G_b \,.$
 It is remarkable the fact that this construction involves only smooth maps between smooth manifolds, by taking into account the smooth structure of the smooth double fibred manifold
$\f G \oset{q}{\to} \f F \oset{p}{\to} \f B \,.$

\myskip

 Actually, we define a tangent vector
$\E X_s$
of
$\f S \,,$
at
$s \in \f S \,,$
as an equivalence class of F--smooth curves
$\wha c : \f I_{\wha c} \to \f S \,,$
such that the induced smooth maps between smooth manifolds
$\wha c^*(\eps) : c^*(\f F) \to c^*(\f G)$
have a 1st order contact in
$s \,.$

 Then, we show that a tangent vector
$\E X_s$
can be naturally represented by a pair
$(u, \Xi_u) \,,$
where
$u \in T_b\f B$
and
$\Xi_u : (T\f F)_u \to (T\f G)_u$
is a suitable smooth section.

 Moreover, we show that the F--smooth tangent space
$\E T\f S$
is equipped with a natural F--smooth structure and exhibit the natural F--smooth maps
\beq
\tau_\f S : \E T\f S \to \f S
\ssep{and}
\E T_1\eps : \E T\f S \car \f M \to T\f N \,.
\eeq

 Indeed, the fibres of
$\tau_\f S : \E T\f S \to \f S$
are naturally equipped with a vector structure.

\myskip

 The procedure followed to achieve the above results is partially similar to that followed for the tangent space of the F--smooth systems of smooth maps between two smooth manifolds
(see
\S\ref{F--smooth tangent prolongation of (S, eps)}).
 However, the present context is more complex and requires an additional care.
 In the present case, the reason of the difficulty and of the consequent complication is due to the fact that the map
$\eps : \f S \ucar{\f B} \f F \to \f G$
acts on a fibred product
$\f S \ucar{\f B} \f F \,,$
not just on a product
$\f S \car \f F \,.$
 Hence, a curve
$\wha c : \f I_{\wha c} \to \f S \,,$
which moves the base points of
$\f S$
in
$\f B \,,$
moves at the same time also the base points of
$\f F$
(and of
$\f G$)
in 
$\f B \,.$
\esm

 Let us consider a smooth double fibred manifold
$\f G \oset{q}{\to} \f F \oset{p}{\to} \f B$
and denote the smooth fibred charts of
$\f G$
by
$(x^\mu, y^i, z^a) \,.$

 Moreover, let us consider an \emp{injective} F--smooth system
$(\f S, \zet, \eps)$
of fibrewisely smooth sections
$\phi \in \utub(\f F, \f G) \,,$
along with the set 
$\C C$
of F--smooth curves
$\wha c: \f I_\wha c \to \f S$
defined in
Theorem \ref{Theorem: F--smooth structure of F--smooth systems of 
smooth sections}.

\myskip

 Then, we define the tangent space of the F--smooth space
$(\f S, \C C) \,,$
via equivalence classes of suitable smooth maps between smooth manifolds, in the following way.

\myskip

\bLm
\label{Lemma: c1*(F) = c2*(F)}
 If
$(c_1, c_2) : (\f I_1, \f I_2) \to \f B$
is a pair of smooth curves and the pair
$(\lam_1, \lam_2) \in \f I_1 \car \f I_2$
is an element, such that
\beq
dc_1 (\lam_1) = 
dc_2 (\lam_2) \in T\f B \,,
\eeq
then we have
\beq
\w(c^*_1(\f F)\w)_{\lam_1} = \w(c^*_2(\f F)\w)_{\lam_2}
\ssep{and}
T\w(c^*_1(\f F)\w)_{\lam_1} = T\w(c^*_2(\f F)\w)_{\lam_2} \,.\END
\eeq
\eLm

\bLm
\label{Lemma: T c * eps}
 If
$\wha c : \f I_{\wha c} \to \f S$
is an F--smooth curve which projects on a smooth curve
$c : \f I_{\wha c} \to \f B \,,$
then we obtain the smooth map
(see
Theorem \ref{Theorem: F--smooth structure of F--smooth systems of smooth sections})
\beq
T \w(\wha c^*(\eps)\w) : T \w(c^*(\f F)\w) \to T\f G \,.\END
\eeq
\eLm

\myskip

 Then, we introduce the concept of 1st order contact for the F--smooth curves of the type
$\wha c : \f I_{\wha c} \to \f S \,.$

\myskip

\bDf
\label{Definition: 1st order contact of F--smooth curves of the F--smooth system of sections}
\index{F--smooth systems!1st order contact of curves}
 We say that two F--smooth curves
(see
Theorem \ref{Theorem: F--smooth structure of F--smooth systems of 
smooth sections})
\beq
\wha c_1 : \f I_1 \to \f S
\ssep{and} 
\wha c_2 : \f I_2 \to \f S \,,
\eeq
which project, respectively, on smooth curves
\beq
c_1 \byd \zet \com \wha c_1 : \f I_1 \to \f B
\ssep{and} 
c_2 \byd \zet \com \wha c_1: \f I_2 \to \f B \,,
\eeq
have a \emp{1st order contact} in
$(\lam_1, \lam_2) \in \f I_1 \car \f I_2$
if they fulfill the following conditions involving smooth manifolds and maps
\beq
dc_1 (\lam_1) = dc_2 (\lam_2) \,,\ltag{1}
\eeq
\beq
T_{\lam_1} \w(\wha c^*_1(\eps)\w) = 
T_{\lam_2} \w(\wha c^*_2(\eps)\w) \,,\ltag{2}
\eeq
i.e., in coordinates,
\beq
c^\mu_1(\lam_1) = x^\mu(b) = c^\mu_2(\lam_2) \,,\ltag{1'a}
\eeq
\beq
(\der_0 c^\mu_1) (\lam_1) 
\eqv \Xi^\mu_0 \eqv
(\der_0 c^\mu_2) (\lam_2) \,,\ltag{1'b}
\eeq
\beq
\der_0 \w(\wha c_1^* (\eps)^a\w)_{|\lam_1} 
\eqv \Xi^a_0 \eqv
\der_0 \w(\wha c_2^* (\eps)^a\w)_{|\lam_2} \,,\ltag{2'a}
\eeq
\beq
\der_i \w(\wha c_1^* (\eps)^a\w)_{|\lam_1}
\eqv \Xi^a_i \eqv
\der_i \w(\wha c_2^* (\eps)^a\w)_{|\lam_2} \,.\ltag{2'b}
\eeq

 Clearly, the above 1st order contact yields equivalence relations
$\sim$
in the sets of pairs
\beq
\w(c : \f I_{\wha c} \to \f B, \;\; \lam \in \f I_{\wha c}\w)
\ssep{and}
\w(\wha c : \f I_{\wha c} \to \f S, \;\; \lam \in \f I_{\wha c}\w) \,,
\eeq
where
$\wha c : \f I_{\wha c} \to \f S$
are F--smooth curves of
$\f S$
and
$c \byd \zet \com \wha c : \f I_{\wha c} \to \f B$
are the associated smooth curves of
$\f B \,.$\END
\eDf

\myskip

 Then, we define the tangent vectors of
$\f S$
via 1st order equivalence classes of F--smooth curves of
$\wha c : \f I_{\wha c} \to \f S$
(see
Theorem \ref{Theorem: F--smooth structure of F--smooth systems of 
smooth sections}).

\myskip

\bDf
\label{Definition: tangent space of the F--smooth space of smooth tubelike sections}
\index{F--smooth systems!tangent vectors}
\index{F--smooth systems!tangent space}
 We define a \emp{tangent vector} of
$\f S \,,$
at
$s \in \f S_b \,,$
with
$b \in \f B \,,$
to be an equivalence class
(see
Definition \ref{Definition: 1st order contact of F--smooth curves of the F--smooth system of sections})
\beq
\E X_s \byd \w[(\wha c, \lam)\w]_\sim
\eeq
where 
$\wha c : I_{\wha c} \to \f S$
are F--smooth curves, such that
\beq
\wha c(\lam) = s \,.
\eeq

 Then, we define:
 
 1) the \emp{tangent space} of
$\f S \,,$
at
$s \in \f S_b \,,$
with
$b \in \f B \,,$
to be the set of tangent vectors of
$\f S$
at
$s$
\beq
\E T_s \f S \byd \wob \E X_s \wcb \,,
\eeq

 2) the \emp{tangent space} of
$\f S$
to be the disjoint union
\beq
\E T\f S \byd \bigsqcup_{s \in \f S}{\E T_s}\f S \,.\END
\eeq
\eDf

\bRm
 In order to mimic the tangent space of standard manifolds, we call the elements
$\E X_s$
``tangent vectors".
 However, so far, we do not know yet that these objects are really elements of a vector space.
 This fact will be proved later in
Theorem \ref{Theorem: vector structure of TS for an F--smooth system of sections}.\END
\eRm

\myskip

 Thus, in virtue of the above 
Definition \ref{Definition: 1st order contact of F--smooth curves of the F--smooth system of sections}
and
Definition \ref{Definition: tangent space of the F--smooth space of smooth tubelike sections}, 
the tangent vectors
$\E X_s \in \E T_s\f S$
can be represented through suitable pairs
$(u, \Xi_u) \,,$
consisting of a base vector
$u \in T_b\f B$
and a smooth map
$\Xi_u : (T\f F)_u \to (T\f G)_u  \,,$
as follows.

\myskip

\bTh
\label{Theorem: representation of tangent space of space of 
F--smooth system of sections}
 Let
$(\f S, \zet, \eps)$
be an F--smooth system of fibrewisely smooth sections.
 
 Then, in virtue of the above 
Definition \ref{Definition: 1st order contact of F--smooth curves of the F--smooth system of sections}
and
Definition \ref{Definition: tangent space of the F--smooth space of smooth tubelike sections},
every tangent vector
\beq
\E X_s \byd [(\wha c, \lam)]_\sim \in \E T_s\f S \,,
\eeq
can be regarded as the pair
$(u, \Xi_u) \,,$
consisting of

 a) the base vector
\beq
u \byd \w[(c,\lam)\w]_\sim \in T_b\f B \,,
\eeq

 b) the smooth map
\beq
\Xi_u : (T\f F)_u \to (T\f G)_u \,,
\eeq
which is a global smooth section of the smooth fibred manifold
$Tq : (T\f G)_u \to (T\f F)_u$
and an affine fibred morphisms over
$s^*(\eps) : \f F_b \to \f G_b \,,$ 
whose smooth fibre derivative
\beq
D\Xi_u : (V\f F)_b \to (V\f G)_b
\eeq
fulfills the equality
\beq
D\Xi_u = T_b\w(s^*(\eps)\w) : T(\f F_b) \to T(\f G_b) \,.
\eeq

 Indeed, the following diagram commutes
\newdiagramgrid{3x2}
{1.3, 1.3, .7}
{.7, .7, .7, .7}
\bdg[grid=3x2]
(T\f F)_u &\rTo^{\id_{(T\f F)_u}} &(T\f F)_u
\\
\uTo^{\id_{T\f F_u}} &&\uTo_{Tq}
\\
(T\f F)_u &\rTo^{\Xi_u} &(T\f G)_u
\\
\dTo^{\tau_\f F} &&\dTo_{\tau_\f G}
\\
\f F_b &\rTo^{s^*(\eps)} &\;\;\f G_b &,
\edg

 We have the coordinate expressions
\bgt
u = u^\mu \, \der_\mu \,,
\\
(x^\mu, \, y^i, \, z^a) \com \Xi_u = 
\w(x^\mu(b), \, y^i, \, s^*(\eps)^a\w) \,,
\qquad
(\dt x^\mu, \, \dt y^i_{|u}, \, \dt z^a_{|u}) \com \Xi_u = 
\w(u^\mu, \, \dt y^i_{|u}, \, \Xi^a\w) \,,
\\
(x^\mu, \, y^i, \, z^a) \com D\Xi_u = 
\w(x^\mu(b), \, y^i, \, s^*(\eps)^a\w) \,,
\qquad
(\dt x^\mu, \, \dt y^i, \, \dt z^a) \com D\Xi_u = 
\w(0, \, \dt y^i, \, (D\Xi)^a\w) \,,
\end{gather*}
where
\bgt
u^\mu = (\der_0 c^\mu)(\lam) \,,
\\
\Xi^a = \Xi^a_0 + \der_i\w(s^*(\eps)^a\w) \, \dt y^i_{|u} \,,
\qquad
(D\Xi)^a = \der_i\w(s^*(\eps)^a\w) \, \dt y^i_{|0} \,,
\end{gather*}
with
\beq
u^\mu \in \Rn \,,
\qquad
\Xi^a_0 \in \map(\f F_b, \Rn) \,.
\eeq

 Thus, in this way, we obtain a natural map
\beq
\E r_s : \E X_s \in \E T_s\f S \mto (u, \Xi_u) \,,
\eeq
where
$u \in T_b\f B$
and
$\Xi_u : (T\f F)_u \to (T\f G)_u$
is a smooth map as above.

 For each
$s \in \f S \,,$
the map
$\E r_s$
turns out to be injective.\END
\eTh

\myskip

\bRm
 With reference to the above
Theorem \ref{Theorem: representation of tangent space of space of F--smooth system of sections},
we stress that a vector
$\E X_s \byd (u, \Xi_u) \in \E T_s\f S$
is characterised, in coordinates, by its components
\beq
u^\mu \in \Rn
\ssep{and}
\Xi^a_0 \in \map(\f F_b, \Rn) \,.
\eeq

 Moreover, the equality
\beq
D\Xi_u = T_b\w(s^*(\eps)\w) : T(\f F_b) \to T(\f G_b) \,.
\eeq
shows that
$D\Xi_u$
depends only on the element
$s \in \f S_b \,.$\END
\eRm

\myskip

 The set
$\E T\f S$
turns out to be equipped with the natural maps
\beq
\tau_\f S : \E T\f S \to \f S \,,
\qquad
\E T\zet : \E T\f S \to T\f B \,,
\qquad
\E T\eps : \E T\f S \ucar{T\f B} T\f F \to T\f G \,.
\eeq

\bPr
\label{Proposition: tau S, T zet, T eps for an F--smooth system of sections}
 We obtain in a natural way the following maps:
 
 1) the natural surjective map
(see
Theorem \ref{Theorem: representation of tangent space of space of F--smooth system of sections})
\beq
\tau_\f S : \E T\f S \to \f S : \E X_s \mto s \,,
\eeq

 2) the natural surjective map
\beq
\E T\zet : \E T\f S \to T\f B \,,
\eeq
given, according to
Theorem \ref{Theorem: representation of tangent space of space of F--smooth system of sections},
by
\beq
\E T\zet : \E X_s \byd (u, \Xi_u) \mto u \,,
\ssep{for each}
s \in \f S_b \,,
\sep{with}
b \in \f B \,,
\eeq

 3) the natural fibred map
\beq
\E T\eps : \E T\f S \ucar{T\f B} T\f F \to T\f G \,,
\eeq
given, according to
Theorem \ref{Theorem: representation of tangent space of space of F--smooth system of sections},
by
\bml
\E T\eps : (\E X_s, Y) \byd \w((u, \Xi_u), Y\w) \mto \Xi_u(Y) \,,
\\
\ssep{for each}
s \in \f S_b \,,
\quad
Y \in (T\f F)_u \,,
\sep{with}
u \in T_b\f B \,.
\end{multline*}

 Indeed, the following natural diagrams commute
\newdiagramgrid{2x2}
{1.2, 1.2, .7}
{.5, .7}
\bdg[grid=2x2]
\E T\f S &\rTo^{\E T\zet} &T\f B
\\
\dTo^{\tau_\f S} &&\dTo_{\tau_\f B}
\\
\f S &\rTo^{\zet} &\;\f B &,
\edg
\newdiagramgrid{3x2}
{1.5, 1.5, .9}
{.8, .8, .8, .8}
\bdg[grid=3x2]
T\f F
&\rTo^{\id_{T\f F}}
&T\f F
\\
\uTo^{\pro_2} &&\uTo_{Tq}
\\
\E T\f S \ucar{T\f B} T\f F
&\rTo^{\E T\eps} &T\f G
\\
\dTo^{\tau_\f S \car \tau_\f F} &&\dTo_{\tau_\f G}
\\
\f S \ucar{\f B} \f F
&\rTo^{\eps}
&\quad\f G &.\END
\edg
\ePr

\myskip

 Moreover, the set
$\E T\f S$
turns out to be equipped with the natural subset
\beq
\E V\f S \sub \E T\f S \,.
\eeq

\myskip

\bPr
\label{Proposition: vertical space of an F--smooth system of smooth sections}
 The elements of the subset
\beq
\E V \f S \byd (T\zet)^{-1} (0) \sub \E T\f S
\eeq
can be regarded as the pairs of the type
\beq
\ul {\E X}_s = (0, \Xi_0) \in V_s \f S \,,
\ssep{for each}
s \in \f S_b \,,
\eeq
where
\beq
\Xi_0 : (V\f F)_b \to (V\f G)_b \,,
\eeq
a) is a global smooth section of the smooth fibred manifold
$Tq : (V\f G)_b \to (V\f F)_b$
and an affine fibred morphisms over
$s^*(\eps) : \f F_b \to \f G_b \,,$ 
according to the following diagram commutative
\newdiagramgrid{3x2}
{1.4, 1.4, .7}
{.7, .7, .7, .7}
\bdg[grid=3x2]
(V\f F)_b &\rTo^{\id_{(V\f F)_b}} &(V\f F)_b
\\
\uTo^{\id_{(V\f F)_b}} &&\uTo_{Tq}
\\
(V\f F)_b &\rTo^{\Xi_0} &(V\f G)_b
\\
\dTo^{\tau_\f F} &&\dTo_{\tau_\f G}
\\
\f F_b &\rTo^{s^*(\eps)} &\;\;\f G_b &,
\edg
b) whose smooth fibre derivative
\beq
D\Xi_u : (V\f F)_b \to (V\f G)_b
\eeq
fulfills the equality
\beq
D\Xi_u = T_b\w(s^*(\eps)\w) : (V\f F)_b \to (V\f G)_b \,.
\eeq

 We have the coordinate expression
\beq
(x^\mu, \, y^i, \, z^a) \com \Xi_0 = 
\w(x^\mu(b), \, y^i, \, s^*_b(\eps)^a\w)
\ssep{and}
(\dt x^\mu_{|0}, \, \dt y^i_{|0}, \, \dt z^a_{|0}) \com \Xi_0 = 
\w(0, \, \dt y^i_{|0}, \, \Xi^a\w) \,,
\eeq
where
\beq
\Xi^a = \Xi^a_0 + \der_i\w(s^*(\eps)^a\w) \, \dt y^i_{|0} \,,
\ssep{with}
\Xi^a_0 \in \map(\f F_b, \Rn) \,.
\eeq
\ePr

\bpf
 The Corollary follows easily from
Theorem \ref{Theorem: representation of tangent space of space of F--smooth system of sections}
and
Proposition \ref{Proposition: tau S, T zet, T eps for an F--smooth system of sections} item 1).

 An alternative direct proof could be obtained by rephrasing the proof of
Theorem \ref{Theorem: representation of tangent space of space of F--smooth system of sections},
taking into account vertical F--smooth curves of
$\f S \,.$\QED
\epf

\myskip

 We can achieve two important simplifications in the representation of vertical elements of
$\E T\f S \,.$

 In fact, we can represent such an element
 
 1) as defined on
$\f F \,,$
instead of
$V\f F \,,$

 2) as valued in
$V_\f F \f G \,,$
instead of
$V_\f B \f G \,.$

\myskip

\bCr
\label{Corollary: alternative representation of the vertical space of an F--smooth system of smooth sections}
 We can equivalently regard the elements of the vertical subset
\beq
\E V \f S \sub \E T\f S
\eeq
as the pairs of the type
\beq
\ul{\E X}_s = (0, \ul\Xi_0) \in \E V_s \f S \,,
\ssep{for each}
s \in \f S_b \,,
\eeq
where
\beq
\ul\Xi_0 \byd \Xi_0 \com 0_\f F : \f F_b \to (V_\f F\f G)_b \,,
\ssep{with}
0_\f F : \f F_b \to (V\f F)_b \,,
\eeq
is a global section of the smooth fibred manifold
$\tau_\f F \com Tq : (V\f G)_b \to \f F_b$
according to the following diagram commutative
\newdiagramgrid{2x2}
{1.2, 1.2, .7}
{.7, .7, .7, .7}
\bdg[grid=2x2]
\f F_b &\lTo^{q} &(\f G)_b
\\
\uTo^{\id_{\f F_b}} &&\uTo_{\tau_\f G}
\\
\f F_b &\rTo^{\ul\Xi_0} &(V_\f F\f G)_b
\\
\dTo^{0} &&\dInto_{\cap}
\\
(V\f F)_b &\rTo^{\Xi_0} &\;(V\f G)_b &.
\edg

 We have the coordinate expression
\beq
(x^\mu, \, y^i, \, z^a) \com \ul\Xi_0 = 
\w(x^\mu(b), \, y^i, \, s^*_b(\eps)^a\w)
\ssep{and}
(\dt x^\mu_{|0}, \, \dt y^i_{|0}, \, \dt z^a_{|0}) \com \ul\Xi_0 = 
\w(0, \, 0, \, \Xi^a\w) \,,
\eeq
where
\beq
\ul\Xi^a = \Xi^a_0 \,,
\ssep{with}
\Xi^a_0 \in \map(\f F_b, \Rn) \,.
\eeq
\eCr

\bpf
 It follows easily from the above
Proposition \ref{Proposition: vertical space of an F--smooth system of smooth sections}.\QED
\epf

\bCr
\label{Corollary: natural splitting of S if G to F is a vector bundle and the system is injective}
 Let us suppose that
$q : \f G \to \f F$
be a vector bundle and consider the natural fibred isomorphism over
$\f F$
\beq
V_\f F\f G \seq \f G \ucar{\f F} \f G \,.
\eeq

 Then, the elemets
\beq
\ul{\E X}_s = (0, \ul\Xi_0) \in \E V_s \f S \,,
\ssep{for each}
s \in \f S_b \,,
\eeq
are characterised by the global smooth sections
\beq
\ul\Xi : \f F_b \to \f G_b \,,
\eeq
according to the following commutative diagram
\newdiagramgrid{2x2}
{1.3, 1.3, .7}
{.7, .7}
\bdg[grid=2x2]
\f F_b &\rTo^{\ul\Xi_0} &(V_\f F \f G)_b
\\
\dTo^{\id_{\f F_b}} &&\dTo_{\pro_2}
\\
\f F_b &\rTo^{\ul\Xi} &\;\f G_b &.
\edg

 Hence, if the F--smooth system
$(\f S, \zet, \eps)$
is injective, then we obtain a natural fibred isomorphism over
$\f B$
\beq
\E V\f S \seq \f S \ucar{\f B} \f S \,. 
\eeq
\eCr

\bpf
 The 1st claim follows directly from the above
Corollary \ref{Corollary: alternative representation of the vertical space of an F--smooth system of smooth sections}.

 Further, the splitting of
$V\f S$
follows from the fact that the smooth sections
$\ul\Xi : \f F_b \to \f G_b$
are just the selected sections of the system
$(\f S, \zet, \eps) \,.$\QED
\epf

\myskip

 The set
$\E T\f S$
and the associated maps
\beq
\tau_\f S : \E T\f S \to \f S \,,
\qquad
\E T\zet : \E T\f S \to T\f B
\qquad
\E T\eps : \E T\f S \ucar{T\f B} T\f F \to T\f G
\eeq
turn out to be F--smooth in a natural way, according to the following
Theorem.

\myskip

\bTh
\label{Theorem: the tangent space of the F--smooth system of smooth sections is F--smooth}
 Let us consider the set
$\E T\C C$
consisting of all curves
\beq
d\wha c : \f I_\wha c \to \E T\f S \,,
\eeq
given, for each
$\wha c \in \C C \,,$ 
according to
Theorem \ref{Theorem: representation of tangent space of space of F--smooth system of sections},
by
\beq
d\wha c : \lam \mto \E X_{\wha c(\lam)} \,,
\eeq

 Indeed, the set
$T\C C$
equips the set
$\E T\f S$
with an F--smooth structure.

 Moreover, the maps
\beq
\tau_\f S : \E T\f S \to \f S \,,
\qquad
T\zet : \E T\f S \to T\f B
\qquad
T\eps : \E T\f S \ucar{T\f B} T\f F \to T\f G
\eeq
turn out to be F--smooth.

 Thus, the 3--plet
$(\E T\f S, \E T\zet, \E T\eps)$
turns out to be an F--smooth system of fibrewisely sections of the smooth double fibred manifold
\newdiagramgrid{1}
{1, 1, 1, 1, .7}
{1}
\bdg[grid=1]
T\f G 
&\rTo^{Tq} &T\f F
&\rTo^{Tp} &T\f B &.\END
\edg
\eTh

\myskip

 The fibred set
$\tau_\f S : \E T\f S \to \f S$
inherits a vector structure in a natural way.
 
 Actually, the proof of this result is not straightforward because we have to achieve invariant linear algebraic operations on bundles which are not vector bundles.

\myskip

\bTh
\label{Theorem: vector structure of TS for an F--smooth system of sections}
 The fibres of the F--smooth fibred set
$\tau_\f S : \E T\f S \to \f S$
inherit a natural vector structure given, for each
$s \in \f S_b \,,$
with
$b \in \f B \,,$
by means of the equalities
\beq
r \, (u, \Xi_u) = (r \, u, r \vp \Xi_u)
\ssep{and}
(u, \Xi_u) \vs (\ac u, \ac\Xi_u) \byd 
\w((u + \ac u), (\Xi_u \vs \ac \Xi_u)\w) \,,
\eeq
where the smooth maps
\beq
r \vp \Xi_u : (T\f F)_{r \, u} \to (T\f G)_{r \, u}
\ssep{and}
\Xi_u \vs \ac\Xi_{\ac u} : (T\f F)_{u + \ac u} \to (T\f G)_{u + \ac u}
\eeq
are defined through the equivariant coordinate equalities
\bgt
(r \vp \Xi)^\mu \byd r \, \Xi^\mu \,,
\qquad
(r \vp \Xi)^i \byd r \, \Xi^i \,,
\qquad
(r \vp \Xi)^a_0 \byd r \, \Xi^a_0 \,,
\\
(\Xi \vs \ac\Xi)^\mu \byd \Xi^\mu + \ac\Xi^\mu \,,
\qquad
(\Xi \vs \ac\Xi)^i \byd \Xi^i + \ac\Xi^i \,,
\qquad
(\Xi \vs \ac\Xi)^a_0 \byd \Xi^a_0 + \ac\Xi^a_0 \,,
\end{gather*}
with
\beq
(r \vp \Xi)^i \byd \dt y^i_{|r u}
\ssep{and}
(\Xi \vs \ac\Xi)^i \byd \dt y^i_{u + \ac u} \,.
\eeq
\eTh

\bpf
 Let us consider the smooth maps
\beq
X : (T\f F)_u \to (T\f G)_u
\ssep{and}
Y : (T\f F)_u \to (T\f G)_u \,,
\eeq
with coordinate expressions
\bat{3}
X^\mu &= u^\mu \,,
\qquad
X^i &&= \dt y^i_{|u} \,,
\qquad
X^a
&&=
X^a_0 + \phi^a_i \, \dt y^i_{|u} \,,
\\
Y^\mu &= v^\mu \,,
\qquad
Y^i &&= \dt y^i_{|v} \,,
\qquad
Y^a
&&=
Y^a_0 + \phi^a_i \, \dt y^i_{|v} \,,
\end{alignat*}
and
\beq
X^a_0,\, Y^a_0, \phi^a_i \in \map(\f F, \Rn) \,.
\eeq

 Then, we consider their local smooth extensions to the total tangent space
$T\f F \,,$
which are induced by the chosen smooth fibred chart
$(x^\lam, y^i, z^a) \,,$
\beq
\wti X : T\f F \to T\f G
\ssep{and}
\wti Y : T\f F \to T\f G \,,
\eeq
with coordinate expressions
\bat{3}
\wti X^\mu &= \dt x^\mu \,,
\qquad
\wti X^i &&= \dt y^i \,,
\qquad
\wti X^a
&&=
X^a_0 + \phi^a_i \, \dt y^i \,,
\\
\wti Y^\mu &= \dt x^\mu \,,
\qquad
\wti Y^i &&= \dt y^i \,,
\qquad
\wti Y^a
&&=
Y^a_0 + \phi^a_i \, \dt y^i \,.
\end{alignat*}

 Next, let us take into account the vector structure of the smooth bundle
$T\f G \to \f G \,,$
by which we can define the algebraic operations
\beq
r \, \wti X : T\f F \to T\f G
\ssep{and}
\wti X + \wti Y : T\f F \to T\f G \,,
\eeq
with coordinate expressions
\bgt
r \, \wti X^\mu = r \, \dt x^\mu \,,
\qquad
r \,\wti X^i = r  \,\dt y^i \,,
\qquad
r \, \wti X^a = r \, X^a_0 + \phi^a_i \, r \, \dt y^i \,,
\\
(\wti X + \wti Y)^\mu = \dt x^\mu + \dt x^\mu \,,
\qquad
(\wti X + \wti Y)^i = \dt y^i + \dt y^i \,,
\qquad
(\wti X + \wti Y)^a =
X^a_0 + Y^a_0 + \phi^a_i \, (\dt y^i + \dt y^i) \,.
\end{gather*}

 Moreover, by restricting the above maps to the zero section
$0 \sub T\f F \,,$
we obtain the smooth maps 
\beq
(r \, \wti X) (0) : \f F \to T\f G
\ssep{and}
(\wti X + \wti Y) (0) : \f F \to T\f G \,,
\eeq
with coordinate expressions
\bgt
r \, \wti X^\mu(0) = 0 \,,
\qquad
r \,\wti X^i(0) = 0 \,,
\qquad
r \, \wti X^a(0) = r \, X^a_0 \,,
\\
(\wti X + \wti Y)^\mu(0) = 0 \,,
\qquad
(\wti X + \wti Y)^i(0) = 0 \,,
\qquad
(\wti X + \wti Y)^a(0) =
X^a_0 + Y^a_0 \,.
\end{gather*}

 Furthermore, by taking into account the affine structure of the smooth bundle
$Tq : T\f G \to T\f F \,,$
we observe that the maps
$\wti X$
and
$\wti Y$
are affine and their derivatives are the maps
\beq
D \wti X = D \wti Y = T\eps_s : T\f F \to V_\f F\f G \,,
\eeq
with coordinate expression
\beq
(D \wti X)^\mu = (D \wti Y)^\mu = \dt x^\mu \,,
\qquad
(D \wti X)^i = (D \wti Y)^i = \dt y^i \,,
\qquad
(D \wti X)^a = (D \wti Y)^a = \phi^a_i \, \dt y^i \,.
\qquad
\eeq

 Then, the following smooth maps are well defined
\bal
r \vp \wti X 
&: 
T\f F \to T\f G : 
v \mto (r \, \wti X)(0) + D \wti X(v) \,,
\\
\wti X \vs \wti Y 
&: 
T\f F \to T\f G : 
v \mto (\wti X + \wti Y)(0) + D \wti X(v)
\end{align*}
and have coordinate expressions
\bgt
r \vp \wti X^\mu = \dt x^\mu \,,
\qquad
r \vp \wti X^i = \dt y^i \,,
\qquad
r \vp \wti X^a = r \, X^a_0 + \phi^a_i \, \dt y^i \,,
\\
(\wti X \vs \wti Y)^\mu = \dt x^\mu \,,
\qquad
(\wti X \vs \wti Y)^i = \dt y^i \,,
\qquad
(\wti X \vs \wti Y)^a =
X^a_0 + Y^a_0 + \phi^a_i \, \dt y^i \,.
\end{gather*}

 Further, we consider the smooth restrictions of the above smooth maps
\beq
r \vp \wti X : T\f F \to T\f G
\ssep{and}
\wti X \vs \wti Y : T\f F \to T\f G
\eeq
to the smooth subbundles
$(T\f F)_{r u}$
and
$(T\f F)_{u + v} \,,$
respectively,
\beq
(r \vp \wti X)_{r u} : (T\f F)_{r u} \to T\f G
\ssep{and}
(\wti X \vs \wti Y)_{u+v} : (T\f F)_{u+v} \to T\f G \,.
\eeq 

 We have the coordinate expressions
\bgt
(r \vp \wti X)^\mu_{r u} = r \, u^\mu \,,
\qquad
(r \vp \wti X)^i_{r u} = \dt y^i \,,
\qquad
(r \vp \wti X)^a_{r u} = r \, X^a_0 + \phi^a_i \, \dt y^i \,,
\\
(\wti X \vs \wti Y)^\mu_{u+v} = r \, u^\mu \,,
\qquad
(\wti X \vs \wti Y)^i_{u+v} = \dt y^i \,,
\qquad
(\wti X \vs \wti Y)^a_{u+v} =
X^a_0 + Y^a_0 + \phi^a_i \, \dt y^i \,.
\end{gather*}

 The above coordinate expressions show that the maps
\beq
(r \vp \wti X)_{r u} : (T\f F)_{r u} \to T\f G
\ssep{and}
(\wti X \vs \wti Y)_{u+v} : (T\f F)_{u+v} \to T\f G
\eeq
factorise through maps (denoted by the same symbols)
\beq
(r \vp \wti X)_{r u} : (T\f F)_{r u} \to (T\f G)_{r u}
\ssep{and}
(\wti X \vs \wti Y)_{u+v} : (T\f F)_{u+v} \to (T\f G)_{u+v} \,,
\eeq 
according to the following commutative diagrams
\newdiagramgrid{2x2-2x2}
{2, 2, 3, 2, 2, .7}
{.7, .7}
\bdg[grid=2x2-2x2]
(T\f F)_{r u} &\rTo^{(r \vp \wti X)_{r u}} &(T\f G)_{r u} 
&(T\f F)_{u+v} &\rTo^{(\wti X \vs \wti Y)_{u+v}} &(T\f G)_{u+v}
\\
\dTo^{\id} &&\dInto_{\sub}
&\dTo^{\id} &&\dInto_{\sub}
\\
(T\f F)_{r u} &\rTo^{(r \vp \wti X)_{r u}} &T\f G
&(T\f F)_{u+v} &\rTo^{(\wti X \vs \wti Y)_{u+v}} &T\f G &.
\edg

 Moreover, the above coordinate expressions show that the maps
\beq
(r \vp \wti X)_{r u} : (T\f F)_{r u} \to (T\f G)_{ru}
\ssep{and}
(\wti X \vs \wti Y)_{u+v} : (T\f F)_{u+v} \to (T\f G)_{u+v}
\eeq
do not depend on the extensions
$\wti X$
and
$\wti Y$
induced by the chart, but depend only by the original maps
$X$
and
$Y \,.$

 For this reason, we can write
\bat{2}
(r \vp X)_{r u} 
&\byd 
(r \vp \wti X)_{r u} 
&&: 
(T\f F)_{r u} \to (T\f G)_{ru} \,,
\\
(X \vs Y)_{u+v} 
&\byd 
(\wti X \vs \wti Y)_{u+v} 
&&: 
(T\f F)_{u+v} \to (T\f G)_{u+v}
\end{alignat*}

 Hence, the above definition of the maps
$(r \vp X)_{r u}$
and
$(X \vs Y)_{u+v}$ 
is coordinate free.\QED
\epf

\myskip

\bRm 
\label{Remark: the algebraic operations are coordinate free}
 The fact that the algebraic operations defined in the above 
Theorem \ref{Theorem: vector structure of TS for an F--smooth system of sections}
be coordinate free can be confirmed by the following explicit check.

 Let
$(\ac x^\mu, \ac y^i, \ac z^a)$
be another fibred chart of
$\f G \,.$

 Then, we have the following transition formulas
\bat{3}
\ac X^i
&=
\ac\der_\mu y^i \, u^\mu + \ac\der_j y^i \, X^j \,,
&&\qquad
\ac X^a
&&=
\ac\der_\mu z^a \, u^\mu + \ac\der_j z^a \, X^j + 
\ac\der_b z^a \, (X^b_0 + \phi^b_j \, \dt y^j_{|u}) \,,
\\
\ac Y^i
&=
\ac\der_\mu y^i \, v^\mu + \ac\der_j y^i \, Y^j \,,
&&\qquad
\ac Y^a
&&=
\ac\der_\mu z^a \, v^\mu + \ac\der_j z^a \, Y^j + 
\ac\der_b z^a \, (Y^b_0 + \phi^b_j \, \dt y^j_{|v}) \,,
\end{alignat*}
which yield also the following equalities
\bat{3}
\ac X^i
&=
\ac\der_j y^i \, X^j + \ac\der_\mu y^i \, u^\mu \,,
&&\qquad
\ac X^a_0
&&=
\ac\der_b z^a \, X^b_0 + \ac\der_\mu z^a \, u^\mu \,,
\\
\ac Y^i
&=
\ac\der_j y^i \, Y^j + \ac\der_\mu y^i \, v^\mu \,,
&&\qquad
\ac Y^a_0
&&=
\ac\der_b z^a \, Y^b_0 + \ac\der_\mu z^a \, v^\mu \,.
\end{alignat*}

 Hence, we obtain
\bal
r \, \ac X^i 
&=
\ac\der_\mu y^i \, (r \, u^\mu) + \ac\der_j y^i \, (r \, X^j) \,,
\\
r \, \ac X^a_0 
&=
\ac\der_b z^a \, (r \, X^b_0) + \ac\der_\mu z^a \, (r \, u^\mu) \,,
\end{align*}
and
\bal
\ac X^i + \ac Y^i 
&=
\ac\der_j y^i \, (X^j + Y^j) + \ac\der_\mu y^i \, (u^\mu + v^\mu) \,,
\\
\ac X^a_0 + \ac Y^a_0 
&=
\ac\der_b z^a \, (X^b_0 + Y^b_0) + \ac\der_\mu z^a \, (u^\mu + v^\mu) \,.
\end{align*}

 Moreover, we have the following transition formulas
\beq
\dtac y^i_{|u} = 
\ac\der_\mu y^i \, u^\mu + \ac\der_j \ac y^i \, \dt y^j_{|u}
\ssep{and}
\dtac y^i_{|v} = 
\ac\der_\mu y^i \, v^\mu + \ac\der_j \ac y^i \, \dt y^j_{|v}
\eeq
and
\beq
\dtac y^i_{|ru} = 
\ac\der_\mu y^i \, (r \, u^\mu) + \ac\der_j \ac y^i \, \dt y^j_{|ru}
\ssep{and}
\dtac y^i_{|u + v} = 
\ac\der_\mu y^i \, (u^\mu + v^\mu) + 
\ac\der_j \ac y^i \, \dt y^j_{|u + v} \,.
\eeq

 Hence, if in one chart we have
\beq
(r \, X)^i \byd \dt y^i_{|r u}
\ssep{and}
(X + Y)^i \byd \dt y^i_{|u+v} \,,
\eeq
then analogous formulas hold in the other chart.\END
\eRm

\bCr
 The fibred subset over
$\f S$
\beq
\E V\f S \sub \E T\f S
\eeq
turns out to be a vector fibred subset.\END
\eCr

\bCr
  The fibres of the F--smooth fibred set
$\E T\zet : \E T\f S \to T\f B$
inherit a natural affine structure, whose associated vector spaces are the fibres of
$\E V\f S \,.$\END
\eCr

\myskip

 We can define the F--smooth tangent prolongation 
$T\sig : T\f B \to \E T\f S$
of F--smooth sections
$\sig \in \Fsec(\f B, \f S)$
analogously to the case of smooth systems of sections
(see
\S\ref{Systems of sections}).

\myskip

\bDf
\label{Definition: tangent prolongation of sections}
 We define the tangent prolongation of an F--smooth section
\beq
\sig \in \Fsec(\f B, \f S)
\eeq
to be the tubelike F--smooth section
\beq
\E T \sig : T\f B \to \E T \f S \,,
\eeq
given, according to
Theorem \ref{Theorem: representation of tangent space of space of F--smooth system of sections},
by
\beq
\E T \sig : u \in T_b\f B \mto 
\E T_u\sig \byd \E X_{\sig(b)} \byd (u, \Xi_u) \in \E T_u \f S \,, 
\eeq
where
\beq
\Xi_u = T_u \w(\sig^*(b)(\eps)\w) : 
(T\f F)_u \to (T\f G)_u \,.\END
\eeq
\eDf

\bPr
 For each F--smooth section
$\sig \in \Fsec(\f B, \f S )\,,$
the following diagram commutes
\newdiagramgrid{3x2}
{1.2, 1.2, .9}
{0.7, 0.7, 0.7, 0.7}
\bdg[grid=3x2]
\f B
&\rTo^{\sig}
&\f S
\\
\uTo^{\tau_\f B} &&\uTo_{\tau_\f S}
\\
T\f B
&\rTo^{\E T\sig} &\E T\f S
\\
\dTo^{\id} &&\dTo_{T\zet}
\\
T\f B
&\rTo^{\id_\f B}
&\;T\f B &.\END
\edg
\ePr

\bNt
 By a certain mild abuse of language, we can write the equality
\beq
\E T\sig = (\E T\sig)^*(\E T\eps) \,.\END
\eeq
\eNt
%--------------------------------------------------------------------%
\subsection{F--smooth differential operators}
\label{F--smooth differential operators}
%--------------------------------------------------------------------%
\bsm
 In view of a discussion on F--smooth connections of an F--smooth system of fibrewisely smooth sections
(see the forthcoming
Section \ref{F--smooth connections}),
we analyse the smooth operators.

\myskip

 Given two injective F--smooth systems 
$(\f S, \zet, \eps)$
and
$(\acf S, \ac\zet, \ac\eps)$
of fibrewisely smooth sections of the smooth double fibred manifolds
$\f G \oset{q}\to \f F \oset{p}\to \f B$
and
$\acf G \oset{\ac q}\to \f F \oset{p}\to \f B$
a sheaf morphism
\beq
\C D : \tub(\f F, \f G) \to \tub(\f F, \acf G) \,,
\eeq
which is compatible with the above F--smooth systems, yields in a natural way a sheaf morphism
\beq
\whaC D : \Fsec(\f B, \f S) \to \Fsec(\f B, \acf S) \,.
\eeq
\esm

 Thus, let us consider two smooth double fibred manifolds
\newdiagramgrid{1}
{.7, .7, .8, .7, 2, 2, .7, .7, .7, .7}
{0}
\bdg[grid=1]
\f G 
&\rTo^{q} &\f F
&\rTo^{p} & \f B
&\text{and}
&\f G' 
&\rTo^{\ac q} &\f F
&\rTo^{p} & \f B \,.
\edg
and denote the smooth fibred charts of
$\f G$
and
$\acf G \,,$
respectively, by
\beq
(x^\mu, y^i, z^a)
\ssep{and}
(x^\mu, y^i, \ac z^a) \,.
\eeq

 Moreover, let us consider two \emp{injective} F--smooth systems of fibrewisely smooth sections of the two smooth fibred manifolds above
\beq
(\f S, \zet, \eps)
\ssep{and}
(\acf S, \ac\zet, \ac\eps) \,.
\eeq

\myskip

\bDf
\label{Definition: compatible smooth tubelike operators}
\index{F--smooth systems!compatible smooth tubelike operator}
 A sheaf morphism
\beq
\C D : \tub(\f F, \f G) \to \tub(\f F, \acf G)
\eeq
is said to be \emp{compatible} with the F--smooth systems of smooth sections 
$(\f S, \zet, \eps)$
and
$(\acf S, \ac\zet, \ac\eps)$
if it restricts to a sheaf morphism
\beq
\C D_\f S : \tub_\f S(\f F, \f G) \to \tub_{\acf S}(\f F, \acf G) \,,
\eeq
according to the following commutative diagram
\newdiagramgrid{2x2}
{1.9, 1.9, 1.5, .5}
{.7, .7}
\bdg[grid=2x2]
\tub(\f F, \f G) &\rTo^{\C D} &\tub(\f F, \acf G)
\\
\uInto &&\uInto
\\
\tub_\f S(\f F, \f G) &\rTo^{\C D_\f S} 
&\tub_{\acf S}(\f F, \acf G) &.\END
\edg
\eDf

\bPr
\label{Proposition: F--smooth operator}
 Let us consider a sheaf morphism of compatible smooth operators
\beq
\C D : \tub(\f F, \f G) \to \tub(\f F, \acf G) \,.
\eeq

 Then, in virtue of
Theorem \ref{Theorem: sections br sig are smooth if and only if sig are F--smooth},
we obtain the sheaf morphism
\beq
\whaC D : \Fsec(\f B, \f S) \to \Fsec(\f B, \acf S) : 
\wha \sig \mto \wha{\C D (\sig)} \,.\END
\eeq
\ePr
%--------------------------------------------------------------------%
%\newpage
\markboth{\rm Chapter \thechapter. Systems of sections}{\rm \thesection. F--smooth vs smooth systems of sections}
\section{F--smooth vs smooth systems of sections}
\label{F--smooth vs smooth systems of sections}
\markboth{\rm Chapter \thechapter. Systems of sections}{\rm \thesection. F--smooth vs smooth systems of sections}
%--------------------------------------------------------------------%
\bsm
 In the above 
Section \ref{F--smooth tangent prolongation of (S, zet, eps)},
we have discussed the F--smooth tangent space
$\E T\f S$
of each F--smooth system
$\f S$
of fibrewisely smooth sections
$\phi : \f F \to \f G$
of a smooth double fibred manifold
$\f G \oset{q}{\to} \f F \oset{p}{\to} \f B$
(see
Definition \ref{Definition: tangent space of the F--smooth space of smooth tubelike sections})
and studied its main properties
(see
Theorem \ref{Theorem: representation of tangent space of space of F--smooth system of sections}
and
Theorem \ref{Theorem: the tangent space of the F--smooth system of smooth sections is F--smooth}).
 
 On the other hand, if the set
$\f S$
is assumed a priory to be a finite dimensional smooth manifold and
$\eps$
to be a smooth map, then we can achieve the tangent space
$T\f S$
directly in terms of the standard differential geometry of smooth manifolds.

 Thus, the need of a comparison between the F--smooth approach to
$\E T\f S$
and the smooth approach to
$T\f S$
arises naturally, having in mind
Theorem \ref{Theorem: relation between smoothness and F--smoothness}.

 Actually, by regarding a smooth system
$(\f S, \zet, \eps)$
as a particular F--smooth system, we obtain a natural map
\beq
\imath : T\f S \to \E T\f S : X \mto \E X
\eeq
and can prove that the representation of
$\E X$
in terms of the smooth map
\beq
\Xi_u \byd \E r(\E X_u) : (T\f F)_u \to (T\f G)_u
\eeq
turns out to be given by the equality
$\Xi_u = T_1\eps_{X_u} \,.$

\myskip

 We leave to the reader the task to develop in detail the above considerations, by rephrasing to the present context of F--smooth systems of fibrewisely smooth sections the considerations that we have discussed in
Section \ref{F--smooth vs smooth systems of maps}, 
which is devoted to such a comparison in the context of systems of maps.

\myskip

 Here, in order to clarify the ``odd terms" appearing in the representation of F--smooth tangent space of an F--smooth system of fibrewisely smooth sections, we just discuss the above items in terms of standard smooth manifolds in the particular case of a finite dimensional smooth system of linear sections.
\esm

\bEx
\label{Example: tangent space of the F--smooth system of linear sections}
 Let us consider two vector bundles
\beq
p : \f F \to \f B
\ssep{and}
p \com q : \f G \to \f B
\eeq
and the system
$(\f S, \zet, \eps)$
of linear sections
$\br\sig \in \tub(\f F, \f G) \,,$
which has been discussed in
Example \ref{Example: system of linear sections},
in the framework of smooth systems of smooth sections.

 In this case,
$\f S$
has been assumed to be a priori a smooth manifold.
 Hence, we can avail of this smooth structure to introduce and analyse the tangent space
$T\f S \,.$

 The smooth double fibred chart
$(x^\lam, y^i, z^a)$
of
$\f G$
induces naturally the smooth fibred charts
$(x^\lam, w^a_i)$
of
$\f S$
and
$(x^\lam, w^a_i; \, \dt x^\lam, \dt w^a_i)$
of
$T\f S \,.$

 Hence, we obtain the following coordinate expressions of the smooth maps
$\zet$
and
$\eps$
\beq
(x^\lam) \com \zet = (x^\lam)
\ssep{and}
(x^\lam, y^i, z^a) \com \eps = (x^\lam, y^i, w^a_j \, y^j) \,.
\eeq

 The coordinate expressions of the smooth maps
\beq
T\zet : T\f S \to T\f B
\ssep{and}
T\eps : T\f S \ucar{T\f B} T\f F \to T\f G
\eeq
are
\bgt
(x^\lam, \dt x^\lam) \com T\zet = (x^\lam, \dt x^\lam) \,,
\\
(x^\lam, y^i, z^a; \, \dt x^\lam, \dt y^i, \dt z^a) \com T\eps = 
(x^\lam, y^i, w^a_j \, y^j; \, 
\dt x^\lam, \dt y^i, \dt w^a_j \, y^j + w^a_j \, \dt y^j) \,.
\end{gather*}

 Thus, the smooth map
$T\eps$
yields a natural representation of tangent vectors of
$\f S$
via the smooth map, which, for each
$X_u \in (T_s\f S)_u \,,$
where
$s \in \f S_b \,,$
$u \in T_b\f B \,,$
$b \in \f B \,,$
yields the smooth section
\beq
\Xi_u \byd (T\eps)_{X_u} : (T\f F)_u \to (T\f G)_u \,,
\eeq
with coordinate expression
\beq
\Xi^\lam = X^\lam = u^\lam \,,
\qquad
\Xi^i = \dt y^i_{|u}\,,
\qquad
\Xi^a = X^a_j \, y^j_{|b} + w^a_j \, \dt y^j_{|u} =
\Xi^a_0 + \der_j \eps^a \, \dt y^j_{|u} \,,
\eeq
where
\beq
X^\lam \in \Rn \,,
\qquad
\Xi^a_0 \byd X^a_j \, y^j_{|b} \,.
\eeq

 Moreover, given
$r \in \Rn$
and
$X_u \in (T_s\f S)_u \,,$
$\ac X_v \in (T_s\f S)_v \,,$ 
the natural vector structure of the smooth vector bundle
$\tau_\f S : T\f S \to \f S$
yields the standard coordinate expressions
\bgt
(r \, X_u)^\lam = r \, u^\lam \,,
\qquad
(r \, X_u)^a_i = r \, X^a_i \,,
\\
(X_u + \ac X_v)^\lam = u^\lam + v^\lam \,,
\qquad
(X_u + \ac X_v)^a_i = X^a_i + \ac X^a_i \,,
\end{gather*}
and the zero vector
$X_0 \byd 0 \in T_s \f S$
has coordinate expression
\beq
X^\lam = 0 \,,
\qquad
X^a_i = 0 \,.
\eeq

 On the other hand, the representation of above formulas in terms of fibred morphisms turns out to read as follows
\bgt
(r \, \Xi_u)^\lam = r \, \Xi^\lam \,,
\qquad
(r \, \Xi_u)^i = r \, \Xi^i \,,
\\
(r \, \Xi_u)^a_i = 
r \, \Xi^a_j \, y^j_{|b} + w^a_j \, \dt y^j_{|ru} \,,
\\
(\Xi_u + \ac\Xi_v)^\lam = 
\Xi^\lam + \ac\Xi^\lam \,,
\qquad
(\Xi_u + \ac\Xi_v)^i = 
\Xi^i + \ac\Xi^i \,,
\\
(\Xi_u + \ac\Xi_v)^a_i = 
(\Xi^a_j + \ac\Xi^a_j) \, y^j_{|b} +  w^a_j \, \dt y^j_{|u+v} \,,
\end{gather*} 
and the representation of the zero vector turns out to be
\beq
\Xi^\lam = 0 \,,
\qquad
\Xi^i = \dt y^i_{|0} \,,
\qquad
\Xi^a_0 = w^a_j \, \dt y^j_{|0} \,.
\eeq

 Thus, in agreement with the general results found for F--smooth systems of sections
(see
Theorem \ref{Theorem: representation of tangent space of space of F--smooth system of sections}
and
Theorem \ref{Theorem: vector structure of TS for an F--smooth system of sections}),
the ``odd" term of the type
$w^a_j \, \dt y^j_{|u}$
appears in the representation of all vectors
$X_u$
and the zero vector is not represented by a vanishing map but by the map
$w^a_j \, \dt y^j_{|0} \,.$\END
\eEx
%--------------------------------------------------------------------%
\chapter{Systems of connections}
\label{Systems of connections}
%--------------------------------------------------------------------%
\bsm
 We start by discussing the \emp{smooth systems} 
$(\f C, \zet, \eps)$
of \emp{smooth connections}
\beq
c : \f F \to T^*\f B \ten T\f F
\eeq
of a smooth fibred manifold
$\f F \to \f B \,.$

 Indeed, such smooth systems can be regarded as be a particular case of smooth systems of smooth sections of the smooth double fibred manifold
$T^*\f B \ten T\f F \to \f F \to \f B \,.$

 Here, the ``space of parameters"
$\zet : \f C \to \f B$
is a \emp{smooth space} and the ``evaluation map"
$\eps : \f C \ucar{\f B} \f F \to T^*\f B \ten T\f F$
a \emp{smooth fibred morphism}
over 
$\f B \,.$

 Moreover, we discuss the \emp{smooth universal connection}
\beq
c\Upa : \f F\Upa \to T^*\f C \ten T\f F\Upa
\eeq
of a smooth system of smooth connections.

\myskip

 Then, we discuss the \emp{F--smooth systems}
$(\f C, \zet, \eps)$
of \emp{fibrewisely smooth connections}
\beq
c : \f F \to T^*\f B \ten T\f F
\eeq
of a smooth double fibred manifold
$\f G \to \f F \to \f B \,.$

 Indeed, such F--smooth systems of connections can be regarded as be a particular case of F--smooth systems of fibrewisely smooth sections of the smooth double fibred manifold
$T^*\f B \ten T\f F \to \f F \to \f B \,.$

 Here, the ``space of parameters"
$\zet : \f C \to \f B$
is an \emp{F--smooth space} and the ``evaluation map"
$\eps : \f C \ucar{\f B} \f F \to \f G$
an \emp{F--smooth fibred morphism}
over 
$\f B \,.$

\myskip
 
 The reader can find further discussions concerning the present subject in
\cite{Cab87,CabKol95,CabKol98,DodMod86,Gar72,Kol87,KolMicSlo93,
KolMod98,KriMic97,ManMod83d,MarMod91,Slo86}.
\esm
%--------------------------------------------------------------------%
%\newpage
\markboth{\rm Chapter \thechapter. Systems of connections}{\rm \thesection. Smooth systems of smooth connections}
\section{Smooth systems of smooth connections}
\label{Section: Smooth systems of smooth connections}
\markboth{\rm Chapter \thechapter. Systems of connections}{\rm \thesection. Smooth systems of smooth connections}
%--------------------------------------------------------------------%
\bsm
 We discuss the \emp{smooth systems of smooth connections}
\beq
c : \f F \to T^*\f B \ten T\f F
\eeq
of a smooth fibred manifold
$\f F \to \f B \,.$

 Moreover, we discuss the \emp{smooth universal connection}
\beq
c\Upa : \f F\Upa \to T^*\f C \ten T\f F\Upa
\eeq
of a smooth system of smooth connections and show that
$c\Upa$
characterises the system. 
\esm
%--------------------------------------------------------------------%
\subsection{Smooth systems of smooth connections}
\label{Smooth systems of smooth connections}
%--------------------------------------------------------------------%
\bsm
 Given a smooth fibred manifold
$p : \f F \to \f B \,,$
a ``\emp{smooth system of smooth connections}" is just a smooth system of smooth sections
(see
Definition \ref{Definition: systems of sections})
of the smooth double fibred manifold
\beq
\f G \byd T^*\f B \ten T\f F \to \f F \to \f B
\eeq
consisting of a selected family of smooth connections 
\beq
c : \f F \to T^*\f B \ten T\f F
\eeq
of the smooth fibred manifold
$p : \f F \to \f B \,.$

 Therefore, all developments discussed in the previous
Section \ref{Smooth systems of smooth sections}
can be applied to the present Section.
\esm

 Let us consider a smooth fibred manifold
$p : \f F \to \f B$
and denote its fibred charts by
$(x^\lam, y^i) \,.$

 Then, we define the fibred manifold
\beq
q : T^*\f B \ten T\f F \to \f F
\eeq
and obtain the smooth double fibred manifold
\newdiagramgrid{1}
{1.2, .7, .5}
{1}
\bdg[grid=1]
T^*\f B \ten T\f F 
&\rTo^{q} &\f F
&\rTo^{p} &\f B \,.
\edg

\myskip

 In the present context, it is convenient to deal with the definition of ``smooth connection" of the smooth fibred manifold
$p : \f F \to \f B$
as a smooth tangent valued form 
\beq
c : \f F \to T^*\f B \ten T\f F \,,
\eeq
which makes the following diagram commutative
\newdiagramgrid{2x2}
{1.2, 1.5, 1.2, 1.2, .7}
{.7, .7}
\bdg[grid=2x2]
\f F &\rTo^{c} &T^*\f B \ten T\f F
\\
\dTo^{p} &&\dTo_{\id_{T^*\f B} \ten Tp}
\\
\f B &\rTo^{\f 1_\f B} &T^*\f B \ten T\f B &.
\edg

 Then, the coordinate expression of
$c$
is of the type
\beq
c = d^\lam \ten (\der_\lam + c^i_\lam \, \der_i) \,,
\ssep{with}
c^i_\lam \in \map(\f F, \Rn) \,.
\eeq

 We denote the subsheaf of smooth tubelike connections of the fibred manifold
$\f F \to \f B$
by
(see
Definition \ref{Definition: tubelike sections})
\beq
\cns\tub(\f F, T^*\f B \ten T\f F) \sub
\tub(\f F, T^*\f B \ten T\f F) \,.
\eeq

\myskip

\bDf
\label{Definition: systems of connections}
\index{smooth systems!systems of connections}
\index{smooth systems!injective system of connections}
 A \emp{smooth system of smooth connections}
is defined to be a 3--plet
$(\f C, \zet, \eps)$
where

 1) $\zet : \f C \to \f B$
is a smooth fibred manifold,

 2) $\eps : \f C \ucar{\f B} \f F \to T^*\f B \ten T\f F$
is a smooth fibred morphism over
\beq
\id_\f F : \f F \to \f F
\ssep{and}
\f 1_\f B : \f B \to T^*\f B \ten T\f B \,,
\eeq
according to the following commutative diagrams
\newdiagramgrid{2x2}
{1.5, 1.5, .5, 3, 1.5, 1.5, 1.5}
{.7, .7}
\bdg[grid=2x2]
\f C \ucar{\f B} \f F &\rTo^{\eps} &T^*\f B \ten T\f F
&&\f C \ucar{\f B} \f F &\rTo^{\eps} &T^*\f B \ten T\f F
\\
\dTo^{\pro_2} &&\dTo_{\tau_\f F}
&&\dTo &&\dTo_{\id_{T^*\f B} \ten Tp}
\\
\f F &\rTo^{\id_\f F} &\;\f F &
&\f B &\rTo^{\1_\f B} &\quad T^*\f B \ten T\f B &.
\edg
 
 We call
$\eps$
the \emp{evaluation map} of the system.
\myskip

 Thus, the evaluation map
$\eps$
yields the sheaf morphism
\beq
\eps_\f C : \sec(\f B, \f C) \to \cns\tub(\f F, T^*\f B \ten T\f F) : 
\gam \mto c \byd \br\gam \,,
\eeq
where, for each
$\gam \in \sec(\f B, \f C) \,,$ 
the tubelike smooth connection
$c \byd \br\gam$
is defined by
\beq
\br\gam : \f F \to T^*\f B \ten \f F : 
f_b \mto \eps \w(\gam(b), f_b\w) \,,
\ssep{for each}
b \in \f B \,.
\eeq

 Therefore, the map
$\eps_\f C : \sec(\f B, \f C) \to \tub(\f F, T^*\f B \ten T\f F)$
provides a \emp{selection} of the tubelike smooth connections
$c : \f F \to T^*\f B \ten T\f F \,,$
given by the subset
\beq
\tub_\f C(\f F, T^*\f B \ten T\f F) \byd 
\eps_\f C\w(\sec(\f B, \f C)\w) \sub 
\cns\tub(\f F, T^*\f B \ten T\f F) \,.
\eeq

 We have the coordinate expressions
\bgt
\eps = d^\lam \ten (\der_\lam + \eps\co\lam i \, \der_i) \,,
\ssep{with}
\eps\co\lam i \in \map(\f F\Upa, \, \Rn) \,,
\\
\br\gam \byd \eps^*(\gam) = 
d^\lam \ten \w(\der_\lam + 
(\eps\co\lam i \com \gam_\f F) \, \der_i\w) \,,
\end{gather*}
where
(see
Definition \ref{Definition: lift fibred manifold})
\beq
\f F\Upa \byd \f C \ucar{\f B} \f F \to \f F 
\ssep{and}
\gam_\f F : \f F \to \f C \ucar{\f B} \f F : 
f_b \mto \w(\gam(b), f_b\w) \,.
\eeq

\myskip

 The smooth system of smooth connections
$(\f C, \zet, \eps)$
is said to be \emp{injective} if the map
$\eps_\f C : \sec(\f B, \f C) \to \cns\tub(\f F, T^*\f B \ten T\f F)$
is injective, i.e. if, for each
$\gam \,, \ac\gam \in \sec(\f B, \f C) \,,$
\beq
\br\gam \eqv \eps^* (\gam) = \brac\gam \eqv \eps^* (\ac\gam)
\Rarr
\gam = \ac\gam \,.
\eeq

 If the system is injective, then we obtain the bijection
\beq
\eps_\f C : 
\sec(\f B, \f C) \to \cns\tub_\f C(\f F, T^*\f B \ten T\f F) :
\gam \mto \br\gam \,,
\eeq
whose inverse is denoted by
\beq
(\eps_\f C)^{-1} : 
\cns\tub_\f C(\f F, T^*\f B \ten T\f F) \to \sec(\f B, \f C) :
c \mto \gam \eqv \wha c \,.\END
\eeq
\eDf

\myskip

 Let us examine a few distinguished examples of injective smooth systems of smooth connections.
 
 Indeed, in the case of smooth systems of linear connections, affine connections and principal connections, our bundle
$\zet : \f C \to \f B$
is just the standard bundle of coefficients of such connections.
 
\myskip

\bEx
\label{Example: system of linear connections}
\index{smooth systems!system of linear connections}
 If
$p : \f F \to \f B$
is a vector bundle, then the linear connections constitute an injective smooth system
$(\f C, \zet, \eps) \,,$
where
$\zet : \f C \to \f B$
is an affine subbundle
\beq
\f C \sub \lin_\f B(\f F, \, T^*\f B \ten T\f F) \,,
\eeq
which is associated with the vector bundle
\beq
\baf C = \lin_\f B(\f F, \, T^*\f B \ten \f F) \,.
\eeq

 The fibred charts induced on
$\f C$
are of the type
$(x^\lam, w\col \lam i j)$
and the coordinate expression of
$\eps$
is
\beq
\eps = d^\lam \ten (\der_\lam + w\col \lam i j \, y^j \, \der_i) \,.\END
\eeq
\eEx

\bEx
\label{Example: system of affine connections}
\index{smooth systems!system of affine connections}
 If
$p : \f F \to \f B$
is an affine bundle, associated with the vector bundle
$\ba p : \baf F \to \f B \,,$
then the affine connections constitute an injective smooth system
$(\f C, \zet, \eps) \,,$
where
$\zet : \f C \to \f B$
is an affine subbundle
\beq
\f C \sub \aff_\f B(\f F, \, T^*\f B \ten T\f F) \,,
\eeq
which is associated with the vector bundle
\beq
\baf C = \aff_\f B(\f F, \, T^*\f B \ten \baf F) \,.
\eeq

 The fibred charts induced on
$\f C$
are of the type
$(x^\lam, w\col \lam i j, w\co \lam i)$
and the coordinate expression of
$\eps$
is
\beq
\eps = 
d^\lam \ten 
\w(\der_\lam + (w\col \lam i j \, y^j + w\co\lam i) \, \der_i\w) \,.\END
\eeq
\eEx

\myskip

\bEx
\label{Example: system of polynomial connections}
\index{smooth systems!system of polynomial connections}
 If
$p : \f F \to \f B$
is an affine bundle, then analogously to the above 
Example \ref{Example: system of affine connections}, 
we can define 

 - the injective smooth system of polynomial connections of degree
$r \,,$
with
$1 \leq r \,,$

 - the injective smooth system of polynomial connections of any degree
$r \,,$
with 
$1 \leq r \leq k \,,$
where
$k$
is a given positive integer.\END
\eEx

\myskip

 All examples above deal with finite dimensional smooth systems of smooth connections, as it is implicitly requested in
Definition \ref{Definition: systems of connections}.

 On the other hand, we can easily extend the concept of smooth system of smooth connections, by considering an infinite dimensional system, which is the direct limit of finite dimensional smooth systems, according to the following
Example \ref{Example: infinite dimensional system of polynomial connections of any degree}.

\myskip

\bEx
\label{Example: infinite dimensional system of polynomial connections of any degree}
\index{smooth systems!system of polynomial connections}
 If
$p : \f F \to \f B$
is an affine bundle, then we obtain the smooth system of all polynomial connections by considering the family of all polynomial connections
$c : \f F \to T^*\f B \ten T\f F$
of any degree 
$r \,,$
with 
$1 \leq r < \infty \,.$

 However, we stress that such a system has a natural infinite dimensional smooth structure.\END
\eEx
 
\myskip

\bEx
\label{Example: system of principal connections}
\index{smooth systems!system of principal connections}
 If
$p : \f F \to \f B$
is a smooth left principal bundle, with structure group
$G \,,$
then the smooth principal connections constitute an injective smooth system
$(\f C, \zet, \eps) \,,$
where
$\zet : \f C \to \f B$
is the quotient bundle with respect to the group of smooth fibred left actions over
$\f B$
\beq
\id \car TL_h : T^*\f B \ten T\f F \to T^*\f B \ten T\f F \,,
\ssep{where}
h \in G \,.\END
\eeq
\eEx

\bDr
\label{Exercise: system of 1--forms}
 Let us consider a smooth manifold
$\f M$
and the trivial smooth principal bundle
$\f F \byd \f M \car \Rn \to \f B \byd \f M$
whose structure group is the abelian group
$\Rn \,.$

 Show that the system of smooth principal connections of this bundle can be naturally identified with the family of 1--forms
$\alp : \f M \to T^*\f M \,.$\END
\eDr
%--------------------------------------------------------------------%
\subsection{Smooth universal connection}
\label{Smooth universal connection}
%--------------------------------------------------------------------%
\bsm
 Eventually, given a smooth fibred manifold
$p : \f F \to \f B \,,$
we discuss the notions of \emp{reducible smooth connection} and \emp{universal smooth connection}.
 Moreover, we discuss the natural bijection between smooth systems of smooth connections and reducible smooth connections.
 
  Namely, a smooth system of smooth connections
$(\f C, \zet, \eps)$
yields a distinguished smooth connection
\beq
\eps\Upa : \f F\Upa \ucar{\f C} T\f C \to T\f F\Upa \,,
\eeq
called \emp{universal}, on the pullback smooth fibred manifold
\beq
p\Upa : \f F\Upa \byd \f C \ucar{\f B} \f F \to \f C \,.
\eeq

 Indeed, this universal connection fulfills a universal property; in fact, all connections of the systems can be obtained from the universal connection by pullback.

\myskip

 This notion was originally introduced by P.L. Garcia \cite{Gar72} in the context of principal connections of a principal bundle.
 Later, this theory has been generalised to any fibred manifold, detached from any structure group of symmetries 
(see, for instance,
\cite{Cab87,CabKol95,ManMod83d}).
 Here, we follow this generalised approach.
\esm

 Let us consider two smooth fibred manifolds
\beq
p : \f F \to \f B
\ssep{and}
\zet : \f C \to \f B \,.
\eeq

 Then, we focus our attention on the fibred manifold
(see
Definition \ref{Definition: lift fibred manifold})
\beq
p\Upa : \f F\Upa \byd \f C \ucar{\f B} \f F \to \f C \,.
\eeq

 The fibred charts
$(x^\lam, y^i)$
and
$(x^\lam, w^A)$
of
$\f F$
and
$\f C \,,$
respectively, yield the fibred chart
$(x^\lam, w^A; y^i)$
of
$\f F\Upa \,.$

 We recall the equality
\beq
T\f F\Upa = T\f C \ucar{T\f B} T\f F \,.
\eeq

 We recall that the smooth connections of the smooth fibred manifolds
\beq
p : \f F \to \f B
\ssep{and}
p\Upa : \f F\Upa \to \f C
\eeq
can be regarded equivalently

 1) as smooth fibred morphisms
\beq
c : 
\f F \to T^*\f B \ten T\f F
\ssep{and}
c\Upa : 
\f F\Upa \to T^*\f C \ten T\f F\Upa
\eeq
whose coordinate expressions are
\bal
c 
&= 
d^\lam \ten (\der_\lam + c\co \lam i \, \der_i) \,,
\\
c\Upa 
&= 
d^\lam \ten (\der_\lam + c\Upa\co \lam i \, \der_i) +
d^A \ten (\der_A + c\Upa\co Ai \, \der_i) \,,
\end{align*}
 
 2) as smooth fibred morphisms
\beq
c : 
\f F \ucar{\f B} T\f B \to T\f F
\ssep{and}
c\Upa : 
\f F\Upa \ucar{\f C} T\f C \to T\f F\Upa \,,
\eeq
whose coordinate expressions are
\bal
(x^\lam, \, y^i; \; \dt x^\lam, \, \dt y^i) \com c 
&=
(x^\lam, \, y^i; \; \dt x^\lam, \, c\co\lam i \, \dt x^\lam) \,,
\\
(x^\lam, \, w^A, \, y^i; \; \dt x^\lam, \, \dt w^A, \, \dt y^i) \com c\Upa 
&=
(x^\lam, \, w^A, \, y^i; \; \dt x^\lam, \, \dt w^A, \, 
c\Upa\co A i \, \dt w^A + c\Upa\co\lam i \, \dt x^\lam) \,,
\end{align*}
where
\beq
c\co \lam i \in \map(\f F, \Rn)
\ssep{and}
c\Upa\co A i, \; c\Upa\co \lam i \in \map(\f F\Upa, \Rn) \,.
\eeq

\myskip

\bDf
\label{Definition: reducible connections}
\index{smooth systems!reducible connection}
 A smooth connection
\beq
c\Upa : \f F\Upa \ucar{\f C} T\f C \to T\f F\Upa
\eeq
of the smooth fibred manifold
$p\Upa : \f F\Upa \to \f C$
is said to be an \emp{upper connection} of the smooth system of smooth connections.

 Moreover, such a smooth connection
$c\Upa$
is said to be \emp{reducible} if it factorises through a smooth system
$(\f C, \zet, \eps)$ 
of smooth connections of the smooth fibred manifold
$p : \f F \to \f B$
according to the following commutative diagram
(see
Definition \ref{Definition: systems of connections})
\newdiagramgrid{2x2}
{1.5, 1.5, .7}
{.5, .7}
\bdg[grid=2x2]
\f F\Upa \ucar{\f C} T\f C &\rTo^{c\Upa} &T\f F\Upa
\\
\dTo^{} &&\dTo
\\
\f C \ucar{\f B} (\f F \ucar{\f B} T\f B) 
&\rTo^{\eps} 
&T\f F &.\END
\edg
\eDf

\myskip

 The above intrinsic condition can be translated in coordinates as follows.
 
\myskip

\bPr
\label{Proposition: reducible connection}
\index{smooth systems!reducible connection}
 A smooth connection
$c\Upa : \f F\Upa \ucar{\f C} T\f C \to T\f F\Upa$
is reducible if and only if, in coordinates,
\beq
c\Upa\co Ai = 0 \,.
\eeq

 Thus, smooth a connection
$c\Upa : \f F\Upa \ucar{\f C} T\f C \to T\f F\Upa$
is reducible if and only if its coordinate expression is of the type
\beq
c\Upa =
d^\lam \ten (\der_\lam + c\Upa\co\lam i \, \der_i) + 
d^A \ten \der_A \,.
\eeq
\ePr

\bpf
 The coordinate expression of a connection
$c\Upa : \f F\Upa \ucar{\f C} T\f C \to T\f F\Upa$
is of the type
\beq
c\Upa =
d^\lam \ten (\der_\lam + c\Upa\co\lam i \, \der_i) +
d^A \ten (\der_A + c\Upa\co Ai \, \der_i) \,,
\ssep{where}
c\Upa\co\lam i, \; c\Upa\co Ai \in\map(\f F\Upa, \Rn) \,.
\eeq

 Hence, the coordinate expression of
$\pro_2 \com c\Upa : \f F\Upa \ucar{\f C} T\f C \to T\f F$
is
\beq
c\Upa =
d^\lam \ten (\der_\lam + c\Upa\co\lam i \, \der_i) +
d^A \ten (c\Upa\co Ai \, \der_i) \,.
\eeq

 Therefore, the composition of smooth maps
\newdiagramgrid{1}
{1.5, 1, 1, 1}
{1}
\bdg[grid=1]
\f F\Upa \ucar{\f C} T\f C 
&\rTo^{c\Upa} &T\f F\Upa
&\rTo^{\pro_2} &T\f F
\edg
factorises through a smooth fibred morphism
\beq
\eps : \f C \ucar{\f B} (\f F \ucar{\f B} T\f B) \to T\f F
\eeq
over
$\f F$
if and only if
\beq
c\Upa\co Ai = 0 \,.
\eeq

 Indeed, the smooth fibred morphism
\beq
\eps : \f C \ucar{\f B} (\f F \ucar{\f B} T\f B) \to T\f F \,,
\eeq
i.e. equivalently, the smooth fibred morphism
\beq
\eps : \f C \ucar{\f B} \f F \to T^*\f B \ten T\f F \,,
\eeq
turns out to be a smooth system of smooth connections.\QED
\epf

\bRm
\label{Remark: reducible connection}
\index{smooth systems!reducible connection}
 Let us consider a generic smooth connection of a generic smooth fibred manifold; if some symbols of the connection vanish in a chart, they need not to vanish in another chart.
 
 On the other hand, for each reducible connection
$c\Upa : \f F\Upa \ucar{\f C} T\f C \to T\f F\Upa$
of the fibred manifold
$p\Upa : \f F\Upa \to \f C \,,$
we have shown the following equality, in any fibred chart,
\beq
c\Upa\co Ai = 0 \,.
\eeq

 Indeed, this unusual vanishing property of some symbols of the connection in any fibred chart is possible because
$\f F\Upa \byd \f C \ucar{\f B} \f F$
is a fibred product of manifolds.

 Accordingly, if
$(x^\lam, w^A, y^i)$ 
and
$(\ac x^\mu, \ac w^B, \ac y^j)$
are two fibred charts of
$\f F\Upa \,,$
then we have
\beq
\der_A \ac y^j = 0
\ssep{and}
\der_i \ac w^B = 0 \,.\END
\eeq
\eRm

\myskip

 We can exhibit a natural bijection between smooth systems of smooth connections
$(\f C, \zet, \eps)$
of the smooth fibred manifold
$p : \f F \to \f B$
and reducible smooth connections of the smooth fibred manifold
$p\Upa : \f F\Upa \to \f C \,,$
according to the following
Proposition \ref{Proposition: bijection between systems of connections and reducible connections}.
 
 Even more, the reducible smooth connections fulfill a ``universal property" with respect to the smooth connections of the associated smooth system of smooth connections, according to the following
Theorem \ref{Theorem: universal connection}.

\myskip

\bPr
\label{Proposition: bijection between systems of connections and reducible connections}
\index{smooth systems!reducible connection}
 We have a natural bijection between smooth systems of smooth connections of the smooth fibred manifold
$p : \f F \to \f B$
and reducible smooth connections of the smooth fibred manifold
$p\Upa : \f F\Upa \to \f C$
in the following way.

\myskip

 1) If
$\eps : \f C \ucar{\f B} (\f F \ucar{\f B} T\f B) \to T\f F$
is a smooth system of smooth connections of the smooth fibred manifold
$p : \f F \to \f B \,,$ 
then the smooth map
\beq
\eps\Upa :
\f F\Upa \ucar{\f C} T\f C \to T\f F\Upa = T\f C \ucar{} T\f F :
(f\Upa \,,\, X) \mto
\W(
X \,,\, 
\w(\eps (f\Upa)\w) 
\w((T\zet) (X)\w)
\W) \,,
\eeq
with coordinate expression
\beq
\eps\Upa =
d^\lam \ten (\der_\lam + \eps\co\lam i \, \der_i) + 
d^A \ten \der_A \,,
\eeq
can be regarded as a reducible smooth connection of the smooth fibred manifold\lnb
$p\Upa : \f F\Upa \to \f C \,.$

\myskip

 2) If
$\eps\Upa : \f F\Upa \ucar{\f C} T\f C \to T\f F\Upa$
is a reducible smooth connection of the smooth fibred manifold
$p\Upa : \f F\Upa \to \f C \,,$
then the factor map
(see
Definition \ref{Definition: reducible connections})
\beq
\eps : \f F\Upa \ucar{\f B} T\f B \to T\f F \,,
\eeq
with coordinate expression
\beq
\eps =
d^\lam \ten (\der_\lam + \eps\col\lam i \, \der_i) \,,
\eeq
turns out to be a smooth system of connections of the smooth fibred manifold
$p : \f F \to \f B \,.$

\myskip

 3) The above coordinate expressions exhibit a natural bijection
\beq
\eps \mto \eps\Upa
\eeq
between the smooth systems of smooth connections of the smooth fibred manifold\lnb
$p : \f F \to \f B$
and the reducible smooth connections of the smooth fibred manifold
$p\Upa : \f F\Upa \to \f C \,.$\END
\ePr

\bTh
\label{Theorem: universal connection}
\index{smooth systems!universal connection}
 Let
$\eps : \f C \ucar{\f B} (\f F \ucar{\f B} T\f B) \to T\f F$
be a smooth system of smooth connections of the smooth fibred manifold
$p : \f F \to \f B \,.$

 Then, the following facts hold.
 
\myskip

 1) The smooth connection
\beq
\eps\Upa : \f F\Upa \ucar{\f C} T\f C \to T\f F\Upa
\eeq
of the smooth fibred manifold
$p\Upa : \f F\Upa \to \f C$
fulfills the following \emp{``universal property"}:

 - all smooth connections
$\br\gam$
of the system can be obtained as pullback of
$\eps\Upa \,,$
through the equality
\beq
\br\gam = \gam^*(\eps\Upa) \,,
\ssep{for each}
\gam \in \sec(\f B, \f C) \,,
\eeq
where
\beq
\gam^*(\eps\Upa) \in \fib(\f F \ucar{\f B} T\f B, \, T\f F)
\eeq
is the smooth connection of the smooth fibred manifold
$p : \f F \to \f B$
defined by the following commutative diagram
\newdiagramgrid{2x3x2}
{1.5, 1.5, 1.5, 1.5, .5}
{.7, .7, .7, .7}
\bdg[grid=2x3x2]
\f F\Upa \ucar{\f C} T\f C
&\rTo^{\eps\Upa}
&&&T\f F\Upa
\\
\dTo &&&&\dTo
\\
\f F \ucar{\f B} T\f B
&\rTo^{\gam_\f C \car \id} &\f F\Upa \ucar{\f B} T\f B &\rTo^{\eps} &T\f F
\\
\uTo^{\id} &&&&\uTo_{\id}
\\
\f F \ucar{\f B} T\f B
&&\rTo^{\gam^*(\eps)\Upa}
&&T\f F &.
\edg

\myskip

 2) The smooth curvature tensor
(see, for instance,
\cite{KolMicSlo93,Mod91})
\beq
R[\eps\Upa] : \f F\Upa \ucar{\f C} \Lam^2T\f C \to V_\f C\f F\Upa
\eeq
fulfills the following \emp{``universal property"}:

 - the curvature tensors
$R[\br\gam]$
of all connections
$\br\gam$
of the smooth system can be obtained as pullback of
$R[\eps\Upa]$
through the equality
\beq
R[\br\gam] = \gam^*\w(R[\eps\Upa]\w) \,,
\ssep{for each}
\gam \in \sec(\f B, \f C) \,,
\eeq
where
\beq
\gam^*\w(R[\eps\Upa]\w) \in 
\fib(\f F \ucar{\f B} \Lam^2T\f B, \, V\f F)
\eeq
is defined by the following commutative diagram
\newdiagramgrid{2x2}
{1.5, 2, .5}
{0.8, 0.8}
\bdg[grid=2x2]
\f F\Upa \ucar{\f C} \Lam^2T\f C &\rTo^{R[\eps\Upa]} &V_\f C\f F\Upa
\\
\uTo^{\gam_\f C \car T\gam} &&\dTo_{}
\\
\f F \ucar{\f B} \Lam^2T\f B &\rTo^{\gam^*\w(R[\eps\Upa]\w)} &V\f F &.
\edg
\eTh

\bpf
 The coordinate expressions of
$\eps\Upa$
and
$\br\gam$
are
\bal
\eps\Upa
&=
d^\lam \ten (\der_\lam + \eps\co\lam i \, \der_i) + d^A \ten \der_A \,,
\\
\br\gam
&=
d^\lam \ten 
\w(\der_\lam + (\eps\co\lam i \com \gam\co\lam i) \, \der_i\w) \,.
\end{align*}

 Hence, the universal property of
$\eps\Upa$
follows from the equality
\bal
\gam^*\eps\Upa 
&=
\gam^*\w(d^\lam \ten (\der_\lam + \eps\co\lam i \, \der_i) + 
d^A \ten \der_A)
\\
&=
d^\lam \ten 
\w(\der_\lam + (\eps\co\lam i \com \gam) \, \der_i\w)
\\
&=
\br\gam \,.
\end{align*}

 The coordinate expressions of
$R[\eps\Upa]$
and
$R[\br\gam]$
are
\bal
R[\eps\Upa] 
&=
- 2 \,
\W(
(\der_\lam \eps\co \mu i + \eps\co \lam j \, \der_j \eps\co \mu i\w) \, 
d^\lam \wed d^\mu +
\der_A \eps\co \mu i \, d^A \wed d^\mu
\W) \ten \der_i \,,
\\
R[\br\gam]
&=
- 2 \,
\w(
\der_\lam (\eps\co \mu i \com \gam) + 
(\eps\co \lam j \com \gam) \, 
\der_j (\eps\co \mu i \com \gam)\w) \, 
d^\lam \wed d^\mu \ten \der_i \,.
\end{align*}

 Hence, by taking into account the equalities
\beq
\der_\lam (\eps\co \mu i \com \gam) =
(\der_\lam \eps\co \mu i ) \com \gam +
(\der_A \eps\co \mu i) \com \gam \, \der_\lam \gam^A
\ssep{and}
\der_j (\eps\co \mu i \com \gam) =
(\der_j \eps\co \mu i ) \com \gam \,,
\eeq
the universal property of
$R[\eps\Upa]$
follows from the equality
\bal
\gam^*R[\eps\Upa] 
&=
- 2 \,
\gam^*\w(
(\der_\lam \eps\co \mu i + \eps\co \lam j \, \der_j \eps\co \mu i) \, 
d^\lam \wed d^\mu +
\der_A \eps\co \mu i \, d^A \wed d^\mu
\w) \ten \der_i
\\[2mm]
&=
- 2 \,
\w(
(\der_\lam \eps\co \mu i) \com \gam + 
(\eps\co \lam j \com \gam) \, (\der_j \eps\co \mu i) \com \gam +
(\der_A \eps\co \mu i ) \com \gam \, \der_\lam \gam^A
\w) \, d^\lam \wed d^\mu \ten \der_i
\\
&=
- 2 \,
\W(
\der_\lam (\eps\co \mu i \com \gam) -
(\der_A \eps\co \mu i) \com \gam \, \der_\lam \gam^A +
(\eps\co \lam j \com \gam) \, \der_j (\eps\co \mu i \com \gam) 
\\
& \quad +
(\der_A \eps\co \mu i) \com \gam \, \der_\lam \gam^A \, 
\W) \, d^\lam \wed d^\mu \ten \der_i
\\
&=
- 2 \,
\w(
\der_\lam (\eps\co \mu i \com \gam) + 
(\eps\co \lam j \com \gam) \, \der_j (\eps\co \mu i \com \gam)
\w) \, d^\lam \wed d^\mu \ten \der_i
\\[2mm]
&=
R[\br\gam] \,.\QED
\end{align*}
\epf

\myskip

 Let us examine a few distinguished examples of universal connections.
 
\myskip

\bEx
\label{Example: universal connection of a system of linear connections}
\index{smooth systems!linear connection}
 Let us refer to the smooth system of linear connections of the vector bundle
$p : \f F \to \f B$
(see
Example \ref{Example: system of linear connections}).

 Then, the associated universal connection of the system has coordinate expression
\beq
\eps\Upa =
d^\lam \ten (\der_\lam + \eps\col\lam ij \, y^j \, \der_i) + 
d\col \lam ij \ten \der^\lam{}_i{}^j \,.\END
\eeq
\eEx

\bEx
\label{Example: universal connection of a system of affine connections}
\index{smooth systems!affine connections}
 Let us refer to the smooth system of affine connections of the affine bundle
$p : \f F \to \f B$
(see
Example \ref{Example: system of affine connections}).

 Then, the associated universal connection of the system has coordinate expression
\beq
\eps\Upa =
d^\lam \ten \w(\der_\lam +
(\eps\col\lam ij \, y^j + \eps\co\lam i) \, \der_i\w) + 
d\col \lam ij \ten \der^\lam{}_i{}^j +
d\co \lam i \ten \der^\lam{}_i \,.\END
\eeq
\eEx

\myskip

 Eventually, we show that the natural Liouville form and symplectic form of a smooth manifold fulfill a well known property, that can be reinterpreted in terms of universal connection and curvature tensor of a smooth system of smooth connections.
 
\myskip

\bDr
\label{Exercise: principal connections}
\index{smooth systems!principal connection}
 Let us refer to the smooth system of smooth principal connections of the trivial principal bundle
$\f M \car \Rn \to \f M$
(see
Exercise \ref{Exercise: system of 1--forms}).

 Then, show the following facts:
 
 - the universal connection of the system can be naturally identified with the Liouville form
$\lam : T^*\f M \to T^*T^*\f M \,,$
with coordinate expression
\beq
\lam = \dt x_\mu \, d^\mu \,.
\eeq

 - the universal curvature of the system can be naturally identified with the symplectic form
$\ome : T^*\f M \to \Lam^2T^*T^*\f M \,,$
with coordinate expression
\beq
\ome \byd - d\lam = d_\mu \wed \dt d^\mu \,.
\eeq

 - the well known universal properties of the 1--form
$\lam$
and of the 2--form
$\ome \byd - d\lam$
(see
\cite{God69})
fit the universal properties of the universal connection and of its curvature.\END
\eDr
%--------------------------------------------------------------------%
%\newpage
\markboth{\rm Chapter \thechapter. Systems of connections}{\rm \thesection. F--smooth systems of connections}
\section{F--smooth systems of connections}
\label{Section: F--smooth systems of connections}
\markboth{\rm Chapter \thechapter. Systems of connections}{\rm \thesection. F--smooth systems of connections}
%--------------------------------------------------------------------%
\bsm
 We discuss the \emp{F--smooth systems of ``fibrewisely smooth connections"} of a smooth manifold
(see
Section \ref{F--smooth systems of smooth sections}).
 
 The concept of \emp{universal connection} that we have discussed for smooth systems of smooth connections can be easily extended to 
F--smooth systems of F--smooth connections.
 The reader who is interested in this subject can refer to
\cite{CabKol98}.
\esm

 Let us consider a smooth fibred manifold
$p : \f F \to \f B$
and denote its fibred charts by
$(x^\lam, y^i) \,.$

 Then, let us consider the smooth double fibred manifold
\newdiagramgrid{1}
{1.3, .8, .8, .8, .4}
{0}
\bdg[grid=1]
T^*\f B \ten T\f F 
&\rTo^{q} &\f F
&\rTo^{p} & \f B &.
\edg

\bDf
\label{Definition: fibrewise smooth tubelike connections}
 We denote by
(see
Definition \ref{Definition: tubelike sections})
\beq
\cns\utub(\f F, T^*\f B \ten T\f F) \sub 
\wob c : \f F \to T^*\f B \ten T\f F \wcb
\eeq
the subsheaf consisting of tubelike sections
$c : \f F \to T^*\f B \ten T\f F \,,$
which fulfill the following condition, \emp{without any further local smoothness requirement},

 - $c_b : \f F_b \to (T^*\f B \ten T\f F)_b$
is global and smooth, for each 
$b \in \f B \,;$

- $c_b$ 
projects over 
$\f 1_b $ 
according to the following commutative diagram
\newdiagramgrid{3x3}
{1.2, 1.5, 1.4, .5}
{.7, .7}
\bdg[grid=3x3]
\f F_b &\rTo^{c_b} &(T^*\f B \ten T\f F)_b&
\\
\dTo &&\dTo
\\
\{b\} &\rTo^{\f 1_b} &(T^*\f B \ten T\f B)_b &.\END
\edg
\eDf

 Thus, let us consider a tubelike connection
$c : \f F \to T^*\f B \ten T\f F \,.$

 We say that it is 
 
 - \emp{fibrewisely smooth} if it is smooth along the fibres
$\f F_b \sub \f F \,,$
for each
$b \in \f B \,,$
 
 - \emp{smooth} if it is smooth in its full domain
$p^{-1}(\f U) \sub \f F \,,$
for each
$\f U \in \f B \,.$

 Therefore, the sheaf of smooth tubelike connections 
\beq
c: \f F \to T^*\f B \ten T\f F
\eeq
of the smooth fibred manifold
$\f F \to \f B$
turns out to be a subsheaf of
\beq
\cns\tub(\f F, T^*\f B \ten T\f F) \sub
\cns\utub(\f F, T^*\f B \ten T\f F) \,.
\eeq

\myskip

 The following Definition is a generalisation of
Definition \ref{Definition: systems of connections},
as here we do not require that
$\f C$
be a smooth finite dimensional manifold
(hence, that the maps
$\zet$
and
$\eps$
be smooth).

\myskip

\bDf
\label{Definition: F--smooth systems of smooth connections}
\index{F--smooth systems!F--smooth system of smooth connections} 
 We define an \emp{F--smooth system of fibrewisely smooth connections} of the smooth fibred manifold
$p : \f F \to \f B$
to be a 3--plet
$(\f C, \zet, \eps) \,,$
where 

 1) $\f C$
is a set,

 2) $\zet : \f C \to \f B$ 
is a surjective map,

 3) $\eps : \f C \ucar{\f B} \f F \to T^*\f B \ten T\f F$
is a fibred map over 
$\f F \,$ and over $ \f 1_\f B: \f B \to T^*\f B \ten T\f B $,
according to the following commutative diagrams
\newdiagramgrid{2x2}
{1.5, 1.5, .5, 3, 1.5, 1.5, 1.5}
{.7, .7}
\bdg[grid=2x2]
\f C \ucar{\f B} \f F &\rTo^{\eps} &T^*\f B \ten T\f F
&&\f C \ucar{\f B} \f F &\rTo^{\eps} &T^*\f B \ten T\f F
\\
\dTo^{\pro_2} &&\dTo_{\tau_\f F}
&&\dTo &&\dTo_{\id_{T^*\f B} \ten Tp}
\\
\f F &\rTo^{\id_\f F} &\;\f F
&&\f B &\rTo^{\1_\f B} &\quad T^*\f B \ten T\f B &,
\edg
which fulfills the following condition:

 *) for each
$c \in \f C_b \,,$
with
$b \in \f B \,,$
the induced section
\beq
\eps_c : \f F_b \to (T^*\f B \ten T\f F)_b
\eeq
of the restricted smooth fibred manifold
$(T^*\f B \ten T\f F)_b \to \f F_b$
is \emp{smooth and globally defined on}
$\f F_b \,.$

 The map
$\eps : \f C \ucar{\f B} \f F \to T^*\f B \ten T\f F$
is called the \emp{evaluation map} of the system.

 We denote by
\beq
\usec(\f B, \f C) \sub \wob \gam : \f B \to \f C \wcb
\eeq
the subsheaf consisting of \emp{local} sections
$\gam : \f B \to \f C \,,$
\emp{without any smoothness requirement}.\END
\eDf

\myskip

 We leave to the reader the easy task to rephrase in the present context the notions and developments that have been previously established for F--smooth systems of smooth sections.
%--------------------------------------------------------------------%
\chapter{F--smooth connections}
\label{Chapter: F--smooth connections}
%--------------------------------------------------------------------%
\bsm
 Given an F--smooth system
$(\f S, \zet, \eps)$
of fibrewisely smooth sections of a smooth double fibred manifold
$\f G \oset{q}{\to} \f F \oset{p}{\to} \f B \,,$
we discuss the ``\emp{F--smooth connections}"
\beq
\E K : \f S \ucar{\f B} T\f B \to \E T\f S
\eeq
of the F--smooth fibred space
$\zet : \f S \to \f B$
(see, for instance, \cite{CabKol95,CabKol98,KolMod98}).

\myskip

 We mention that the curvature of an F--smooth connection
$\E K$
as above can be defined via the generalised Fr\"olicher--Nijenhuis bracket on F--smooth spaces in a way analogous to the curvature of a smooth connection on a smooth fibred manifold.
 The reader who is interested in this subject can refer, for instance, to
\cite{KolMod98,Mod91}.
\esm
%--------------------------------------------------------------------%
%\newpage
\markboth{\rm Chapter \thechapter. F--smooth connections}{\rm \thesection. F--smooth connections}
\section{F--smooth connections}
\label{F--smooth connections}
\markboth{\rm Chapter \thechapter. F--smooth connections}{\rm \thesection. F--smooth connections}
%--------------------------------------------------------------------%
\bsm
  Given an F--smooth system
$(\f S, \zet, \eps)$
of fibrewisely smooth sections of a smooth double fibred manifold
$\f G \oset{q}{\to} \f F \oset{p}{\to} \f B \,,$
we define the ``\emp{F--smooth connections}"
\beq
\E K : \f S \ucar{\f B} T\f B \to \E T\f S
\eeq
of the F--smooth fibred space
$\zet : \f S \to \f B$
and show that such a 
$\E K$
is characterised by a smooth section of a smooth bundle of the type
\beq
\Xi_{(s,u)} : (T\f F)_u \to (T\f G)_u \,,
\ssep{for each}
s \in \f S \,,
\quad
u \in T_{\zet(s)} \f B \,.
\eeq
\esm

 Thus, let us consider a smooth double fibred manifold
$\f G \oset{q}{\to} \f F \oset{p}{\to} \f B$
and denote the typical smooth fibred chart of
$\f G$
by
$(x^\mu, y^i, z^a) \,.$

 Moreover, let us consider an F--smooth system
$(\f S, \zet, \eps)$
of fibrewisely smooth sections of the above smooth double fibred manifold.

\myskip

\bDf
\label{Definition: F--smooth connections}
\index{F--smooth systems!F--smooth connection}
 We define an \emp{F--smooth connection} of the F--smooth fibred space
$\zet : \f S \to \f B$
to be an F--smooth fibred morphism over
$\f S$
and over
$T\f B \,,$
\beq
\E K : \f S \ucar{\f B} T\f B \to \E T\f S \,,
\eeq
which is linear with respect to the 2nd factor
$T\f B \,,$
according to the following commutative diagram
\newdiagramgrid{2x2}
{.7, .7, .7, .7, .5}
{.7, .7, .7, .7}
\bdg[grid=2x2]
\f S &&\rTo^{\id_\f S} &&\f S
\\
\uTo^{\pro_1} &&&&\uTo_{\tau_\f S}
\\
\f S \ucar{\f B} T\f B &&\rTo^{\E K} &&\E T\f S
\\
\dTo^{\pro_2} &&&&\dTo_{T\zet}
\\
T\f B &&\rTo^{\id_\f B} &&T\f B &,
\edg
or, equivalently, to be an F--smooth tangent valued
1--form
\beq
\E K : \f S \to T^*\f B \ten \E T\f S \,,
\eeq
which projects on
$\1_\f B : \f B \to T^*\f B \ten T\f B \,,$
according to the commutative diagram
\newdiagramgrid{2x2}
{.6, .6, .6, 1, 1, .7}
{.7, .7}
\bdg[grid=2x2]
\f S &&\rTo^{\E K} &&T^*\f B \ten \E T\f S
\\
\dTo^{\zet} &&&&\dTo_{\id \ten \E T\zet}
\\
\f B &&\rTo^{\1_\f B} &&T^*\f B \ten T\f B &.
\edg

 By recalling the representation of
$\E T \f S$
provided by 
Theorem \ref{Theorem: representation of tangent space of space of F--smooth system of sections},
the F--smooth connection
$\E K$
is characterised by a map of the type
\beq
\E K : \f S \ucar{\f B} T\f B \to \E T \f S :
(s, u) \mto \Xi_{(s,u)} \,,
\eeq
where
\beq
\Xi_{(s,u)} : (T\f F)_u \to (T\f G)_u
\eeq
is a smooth section
(see
Theorem \ref{Theorem: vector structure of TS for an F--smooth system of sections}).\END
\eDf

\myskip

 Let us consider an F--smooth connection
$\E K : \f S \to T^*\f B \ten \E T\f S \,.$

\myskip

\bDf
\label{Definition: F--smooth covariant differentials}
\index{F--smooth systems!F--smooth covariant differential}
 We define the \emp{F--smooth covariant differential} of an F--smooth section
$\sig \in \Fsec(\f B, \f S) \,,$ 
with respect to the F--smooth connection
$\E K \,,$
to be the F--smooth section
\beq
\nab \sig \byd \E T\sig - \E K \com \sig : 
\f B \to T^*\f B \ten \E V \f S \,,
\eeq
according to the commutative diagram
\newdiagramgrid{2x2}
{1.3, 3, 1.3}
{.7, .7}
\bdg[grid=2x2]
\f S &\rTo^{(\E T\sig, \E K \com \sig)} 
&(T^*\f B \ten \E T\f S) \ucar{\f S} (T^*\f B \ten \E T\f S)
\\
\uTo^{\sig} &&\dTo_{-_\f S}
\\
\f B &\rTo^{\nab \sig} &T^*\f B \ten \E V\f S &.
\edg
\eDf
%--------------------------------------------------------------------%
%\newpage
\markboth{\rm Chapter \thechapter. F--smooth connections}{\rm \thesection. F--smooth connections in the linear case}
\section{F--smooth connections in the linear case}
\label{F--smooth connections in the linear case}
\markboth{\rm Chapter \thechapter. F--smooth connections}{\rm \thesection. F--smooth connections in the linear case}
%--------------------------------------------------------------------%
\bsm
 Next, let us further suppose that
$q : \f G \to \f F$
be a \emp{vector bundle} and that the system
$(\f S, \zet, \eps)$
be injective.

 Then, we show a natural bijection
\beq
\E K \mto \C D_\E K
\eeq
between F--smooth connections
$\E K$
of the F--smooth fibred space
$\zet : \f S \to \f B$
and certain smooth differential operators
$\C D_\E K$
between finite dimensional smooth manifolds
(see
Definition \ref{Definition: F--smooth systems of smooth sections})
\beq
\C D_\E K : 
\tub_\f S(\f F, \f G) \to \tub(\f F, T^*\f B \ten \f G) \,.
\eeq

\myskip

 We stress that the above smooth differential operators
$\C D_\E K \,,$
play a role analogous to the matrix of symbols
$(K^i_\lam)$
of a standard smooth connection
$K$
of a standard smooth fibred manifold.
\esm

 Thus, let us consider a smooth vector bundle
$\f G \to \f F \,,$
an injective F--smooth system 
$(\f S, \zet, \eps)$
of fibrewisely smooth sections of the smooth double fibred manifold
$\f G \to \f F \to \f B \,.$

  We recall that, in the linear case, the F--smooth fibred space
$\zet : \f S \to \f B$
inherits naturally a vector structure
(see
Proposition \ref{Proposition: vector structure of S in the case when G to F is a vector bundle})
and that there is a natural F--smooth linear fibred isomorphism
$\E V\f S \to \f S \ucar{\f B}\f S$
over
$\f S$
(see
Corollary \ref{Corollary: natural splitting of S if G to F is a vector bundle and the system is injective}).
 
\myskip

\bNt
\label{Note: F--smooth covariant differential}
 We can regard the covariant differential of a section
$\sig \in \Fsec(\f B, \f S) \,,$
with respect to an the F--smooth connection
$\E K \,,$
as an F--smooth section
\beq
\nab \sig : \f B \to T^*\f B \ten \f S \,.\END
\eeq
\eNt

\myskip

 The covariant differential
$\nab_\E K$
associated with the F--smooth connection
$\E K$
is a differential operator
$\C D_\E K$
of a certain type.
 Indeed, there is a natural bijection between these objects.
 
 This result extends to the present F--smooth framework an analogous result holding for standard smooth connections of smooth fibred manifolds.

\myskip

\bPr
\label{Proposition: bijection between covariant differentials and differential operators}
 The following facts hold.
 
\myskip

 1)  Let us consider an F--smooth connection
$\E K \,.$

 Then, there exists a unique F--smooth 
sheaf morphism
\beq
\C D \eqv \C D[\E K] : 
\tub_\f S(\f F, \f G) \to \tub(\f F, T^*\f B \ten\f G) :
\phi \mto \C D \phi \,,
\eeq
such that, for each
$\phi \in \tub_\f S(\f F, \f G)$
and
$u \in T\f B \,,$
\beq
\wha{(\C D \phi) (u)} = \nab_u \, \wha\phi \,.
\eeq

 The sheaf morphism
$\C D$
turns out to be a differential  operator of horizontal order 1,
which, for each
$b \in \f B \,,$ 
factorises fibrewisely, through a smooth sheaf morphism
\beq
\chC D_b : 
\sec(\f F_b, \f G_b) \to \sec(\f F_b, T^*_b\f B \ten \f G_b) \,,
\eeq
according to the following commutative diagram
\newdiagramgrid{2x2}
{2, 2, 1.7}
{.7, .7}
\bdg[grid=2x2]
\tub_\f S(\f F, \f G) &\rTo^{\C D} 
&\tub(\f F, \; T^*\f B \ten \f G)
\\
\dTo^{} &&\dTo_{}
\\
\sec_\f S(\f F_b, \f G_b) &\rTo^{\chC D_b} 
&\sec(\f F_b, \; T^*_b\f B \ten \f G_b) &.
\edg

 The coordinate expression of the sheaf morphism
$\C D$
is of the type
\beq
(\C D \phi)^a_\lam = 
\der_\lam \phi^a - \chC D{}^a_\lam (\phi) \,,
\eeq
where
$\chC D{}^a_\lam$
are smooth sheaf morphisms
\beq
\chC D{}^a_\lam : \tub_\f S(\f F, \f G) \to \Rn \,,
\eeq
which, for each
$b \in \f B \,,$ 
factorise fibrewisely, through smooth sheaf morphisms
\beq
\chC D_b{}^a_\lam : 
\sec_\f S(\f F_b, \f G_b) \to \Rn \,,
\eeq
according to the following commutative diagram
\newdiagramgrid{2x2}
{2, 2, .5}
{.7, .7}
\bdg[grid=2x2]
\tub_\f S(\f F, \f G) &\rTo^{\chC D{}{}^a_\lam} &\Rn
\\
\dTo^{} &&\dTo_{\id_\Rn}
\\
\sec_\f S(\f F_b, \f G_b) &\rTo^{\chC D{}_b{}^a_\lam} &\Rn &.
\edg

\myskip

 2) Conversely, let us consider an F--smooth sheaf morphism
\beq
\C D : \tub_\f S(\f F, \f G) \to \tub(\f F, T^*\f B \ten\f G) \,,
\eeq
whose local coordinate expression is of the type
\beq
(\C D \phi)^a_\lam = 
\der_\lam \phi^a - \chC D{}^a_\lam (\phi) \,,
\eeq
as in the above item 1).

 Then, there exists a unique F--smooth connection
$\E K$
of
$\f S \,,$
such that, for each\lnb
$\sig \in \sec(\f B, \f S) \,,$
the associated sheaf morphism 
\beq
\C D \eqv \C D[\E K] : 
\tub_\f S(\f F, \f G) \to \tub(\f F, T^*\f B \ten\f G)
\eeq
be given by
\beq
\nab[\E K] \, \sig = \wha{\C D \, \br\sig} \,.
\eeq

 Indeed, we obtain, for each
$\sig \in \sec(\f B, \f S) \,,$
\beq
\E K \com \sig = d\sig - \wha{\C D \, \br\sig} \,.\END
\eeq
\ePr

\myskip

\bDf
\label{Definition: linear F--smooth connections}
  Let us suppose that the fibred manifold
$q : \f G \to \f F$
be a vector bundle.
 Then, the F--smooth connection
$\E K$
is said to be \emp{linear} if it is a linear fibred morphism over
$\1_\f B : \f B \to T^*\f B \ten T\f B \,,$
according to the following commutative diagram
\newdiagramgrid{2x2}
{1.3, 1.3, 1.3}
{.7, .7}
\bdg[grid=2x2]
\f S &\rTo^{\E K} &T^*\f B \ten \E T\f S
\\
\dTo^{\zet} &&\dTo_{\id_{T^*\f B} \ten \E T\tau_\f S}
\\
\f B &\rTo^{\f 1_\f B} &T^*\f B \ten T\f B &.\END
\edg
\smallskip
\eDf
%--------------------------------------------------------------------%
{\fz\fz
%--------------------------------------------------------------------%
\thelistofsymbols
\markboth{\rm List of Symbols}{\rm List of Symbols}
%--------------------------------------------------------------------%
{\bf Introduction}

\myskip

$\Map(\f M, \f N)$
\qquad
set of global smooth maps
$f : \f M \to \f N$
\qquad
\S Introduction

$\Sec(\f B, \f F)$
\qquad
set of global smooth sections
$s : \f B \to \f F$
\qquad
\S Introduction

$\sec(\f B, \f F)$
\qquad
sheaf of local smooth sections
$s : \f B \to \f F$
\qquad
\S Introduction

\myskip

\noindent {\bf Smooth manifolds and F--smooth spaces}

\myskip

$(x^i) : \f M \to \Rn^m$
\qquad
smooth chart of a smooth manifold
\qquad
Def \ref{Definition: smooth manifolds}

$(\f S, \C C)$
\qquad
F--smooth space
\qquad
Def \ref{Definition: F--smooth spaces}

$c : \f I_c \to \f S$
\qquad
basic curve
\qquad
Def \ref{Definition: F--smooth spaces}

\myskip

\noindent {\bf Systems of maps}

\myskip

$(\f S, \eps)$
\qquad
smooth system of smooth maps
\qquad
Def \ref{Definition: smooth systems of smooth maps}

$\eps : \f S \car \f M \to \f N$
\qquad
evaluation map
\qquad
Def \ref{Definition: smooth systems of smooth maps}

$\eps_\f S : \f S \to \Map(\f M, \f N) : 
s \mto \br s$
\qquad
induced map
\qquad
Def \ref{Definition: smooth systems of smooth maps}

$(\eps_\f S)^{-1} : 
\Map_\f S(\f M, \f N) \to \f S : f \mto \wha f$
\qquad
inverse map
\qquad
Def \ref{Definition: smooth systems of smooth maps}

$\br s : \f M \to \f N : m \mto \eps (s, m)$
\qquad
selected map
\qquad
Def \ref{Definition: smooth systems of smooth maps}

$\Map_\f S(\f M, \f N) \byd \eps_\f S(\f S) \sub \Map(\f M, \f N)$
\qquad
subset of selected maps
\qquad
Def \ref{Definition: smooth systems of smooth maps}

$T\eps : T\f S \car T\f M \to T\f N$
\qquad
tangent prolongation
\qquad
Pro \ref{Proposition: tangent prolongations of systems of maps}

$T_1\eps : T\f S \car \f M \to T\f N$
\qquad
tangent prolongation
\qquad
Pro \ref{Proposition: tangent prolongations of systems of maps}

$T_2\eps : \f S \car T\f M \to T\f N$
\qquad
tangent prolongation
\qquad
Pro \ref{Proposition: tangent prolongations of systems of maps}

\myskip

$(\f S, \eps)$
\qquad
F--smooth system of smooth maps
\qquad
Def \ref{Definition: F--smooth systems of smooth maps}

$\eps : \f S \car \f M \to \f N$
\qquad
evaluation map
\qquad
Def \ref{Definition: F--smooth systems of smooth maps}

$\eps_\f S : \f S \to \Map(\f M, \f N) : 
s \mto \br s$
\qquad
induced map
\qquad
Def \ref{Definition: F--smooth systems of smooth maps}

$(\eps_\f S)^{-1} : 
\Map_\f S(\f M, \f N) \to \f S : f \mto \wha f$
\qquad
inverse map
\qquad
Def \ref{Definition: F--smooth systems of smooth maps}

$\br s : \f M \to \f N : m \mto \eps (s, m)$
\qquad
selected map
\qquad
Def \ref{Definition: F--smooth systems of smooth maps}

$\Map_\f S(\f M, \f N) \byd \eps_\f S(\f S) \sub \Map(\f M, \f N)$
\qquad
subset of selected maps
\qquad
Def \ref{Definition: F--smooth systems of smooth maps}

$c^*(\eps) : \f I_c \car \f M \to \f N :
(\lam, m) \mto \eps\w(c(\lam),m\w)$
\qquad
pullback of a curve
\qquad
The \ref{Theorem: F--smooth structure of F--smooth systems of smooth maps}

$\E X_s \byd \w[(\wha c_s, \lam)\w]_\sim$
\qquad
tangent vector
\qquad
Def \ref{Definition: F--smooth tangent space of the F--smooth space of smooth maps}

$\Xi_s \byd \w(T_1(\wha c^*(\eps))\w)_{|(\lam, 1)} : \f M \to T\f N$
\qquad
representation of tangent vector
\qquad
The \ref{Theorem: representation of tangent space of space of maps}

$\ol\Xi_s : T\f M \to T\f N$
\qquad
representation of tangent vector
\qquad
Cor \ref{Corollary: equivalent representation of the tangent space of space of maps}

$\tau_\f S : \E T\f S \to \f S : \Xi_s \mto s$
\qquad
natural projection
\qquad
Lem \ref{Lemma: natural maps tau and T eps}

$\E T_1\eps : \E T\f S \car \f M \to T\f N :
(\Xi_s, m) \mto \Xi_s(m)$
\qquad
natural evaluation map
\qquad
Lem \ref{Lemma: natural maps tau and T eps}

\myskip

\noindent {\bf Systems of sections}

\myskip

$\tub(\f F, \f G) \sub \sec(\f F, \f G)$
\qquad
subsheaf of tubelike sections
\qquad
Def \ref{Definition: tubelike sections}

$(\f S, \zet, \eps)$
\qquad
smooth system of smooth sections
\qquad
Def \ref{Definition: systems of sections}

$\zet : \f S \to \f B$
\qquad
smooth projection
\qquad
Def \ref{Definition: systems of sections}

$\eps : \f S \ucar{\f B} \f F \to \f G$
\qquad
smooth evaluation map
\qquad
Def \ref{Definition: systems of sections}

$\eps_\f S : \sec(\f B, \f S) \to \tub(\f F, \f G) : 
\sig \mto \br\sig$
\qquad
sheaf morphism
\qquad
Def \ref{Definition: systems of sections}

$(\eps_\f S)^{-1} : 
\tub_\f S(\f F, \f G) \to \sec(\f B, \f S) :
\phi \mto \wha\phi$
\qquad
inverse map
\qquad
Def \ref{Definition: systems of sections}

$\br\sig : \f F \to \f G : f_b \mto \eps (\sig(b), f_b)$
\qquad
selected section
\qquad
Def \ref{Definition: systems of sections}

$\tub_\f S(\f F, \f G) \byd \eps_\f S\w(\sec(\f B, \f S)\w) \sub \tub(\f F, \f G)$
subset of selected sections
\qquad
Def \ref{Definition: systems of sections}

$\f F\Upa \byd \f S \ucar{\f B} \f F$
\qquad
lifted fibred manifold
\qquad
Def \ref{Definition: lift fibred manifold}

$p\Upa : \f F\Upa \to \f S : (s_b, f_b) \to s_b$
\qquad
natural projection
\qquad
Def \ref{Definition: lift fibred manifold}

$\eps : \f F\Upa \to \f G$
\qquad
evaluation map
\qquad
Def \ref{Definition: lift fibred manifold}

\myskip

$\utub(\f F, \f G) \sub \wob s : \f F \to \f G \wcb$
sheaf of fibrewisely smooth tubelike sections
\qquad
Def \ref{Definition: fibrewisely smooth tubelike sections}

$\tub(\f F, \f G) \sub \utub(\f F, \f G)$
subsheaf of smooth tubelike sections
\qquad
Def \ref{Definition: fibrewisely smooth tubelike sections}

$(\f S, \zet, \eps)$
\qquad
F--smooth system of fibrewisely smooth sections
\qquad
Def \ref{Definition: F--smooth systems of smooth sections}

$\zet : \f S \to \f B$
\qquad
smooth projection
\qquad
Def \ref{Definition: F--smooth systems of smooth sections}

$\eps : \f S \ucar{\f B} \f F \to \f G$
\qquad
smooth evaluation map
\qquad
Def \ref{Definition: F--smooth systems of smooth sections}

$\eps_\f S : \usec(\f B, \f S) \to \utub(\f F, \f G) : 
\sig \mto \br\sig$
\qquad
sheaf morphism
\qquad
Def \ref{Definition: F--smooth systems of smooth sections}

$(\eps_\f S)^{-1} : 
\utub_\f S(\f F, \f G) \to \usec(\f B, \f S) :
\phi \mto \wha\phi$
\qquad
inverse map
\qquad
Def \ref{Definition: F--smooth systems of smooth sections}

$\br\sig : \f F_b \to \f G_b : f_b \mto \eps \w(\sig(b), f_b\w)$
\qquad
selected section
\qquad
Def \ref{Definition: F--smooth systems of smooth sections}

$\tub_\f S(\f F, \f G) \byd \eps_\f S\w(\sec(\f B, \f S)\w) \sub \tub(\f F, \f G)$
subset of selected sections
\qquad
Def \ref{Definition: F--smooth systems of smooth sections}

$c^* (\f F) \byd
\{(\lam, f) \in \f I_c \car \f F \sst c(\lam) = p(f)\} \sub
\f I_c \car \f F$
\qquad
pullback space
\qquad
Lem \ref{Lemma: c*(F)}

$c^*(p) : c^*(\f F) \to \f I_c : (\lam, f) \mto \lam$
\qquad
pullback map
\qquad
Lem \ref{Lemma: c*(F)}

$c^*_\f F : c^*(\f F) \to \f F : (\lam, f) \mto f$
\qquad
pullback map
\qquad
Lem \ref{Lemma: c*(F)}

$\wha c^*(\eps) : c^*(\f F) \to \f G :
(\lam, f) \mto \eps\w(\wha c(\lam), f\w)$
\qquad
pullback map
\qquad
Lem \ref{Lemma: c* : c*(F) to G}

$c^*(\br\sig) : c^*(\f F) \to \f G$
\qquad
pullback section
\qquad
Lem \ref{Lemma: pullback section induced by a curve}

$T \w(\wha c^*(\eps)\w) : T \w(c^*(\f F)\w) \to T\f G$
\qquad
tangent map
\qquad
Lem \ref{Lemma: T c * eps}

$\E X_s \byd \w[(\wha c, \lam)\w]_\sim$
\qquad
tangent vector
\qquad
Def \ref{Definition: tangent space of the F--smooth space of smooth tubelike sections}

$\Xi_u : (T\f F)_u \to (T\f G)_u$
\qquad
representative of a tangent vector
\qquad
The \ref{Theorem: representation of tangent space of space of 
F--smooth system of sections}

$\tau_\f S : \E T\f S \to \f S : \E X_s \mto s$
\qquad
projection
\qquad
Pro \ref{Proposition: tau S, T zet, T eps for an F--smooth system of sections}

$\E T\zet : \E T\f S \to T\f B$
\qquad
projection
\qquad
Pro \ref{Proposition: tau S, T zet, T eps for an F--smooth system of sections}

$\E V \f S \byd (T\zet)^{-1} (0) \sub \E T\f S$
\qquad
vertical subspace
\qquad
Pro \ref{Proposition: vertical space of an F--smooth system of smooth sections}

$\E T \sig : T\f B \to \E T \f S$
\qquad
tangent prolongation of a section
\qquad
Def \ref{Definition: tangent prolongation of sections}

$\C D : \tub(\f F, \f G) \to \tub(\f F, \acf G)$
\qquad
differential operator
\qquad
Def \ref{Definition: compatible smooth tubelike operators}

$\C D : \utub(\f F, \f G) \to \utub(\f F, \acf G)$
\qquad
compatible differential operator
\qquad
Def \ref{Definition: compatible smooth tubelike operators}

$\whaC D : \Fsec(\f B, \f S) \to \Fsec(\f B, \acf S) : 
\wha \sig \mto \wha{\C D (\sig)}$
\qquad
F--smooth differential operator
\qquad
Pro \ref{Proposition: F--smooth operator}

\myskip

\noindent {\bf Systems of connections}

\myskip

$\cns\tub(\f F, T^*\f B \ten T\f F) \sub
\tub(\f F, T^*\f B \ten T\f F)$
\qquad
subsheaf 0f smooth tubelike connections
\qquad
\$\ref{Smooth systems of smooth connections}

$(\f C, \zet, \eps)$
\qquad
smooth system of smooth connections
\qquad
Der \ref{Definition: systems of connections}

$\zet : \f C \to \f B$
\qquad
projection
\qquad
Def \ref{Definition: systems of connections}

$\eps : \f C \ucar{\f B} \f F \to T^*\f B \ten T\f F$
\qquad
evaluation map
\qquad
Def \ref{Definition: systems of connections}

$\eps_\f C : \sec(\f B, \f C) \to \cns\tub(\f F, T^*\f B \ten T\f F) : 
\gam \mto \br\gam$
\qquad
sheaf morphism
\qquad
Def \ref{Definition: systems of connections}

$(\eps_\f C)^{-1} : 
\cns\tub_\f C(\f F, T^*\f B \ten T\f F) \to \sec(\f B, \f C) :
c \mto \wha c$
\qquad
inverse sheaf morphism
\qquad
Def \ref{Definition: systems of connections}

$c\Upa : \f F\Upa \ucar{\f C} T\f C \to T\f F\Upa$
\qquad
upper connection
\qquad
Def \ref{Definition: reducible connections}

\myskip

$\cns\utub(\f F, T^*\f B \ten T\f F) \sub 
\wob c : \f F \to T^*\f B \ten T\f F \wcb$
\qquad
s. of fib.ly smooth connections
\qquad
Def \ref{Definition: fibrewise smooth tubelike connections}

$\cns\tub(\f F, T^*\f B \ten T\f F) \sub
\cns\utub(\f F, T^*\f B \ten T\f F)$
subsheaf of smooth connections
\qquad
Def \ref{Definition: fibrewise smooth tubelike connections}

$(\f C, \zet, \eps)$
\qquad
system of fibrewisely smooth connections
\qquad
Def \ref{Definition: F--smooth systems of smooth connections}

$\zet : \f C \to \f B$
\qquad
surjective map
\qquad
Def \ref{Definition: F--smooth systems of smooth connections}

$\eps : \f C \ucar{\f B} \f F \to T^*\f B \ten T\f F$
\qquad
evaluation map
\qquad
Def \ref{Definition: F--smooth systems of smooth connections}

$\usec(\f B, \f C) \sub \wob \gam : \f B \to \f C \wcb$
\qquad
subsheaf of sections
\qquad
Def \ref{Definition: F--smooth systems of smooth connections}

\myskip

\noindent {\bf F--smooth connections}

\myskip

$\E K : \f S \ucar{\f B} T\f B \to \E T\f S$
\qquad
F--smooth connections
\qquad
Def \ref{Definition: F--smooth connections}

$\C D \eqv \C D[\E K] : 
\tub_\f S(\f F, \f G) \to \tub(\f F, T^*\f B \ten\f G) :
\phi \mto \C D \phi$
\qquad
differential operator
\qquad
Pro \ref{Proposition: bijection between covariant differentials and differential operators}
%--------------------------------------------------------------------%
\newpage
\addcontentsline{toc}{chapter}{\hspace*{1.25em}\indexname}
\printindex
\markboth{\rm Index}{\rm Index}

%--------------------------------------------------------------------%

}

%--------------------------------------------------------------------%

\begin{thebibliography}{111}
\markboth{\rm Bibliorgaphy}{\rm Bibliorgaphy}
\addcontentsline{toc}{chapter}{\hspace*{1.25em}\bibliographyname}

\fz
%--------------------------------------------------------------------%
\bibitem{Bom67}
{\sc J. Boman}: 
\emp{Differentiability of a function and of its composition with
functions of one variable}, 
Math. Scand. {\bf 20} (1967), 249--268.
%--------------------------------------------------------------------%
\bibitem{Cab87}
{\sc A. Cabras}:
\emp{On the universal connection of a system of connections},
Note Mat. {\bf 7} (1987), 173-209.
%--------------------------------------------------------------------%
\bibitem{CabKol95}
{\sc A. Cabras, I. Kol\'a\v r}:
\emp{Connections on some functional bundles},
Czechoslovak Math. J. {\bf 45}(120) (1995), 529--548.
%--------------------------------------------------------------------%
\bibitem{CabKol98}
{\sc A. Cabras, I. Kol\'a\v r}:
\emp{The universal connection of an arbitrary system},
Ann. Mat. Pura Appl. {\bf 4}(174) (1998), 1--11.
%--------------------------------------------------------------------%
\bibitem{DodMod86}
{\sc C. T. J. Dodson, M. Modugno}: 
{\em Connections over connections and a universal calculus\/}, 
in ``Relativit\`a Generale e Fisica della Gravitazione", 
Eds.: R. Fabbri, M. Modugno, 
Proc. of VI Nat. Conf., Florence, 10-13 Oct. 1984,
Pitagora Editrice, Bologna, 1986, 89--97.
%--------------------------------------------------------------------%
\bibitem{Fro82}
{\sc A. Fr\"olicher}:
\emp{Smooth structures},
in Lecture Notes in Math. {\bf 962} (1982), Springer-Verlag, 69-81.
%--------------------------------------------------------------------%
\bibitem{Gar72}
{\sc P.L. Garc\'\i a}:
\emp{Connections and 1-jet fibre bundle},
Rend. Sem. Mat. Univ. Padova {\bf 47} (1972), 227--242.
%--------------------------------------------------------------------%
\bibitem{God69}
{\sc C. Godbillon}:
G\'eom\'etrie diff\'erentielle et m\'ecanique analytique,
Hermann, Paris, 1969.
%--------------------------------------------------------------------%
\bibitem{JanMod02c}
{\sc J. Jany\v ska, M. Modugno}:
{\em Covariant Schr\"odinger operator\/}, 
Jour. Phys.: A, Math. Gen, {\bf 35}, (2002), 8407--8434.
%--------------------------------------------------------------------%
\bibitem{JanMod20}
{\sc J. Jany\v ska, M. Modugno}:
{An Introduction to Covariant Quantum Mechanics}, 
Book in preparation, 2020.
%--------------------------------------------------------------------%
\bibitem{Kol87}
{\sc I.~Kol\'a\v{r}}: 
\emp{Some natural operators in differential geometry}, in Proc. Conf.
Diff. Geom. and Its Appl., Brno 1986, D. Reidel, 1987, 91--110.
%--------------------------------------------------------------------%
\bibitem{KolMicSlo93}
{\sc I. Kol\'a\v r, P. Michor, J. Slov\'ak}: {Natural operators in differential geometry},
Springer-Verlag, Berlin, 1993.
%--------------------------------------------------------------------%
\bibitem{KolMod98}
{\sc I. Kol\'{a}\v{r}, M. Modugno}: 
\emp{The Fr\"olicher-Nijenhuis bracket on some functional spaces},
Ann. Polon. Math. {\bf 68} (1998), 97--106.
%--------------------------------------------------------------------%
\bibitem{KriMic97}
{\sc A. Kriegl, P.W. Michor}: 
The Convenient Setting of Global Analysis,
in: "Mathematical Surveys and Monographs", {\bf 3}, Shiva Publ. Orpington, 1997.
%--------------------------------------------------------------------%
\bibitem{ManMod83d}
{\sc L. Mangiarotti, M. Modugno}: 
\emp{Fibred spaces, jet spaces and connections for field theories},
in "Geometry and Physics",
Edr.: M. Modugno, 
Proc. Intern. Meet., Florence 12-15 October 1982, 
Pitagora Editrice, Bologna, 1983, 135--165.
%--------------------------------------------------------------------%
\bibitem{MarMod91}
{\sc K. B. Marathe, M. Modugno}: 
{\em Polynomial connections on affine bundles\/}, 
Tensor N. S. {\bf 50}, 1 (1991), 35--49.
%--------------------------------------------------------------------%
\bibitem{Mic80}
{\sc P.W. Michor}: 
Manifolds of differentiable mappings,
Shiva Mathematics Series, {\bf 3}, 1980. 
%--------------------------------------------------------------------%
\bibitem{Mic93}
{\sc P.W. Michor}: 
\emp{The relation between systems and associated bundles},
Annali di Matematica Pura ed Applicata (IV), Vol. CLXIII (1993), 385--399.
%--------------------------------------------------------------------%
\bibitem{Mod91}
{\sc M. Modugno}: 
\emp{Torsion and Ricci tensor for non linear connections},
Diff. Geom. Appl. {\bf 1} (1991), 177--192.
%--------------------------------------------------------------------%
\bibitem{Sch67}
{\sc L. Schwartz}: 
Course d'analyse, Vol. I, Vol. II,
Hermann, Paris, 1967. 
%--------------------------------------------------------------------%
\bibitem{Slo86}
{\sc J. Slov\'ak}: 
\emp{Smooth structures on fibre jet spaces}, 
Czechoslovak Math. J. {\bf 36}(111) (1986), 358-375.
%--------------------------------------------------------------------%
\end{thebibliography}
\end{document}